\documentclass[a4paper,10pt]{article}
\usepackage{amsmath, amssymb, amsthm, mathtools}

\usepackage{bm}
\usepackage[margin=1in]{geometry}

\usepackage{dsfont}
\usepackage{bbm}
\usepackage{graphicx}
\usepackage{float}       
\usepackage[colorlinks=true, linkcolor=blue, citecolor=red, urlcolor=blue]{hyperref}
\usepackage{empheq}

\usepackage{cite}
\usepackage{bibentry}
\usepackage{cancel}

\usepackage{alphalph}

\title{\textbf{From many-body Schrödinger dynamics to magnetic Vlasov equation}}
\author{\textbf{Li Chen, Sebastian Stegmüller}\\[1ex]
\normalsize University of Mannheim, Department of Mathematics, 68131 Mannheim, Germany\\
\normalsize \texttt{chen@math.uni-mannheim.de}, \texttt{sebastian.stegmueller@students.uni-mannheim.de}}
\date{November 7, 2026}

\DeclareMathSizes{11}{11}{8}{6}

\theoremstyle{plain}
\newtheorem{theorem}{Theorem}[section]
\newtheorem{lemma}{Lemma}[section]
\newtheorem{proposition}{Proposition}[section]
\newtheorem{corollary}{Corollary}[section]

\theoremstyle{definition}
\newtheorem{definition}{Definition}[section]

\newtheorem{remark}{Remark}[section]

\newtheorem{assumptionA}{Assumption}

\newtheorem{assumptionI}{Assumption}

\usepackage{tikz}
\usetikzlibrary{positioning}

\usetikzlibrary{arrows.meta, shapes.geometric}

\usepackage{xcolor}

\usepackage{booktabs}

\newcommand{\C}{\mathbf{C}}

\newcommand{\rate}{\vartheta}
\newcommand{\mrate}{\eta}

\newcommand{\dmat}{\gdef\dmat##1(##2;##3){\big\langle \Psi_{N,t}, a_{##2}^\dagger a_{##3} \Psi_{N,t} \big\rangle}}\dmat

\newcommand{\dZk}{(dq \, dv)^{\otimes k}}

\newcommand{\dXkYk}{(dx \, dy)^{\otimes k}}
\newcommand{\dXkYkplus}{(dx \, dy)^{\otimes (k+1)}}

\newcommand{\mpi}{\Pi}

\newcommand{\Vext}{V}

   \newcommand{\mhus}{\mathfrak{w}}

\mathtoolsset{showonlyrefs=true}

\renewcommand{\theenumi}{(\roman{enumi})}
\renewcommand{\labelenumi}{\theenumi}

\newcommand{\brkminrows}{4}
\makeatletter
\newcount\brk@rows
\newif\ifbrk@ok
\let\brk@measure\measure@
\def\measure@#1{\brk@measure{#1}\global\brk@rows\row@}
\let\brk@cr\math@cr@@
\def\math@cr@@[#1]{%
  \ifx\math@cr@@@\math@cr@@@align
    \expandafter\@firstoftwo
  \else
    \expandafter\@secondoftwo
  \fi
  {\gdef\brk@skip{#1}\global\let\brk@held\@empty \global\brk@okfalse \brk@peek}%
  {\brk@cr[#1]}}
\def\brk@eatspace#1{\afterassignment#1\let\brk@tok= }
\def\brk@peek{\futurelet\brk@tok\brk@peek@}
\def\brk@peek@{%
  \let\brk@next\brk@finish
  \ifx\brk@tok\@sptoken \def\brk@next{\brk@eatspace\brk@peek}\fi
  \ifx\brk@tok\label    \let\brk@next\brk@eatlabel \fi
  \ifx\brk@tok&         \let\brk@next\brk@eatamp \fi
  \brk@next}
\def\brk@eatlabel\label#1{\g@addto@macro\brk@held{\label{#1}}\brk@peek}
\def\brk@eatamp&{\g@addto@macro\brk@held{&}\brk@rel}
\def\brk@rel{\futurelet\brk@tok\brk@rel@}
\def\brk@rel@{%
  \let\brk@next\brk@finish
  \ifx\brk@tok\@sptoken \def\brk@next{\brk@eatspace\brk@rel}\fi
  \ifx\brk@tok\quad     \let\brk@next\brk@eatquad \fi
  \ifx\brk@tok\qquad    \let\brk@next\brk@eatqquad \fi
  \ifx\brk@tok\hspace   \let\brk@next\brk@eathspace \fi
  \ifx\brk@tok=    \global\brk@oktrue \fi
  \ifx\brk@tok\le  \global\brk@oktrue \fi
  \ifx\brk@tok\leq \global\brk@oktrue \fi
  \ifx\brk@tok\ge  \global\brk@oktrue \fi
  \ifx\brk@tok\geq \global\brk@oktrue \fi
  \ifx\brk@tok<    \global\brk@oktrue \fi
  \ifx\brk@tok>    \global\brk@oktrue \fi
  \brk@next}
\def\brk@eatquad\quad{\g@addto@macro\brk@held{\quad}\brk@rel}
\def\brk@eatqquad\qquad{\g@addto@macro\brk@held{\qquad}\brk@rel}
\def\brk@eathspace\hspace{\futurelet\brk@tok\brk@hsp@}
\def\brk@hsp@{\ifx\brk@tok*\expandafter\brk@hsp@s\else\expandafter\brk@hsp@n\fi}
\def\brk@hsp@s*#1{\g@addto@macro\brk@held{\hspace*{#1}}\brk@rel}
\def\brk@hsp@n#1{\g@addto@macro\brk@held{\hspace{#1}}\brk@rel}
\def\brk@finish{%
  \ifx\brk@tok\black@ \global\@eqpen\z@ \fi  
  \ifbrk@ok
    \ifnum\brk@rows<\brkminrows\relax\else
      \ifnum\@eqpen=\interdisplaylinepenalty \global\@eqpen\z@ \fi
    \fi
  \fi
  \ifnum0=`{\fi \iffalse}\fi
  \math@cr@@@\noalign{\vskip\brk@skip\relax}\brk@held}
\makeatother

\usepackage{etoolbox}
\makeatletter
\def\tagfix@label#1{\gdef\df@label{#1}}
\def\tagfix@err#1{\PackageError{preamble}{Tag fix, step #1, failed: place the
  block after loading amsmath and mathtools}{}}
\def\tagfix@do#1{%
  \patchcmd{#1}{\let\label\@gobble}{\let\label\tagfix@label}{}{\tagfix@err1}%
  \apptocmd{#1}{\global\let\df@label\@empty}{}{\tagfix@err3}}
\@ifundefined{brk@measure}{\tagfix@do\measure@}{\tagfix@do\brk@measure}
\patchcmd{\math@cr@@@align@measure}{\tagform@\theequation}%
  {\tagform@\theequation\global\let\df@label\@empty}{}{\tagfix@err2}
\makeatother

\begin{document}
\newcommand{\One}{\mathbbm{1}}

\maketitle

\begin{abstract}
We analyse the combined mean-field and semiclassical limit from the time-dependent Schrödinger dynamics of $N$ spinless fermions in a non-constant magnetic field to the magnetic Vlasov equation with scaling $\hbar = N^{-1/3}$. With the velocity-shifted coherent states and the induced magnetic Husimi measure, we derive an exact BBGKY hierarchy with a remainder of order $\hbar^{1/2-}$.  The main result is the weak $L^1$ convergence of the $k$-particle magnetic Husimi measures to the $k$-copies of the solution to the magnetic Vlasov equation. Furthermore, the convergence holds in the Wasserstein distance $\mathbb{W}_2$ when the second moment of initial Husimi measure correspondingly converges.
\end{abstract}

\vspace{0.5cm}
\noindent \textbf{Keywords:} Magnetic Vlasov equation $\cdot$ Many-body Schrödinger equation $\cdot$ Semiclassical limit $\cdot$ Mean-field limit $\cdot$ Magnetic Husimi measure $\cdot$ Coherent states $\cdot$ BBGKY hierarchy $\cdot$ Interacting fermions

\section{Introduction}
We consider a system of $N$ fermions in the antisymmetric subspace
\begin{align}
\label{eq:antisymmetric_subspace_def}
    L_a^2(\mathbb{R}^{3N}) &\coloneqq \Big\{ \Psi \in L^2(\mathbb{R}^{3N}) :
        \Psi(x_{\pi(1)}, \dots, x_{\pi(N)}) = \operatorname{sgn}(\pi) \, \Psi(x_1, \dots, x_N) \ \text{for all } \pi \in S_N \Big\},
\end{align}
governed by the many-body Hamiltonian
\begin{align}
\label{eq:HN_magnetic}
    H_N &= \sum_{j=1}^N \Big( \tfrac{1}{2} \big( -i\hbar \nabla_{x_j} + A(x_j) \big)^2 + \Vext(x_j) \Big) + \tfrac{1}{2N} \sum_{i \neq j}^N W(x_i - x_j),
\end{align}
with a magnetic vector potential $A$, an external potential $\Vext$ and a two-body interaction potential $W$. The time evolution of a normalized initial state $\Psi_N \in L_a^2(\mathbb{R}^{3N})$ is governed by the Schrödinger equation
\begin{align}
\label{eq:Schrodinger_scaled}
    i\hbar \partial_t \Psi_{N,t} &= H_N \Psi_{N,t}, \qquad \Psi_{N,0} = \Psi_N.
\end{align}

The corresponding classical equation is the magnetic Vlasov equation for an unknown  $\mathbf{f}_t$ on phase space in the position $q$ and the kinetic velocity $v$, named after Vlasov \cite{vlasov1938vibration},
\begin{align}
\label{eq:magnetic_vlasov_intro}
    &\partial_t \mathbf{f}_t(q,v) + v \cdot \nabla_q \mathbf{f}_t(q,v)
        - \Big( v \times B(q) + \nabla_q \Vext(q) + \nabla_q (W * \rho_t)(q) \Big) \cdot \nabla_v \mathbf{f}_t(q,v) = 0,
\end{align}
where the magnetic field $B \coloneqq \nabla \times A$ and the spatial density
\begin{align}
\label{eq:spatial_density_intro}
    \rho_t(q) &\coloneqq \int_{\mathbb{R}^3} dv \, \mathbf{f}_t(q,v).
\end{align}

In this paper, we study the combined mean-field and semiclassical limit $N \to \infty$ of \eqref{eq:Schrodinger_scaled}, with $\hbar = N^{-1/3}$. For the understanding of this scaling for fermions, we refer to references \cite{elgart2004nonlinear, benedikter2014meanfield, fournais2018semiclassical, golse2016mean}.

\subsection{Literature Review}
Many works have been done for the mean-field limit of classical particle system with magnetic field. For example, Ben Porat, Chen and Yuan derived the Euler equations with Lorentz forcing \cite{porat2026supercriticalmeanfieldlimitmagnetized}. For systems coupled with Maxwell equation, we refer to the works of Golse \cite{golse2012mean} and Chen, Li, Pickl, and Yin \cite{chen2020combined}. Further references on this topic were reviewed in \cite{porat2026supercriticalmeanfieldlimitmagnetized} .

On the quantum level, the mean-field limit and semiclassical limit of Schrödinger equations were investigated intensively in the last decades, due to its importance in mathematical physics.
For $A = 0$ the derivation of the Vlasov equation from the many-body Schrödinger dynamics was established through two routes.

The two-step route passes through the Hartree equation. Elgart, Erd\H{o}s, Schlein and Yau derived, in the same scaling $\hbar = N^{-1/3}$, the nonlinear Hartree equation from the $N$-body dynamics by a BBGKY hierarchy with an iteration scheme, for real-analytic interaction potentials, on a fixed time interval and with a rate of order $1/N$ \cite{elgart2004nonlinear}. Benedikter, Porta and Schlein applied the second-quantized framework to fermionic mean-field dynamics \cite{benedikter2014meanfield}, and the semiclassical path from the Hartree dynamics to the Vlasov equation was carried out in \cite{benedikter2016hartree}. Following this path, inverse power laws $|x|^{-\alpha}$, $\alpha \in (0,1/2)$, are reached from the Hartree equation \cite{saffirio2020semiclassical} to the Vlasov equation and from Hartree \cite{saffirio2019hartree} and, for general singular potentials, from Hartree--Fock \cite{lafleche2020strong} to the Vlasov--Poisson equation. For the dynamics \eqref{eq:Schrodinger_scaled} with a non-constant magnetic field, the two-step route is open. Benedikter, Boccato, Monaco and Nguyen derived the Hartree--Fock dynamics from the $N$-body dynamics for a constant magnetic field \cite{benedikter2025magnetic}. For a nonuniform, purely space-dependent magnetic field the second step from the magnetic Hartree--Fock equation to the magnetic Vlasov--Poisson equation was done by Chong and You \cite{chong2026gauge}. Ben Porat obtained the magnetic Liouville equation as the semiclassical limit of the magnetic von Neumann equation \cite{benporat2024liouville}, and the magnetized pressureless Euler--Poisson equation as the semiclassical limit of the magnetized Schr\"odinger--Poisson equation and as the mean-field semiclassical limit of the magnetized von Neumann equation \cite{porat2026magnetizedcoulombicquantumdynamics}.  

The one-step route is due to Narnhofer and Sewell for analytic interaction potentials \cite{narnhofer1981vlasov} and to Spohn for $C^2$ potentials \cite{spohn1981vlasov}. Golse, Mouhot and Paul proved an inequality that renders the mean-field limit uniform in the classical limit \cite{golse2016mean, golse2017schrodinger}. They introduced a quantum analogue of
the Wasserstein distance $\mathbb{W}_2$, which is not a distance on
the density operators, since it is bounded below by $2\,d\,\hbar$. With that pseudo-distance, they get a
mean-field convergence to the Hartree equation with a rate $1/N$, uniform in $\hbar$. Golse and Paul derived in a related pseudo-distance the Vlasov equation
from the $N$-body Schrödinger dynamics \cite{golse2017schrodinger}. With that pseudo-distance Ben
Porat derived the magnetic Liouville equation from the magnetic von Neumann equation
\cite{benporat2024liouville}. 
Chen, Lee and Liew reformulated the many-body Schrödinger equation as a BBGKY-type hierarchy for the $k$-particle Husimi measure and, with estimates on the remainders and compactness argument, obtained the combined limit without intermediate approximations such as the Hartree or Liouville equations \cite{ChenLeeLiew2021}. They proved the
convergence of the $N$-body Schrödinger dynamics to the Vlasov equation in Wasserstein distance $\mathbb{W}_1$. 

We follow the direct path of Chen, Lee and Liew for a problem with a non-constant magnetic field, where the magnetic Husimi measure, constructed from coherent states with a shifted kinetic velocity, plays an important role.      
In contrast to the Wigner transform \cite{lions1993sur, hudson1974when}, the magnetic Husimi measure is a probability density on phase space. The derivation of the magnetic Vlasov equation from the magnetic $N$-body dynamics \eqref{eq:Schrodinger_scaled} in one direct path, for a non-constant magnetic field, is the gap the present work closes.
We prove that the magnetic Husimi measure converges weakly in $L^1$ to the solution of magnetic Vlasov equation \eqref{eq:magnetic_vlasov_intro}. As a corollary, the convergence is valid in Wasserstein distance
$\mathbb{W}_p$ with $1\leq p < 2$. Furthermore, we also obtain the convergence in $\mathbb{W}_2$, with an additional proof of the convergence
of the second moment \cite[Thm.~6.9]{villani2008optimal}. 
\begin{center}
    \definecolor{routeTeal}{RGB}{38,110,110}
    \definecolor{mattRed}{RGB}{176,72,60}
    \hypersetup{citecolor=routeTeal}
    \resizebox{\textwidth}{!}{%
    \begin{tikzpicture}[>=Stealth,
        lbox/.style={draw=routeTeal, thick, rounded corners, fill=routeTeal!8, align=center,
                     text width=3.0cm, minimum height=1.2cm, font=\small},
        arr/.style={->, thick, routeTeal},
        lab/.style={font=\footnotesize, text=routeTeal, align=center},
        ]
        \node[font=\bfseries, text=routeTeal] at (3.5, 1.7) {No magnetic field, $A = 0$};
        \node[lbox] (NB) at (0, 0) {$N$-body\\ Schrödinger};
        \node[lbox] (VL) at (7, 0) {Vlasov};
        \node[lbox] (HF) at (3.5, -2.2) {Hartree,\\ Hartree--Fock};
        \draw[arr, line width=1.6pt] (NB) -- node[lab, above] {\cite{narnhofer1981vlasov, spohn1981vlasov, golse2016mean, ChenLeeLiew2021}} node[lab, below] {one step} (VL);
        \draw[arr] (NB) -- node[lab, left=2pt] {\cite{elgart2004nonlinear, benedikter2014meanfield}} (HF);
        \draw[arr] (HF) -- node[lab, anchor=north west, pos=0.55] {\cite{benedikter2016hartree}\\ singular $W$:\\ \cite{saffirio2020semiclassical, lafleche2020strong, saffirio2019hartree}} (VL);
        \node[font=\footnotesize, text=routeTeal, align=center] at (7, 1.05) {well-posedness: \cite{dobrushin1979vlasov}};

        \node[font=\bfseries, text=mattRed] at (14.5, 1.7) {Non-constant magnetic field};
        \node[lbox] (MNB) at (11, 0) {magnetic $N$-body\\ Schrödinger};
        \node[lbox, draw=mattRed, fill=mattRed!10, line width=1.4pt] (MVL) at (18, 0) {magnetic\\ Vlasov};
        \node[lbox] (MHF) at (14.5, -2.2) {magnetic\\ Hartree--Fock};
        \draw[->, line width=2.6pt, mattRed] (MNB) -- node[fill=mattRed, text=white, font=\small\bfseries, rounded corners=2pt, inner sep=3pt, above=3pt] {This Work} node[font=\footnotesize\bfseries, text=mattRed, below=2pt] {one step} (MVL);
        \draw[->, thick, dashed, black!50] (MNB) -- node[font=\footnotesize, text=black!60, anchor=north east, pos=0.45, align=right] {\hypersetup{citecolor=black!60}\cite{benedikter2025magnetic}\\ constant $B$ only} (MHF);
        \draw[arr] (MHF) -- node[lab, anchor=north west, pos=0.55] {\cite{chong2026gauge}\\ to Vlasov--\\ Poisson} (MVL);
        \node[font=\footnotesize, text=mattRed, align=center] at (18, 1.15) {well-posedness: \cite{rege2023magnetized, rege2024stability, benporat2025moments}\\ Vlasov--Poisson only};
    \end{tikzpicture}%
    }
    
    {\small TikZ code drafted with the assistance of AI.}
\end{center}

\subsection{Setting and the Main Result}
 On the fermionic Fock space $\mathcal{F}_a = \bigoplus_{n \ge 0} L_a^2(\mathbb{R}^{3n})$, with operator-valued distributions $a_x, a_x^\dagger$ satisfying $\{a_x, a_y^\dagger\} = \delta(x-y)$ and number operator $\mathcal{N} = \displaystyle \int_{\mathbb{R}^3} dx \, a_x^\dagger a_x$, the Hamiltonian $H_N$ is the restriction to the $N$-particle sector of
\begin{align}
\label{eq:H_fock_def}
    \mathcal{H} &= \mathcal{K}_A + \mathcal\Vext + \mathcal{W}, \qquad
    \mathcal{K}_A = \tfrac{1}{2} \int_{\mathbb{R}^3} dx \, a_x^\dagger \mpi_x^2 a_x, \qquad
    \mpi_x \coloneqq -i\hbar\nabla_x + A(x), \nonumber \\
    \mathcal\Vext &= \int_{\mathbb{R}^3} dx \, \Vext(x) \, a_x^\dagger a_x, \qquad
    \mathcal{W} = \tfrac{1}{2N} \iint_{\mathbb{R}^3 \times \mathbb{R}^3} dx \, dy \, W(x-y) \, a_x^\dagger a_y^\dagger a_y a_x.
\end{align}

\vskip3mm
The potentials $A$, $\Vext$ and $W$ are assumed to satisfy the following assumptions.

\begin{assumptionA}
\label{assump:A1}
The magnetic vector potential $A \in C^1(\mathbb{R}^3; \mathbb{R}^3)$ satisfies $\nabla A \in W^{1,\infty}(\mathbb{R}^3; \mathbb{R}^{3 \times 3})$.
\end{assumptionA}

\begin{assumptionA}
\label{assump:A2}
The external potential $\Vext \in C^2(\mathbb{R}^3; \mathbb{R})$ satisfies $\nabla^2 \Vext \in L^\infty(\mathbb{R}^3; \mathbb{R}^{3 \times 3})$ and $\Vext \ge -c_1$ for a constant $c_1 \ge 0$.
\end{assumptionA}

\begin{assumptionA}
\label{assump:A3}
The two-body interaction potential $W \in W^{2,\infty}(\mathbb{R}^3; \mathbb{R})$ satisfies $W(-x) = W(x)$.
\end{assumptionA}

Assumption \ref{assump:A2} admits $\Vext \equiv 0$ and the harmonic trap $\Vext(x) = \tfrac{1}{2} \omega^2 |x|^2$,
which does not lie in $W^{2,\infty}(\mathbb{R}^3)$.

The semiclassical coherent states given in \cite{fournais2018semiclassical} are centred at a position $q_j$ and a canonical momentum $p_j$. In the magnetic Vlasov equation \eqref{eq:magnetic_vlasov_intro} is written in the kinetic velocity, so the coherent states are centred at $v_j = p_j + A(q_j)$. Therefore, it is natural to use the velocity-shifted coherent state of a particle $j$. Hence, the $k$-particle magnetic Husimi measure is given with the phase-space coordinate $(q_j, v_j)$ by
\begin{align}
\label{eq:coherent_shifted}
    \mathfrak{f}_j (x) &\coloneqq f^\hbar_{q_j,v_j-A(q_j)}(x) = \hbar^{-3/4} f \Big( \tfrac{x - q_j}{\sqrt{\hbar}} \Big) e^{\frac{i}{\hbar} (v_j - A(q_j)) \cdot x}.
\end{align}
where we assume that the profile $f$ satisfies:
\begin{assumptionA}
\label{assump:A4}
The profile $f \in H^2(\mathbb{R}^3; \mathbb{R})$ satisfies
\begin{align}
\label{eq:assump_A4_profile}
    \|f\|_{L^2(\mathbb{R}^3)} &= 1 \quad \text{and} \quad \big\| |\cdot|^2 f \big\|_{L^2(\mathbb{R}^3)} + \big\| |\cdot|^2 \nabla f \big\|_{L^2(\mathbb{R}^3)} < \infty .
\end{align}
\end{assumptionA}
The profile of these coherent states needs no compact support, in contrast to \cite[Assumption~A1]{ChenLeeLiew2021}, so Gaussian coherent states are admitted.

It is obvious that the resolution of the identity on $L^2(\mathbb{R}^3)$ is preserved under the velocity shift:
\begin{align}
\label{eq:resolution_identity}
    \tfrac{1}{(2\pi\hbar)^3} \iint_{\mathbb{R}^6} dq_j \, dv_j \, \left| \mathfrak{f}_j \right\rangle \left\langle \mathfrak{f}_j \right| &= \mathbb{I}_{L^2(\mathbb{R}^3)}.
\end{align}

The gradients of the velocity-shifted coherent state are:
\begin{align}
\label{eq:state_level_chain_rule}
    \nabla_{v_j} \mathfrak{f}_j(x_j) &= \tfrac{i x_j}{\hbar} \mathfrak{f}_j(x_j), \quad
    \nabla_{q_j}\Big|_{p_j} \mathfrak{f}_j(x_j) = \left( \nabla_{q_j} + (\nabla_{q_j} A(q_j))^T \nabla_{v_j} \right) \mathfrak{f}_j(x_j) .
\end{align}
where the notation $\nabla_{q_j}\big|_{p_j}$ means the gradient in $q_j$ at fixed $p_j = v_j - A(q_j)$.

We use the standard notation for the annihilation operator $a(g) = \displaystyle \int_{\mathbb{R}^3} dx \, \overline{g(x)} \, a_x$ for $g \in L^2(\mathbb{R}^3)$. It follows from \eqref{eq:resolution_identity} that
\begin{align}
\label{eq:phase_space_number_operators}
       \tfrac{1}{(2\pi\hbar)^3} \iint_{\mathbb{R}^6} dq_j \, dv_j \, a^\dagger(\mathfrak{f}_j) \, a(\mathfrak{f}_j) = \mathcal{N}.
\end{align}
The $k$-particle annihilation operator and its variant with the $j$-th factor replaced by $a(g)$ read
\begin{align}
\label{eq:operator_strings}
    \mathfrak{A}_k[\mathfrak{f}] &\coloneqq a(\mathfrak{f}_k) \cdots a(\mathfrak{f}_1), \quad
    \mathfrak{A}_{k,j}[g] \coloneqq a(\mathfrak{f}_k) \cdots a(\mathfrak{f}_{j+1}) \, a(g) \, a(\mathfrak{f}_{j-1}) \cdots a(\mathfrak{f}_1).
\end{align}
The $k$-particle magnetic Husimi measure takes the velocity-shifted coherent states \eqref{eq:coherent_shifted} in their place and is normalized in the $k$-particle sector by the sector ratio
\begin{align}
\label{eq:kinetic_husimi_M_def}
  \mhus^{(k)}_{N,t}(Z_k) &\coloneqq \tfrac{c_{N,k}}{(2\pi)^{3k}} \left\langle \Psi_{N,t}, \mathfrak{A}_k^\dagger[\mathfrak{f}] \, \mathfrak{A}_k[\mathfrak{f}] \Psi_{N,t} \right\rangle, \qquad  c_{N,k} \coloneqq \tfrac{N^k \, (N-k)!}{N!},
\end{align}
with $Z_k \coloneqq (q_1, v_1, \dots, q_k, v_k)$.
Up to the factor $(2\pi)^{-3k}c_{N,k}$, the magnetic Husimi measure $\mhus^{(k)}_{N,t}$ coincides with \cite[Eq.~(1.15)]{fournais2018semiclassical} at $p_j = v_j - A(q_j)$, and \cite[Lemma~2.2, p.~14]{fournais2018semiclassical} implies that: $\mhus^{(k)}_{N,t}$ is symmetric and, for every $N$, $t$ and $1 \le k \le N$,
\begin{align}
\label{eq:normalized_husimi_measure}
    &0 \le \mhus_{N,t}^{(k)}(Z_k) \le \tfrac{c_{N,k}}{(2\pi)^{3k}} \,\,\text{a.e.} ,\quad \int_{\mathbb{R}^{6k}} \dZk \, \mhus_{N,t}^{(k)}(Z_k) = 1 ,\quad  \mhus_{N,t}^{(k)} \in \mathcal{P}(\mathbb{R}^{6k}) , \notag\\
    &\iint_{\mathbb{R}^6} dq_{k+1} \, dv_{k+1} \, \mhus_{N,t}^{(k+1)}(Z_k, q_{k+1}, v_{k+1})
    = \mhus_{N,t}^{(k)}(Z_k) , \qquad k \le N-1.
\end{align}
The main purpose of this paper is to prove that for any fixed $k$, the $k$-particle magnetic Husimi measure $\mhus_{N,t}^{(k)}$ converges to $\mathbf{f}_t^{\otimes k}$, the solution of \eqref{eq:magnetic_vlasov_intro}.

The well-posedness of \eqref{eq:magnetic_vlasov_intro} can be obtained by direct application of the classical method \cite{golse2013empirical, dobrushin1979vlasov} by using following assumptions for initial data:
\begin{assumptionI}[Classical Initial Data]
\label{assump:I0}
The initial density $\mathbf{f}_0 \in L^\infty(\mathbb{R}^6)$ is a probability density with compact support.
\end{assumptionI}

\begin{proposition}[Well-Posedness of \eqref{eq:magnetic_vlasov_intro}] \label{thm:vlasov_wellposedness}
Let Assumptions \ref{assump:A1}--\ref{assump:A3} and \ref{assump:I0} hold, and let $T > 0$. Then
\eqref{eq:magnetic_vlasov_intro} has a unique weak solution $(\mathbf{f}_t)_{t \in [0,T]}$ with initial data
$\mathbf{f}_0$, $\|\mathbf{f}_t\|_{L^\infty} = \|\mathbf{f}_0\|_{L^\infty}$ and with compact support for every $t \in [0,T]$. Namely, there is $R(t) > 0$, nondecreasing in $t$, such that
\begin{align}
\label{eq:vlasov_support}
    \operatorname{supp}\mathbf{f}_t &\subset \big\{ |q|^2 + |v|^2 \le R(t)^2 \big\} , \qquad t \in [0,T] .
\end{align}
\end{proposition}

Further assumptions on initial data for the proof of the main results are: the boundedness of the total energy and the second spatial moment,
\begin{assumptionI}[Initial Many-Body State]
\label{assump:I1}
The initial states $\Psi_{N,0} \in L^2_a(\mathbb{R}^{3N})$ satisfy, with a constant $\C > 0$ independent of $N$,
\begin{align}
\label{eq:assumption_initial_data}
    \tfrac{1}{N} \big\langle \Psi_{N,0}, \mathcal{H} \Psi_{N,0} \big\rangle
    &\le\C,\qquad \qquad
    \iint_{\mathbb{R}^6} dq \, dv \, |q|^2 \, \mhus_{N,0}^{(1)}(q,v)
    \le \C.
\end{align}
\end{assumptionI}
and the asymptotic factorization of the initial magnetic Husimi measures in every sector, 
\begin{assumptionI}[Factorized Initial Data]
\label{assump:I2}
With $\mathbf{f}_0$ of Assumption \ref{assump:I0}, it holds
\begin{align}
\label{eq:initial_condition_tensor}
    \mhus_{N,0}^{(k)}
    &\; \xrightharpoonup[\hbar = N^{-1/3}]{\hbar \to 0} \; \mathbf{f}_0^{\otimes k}
    \quad \text{weakly in } L^1(\mathbb{R}^{6k}),
     \qquad k \ge 1 .
\end{align}
\end{assumptionI}

\vskip3mm
The main result of this paper reads:
\begin{theorem}[Convergence to the Magnetic Vlasov Equation]
\label{thm:magnetic_vlasov_limit}
Let Assumptions \ref{assump:A1}--\ref{assump:A4} and \ref{assump:I0}--\ref{assump:I2} hold, and let $T > 0$. Then it holds, for every $k \ge 1$ and $t \in [0,T]$,
\begin{align}
\label{eq:full_sequence_convergence}
    \mhus_{N,t}^{(k)} &\xrightharpoonup[\hbar = N^{-1/3}]{N \to \infty}
    \mathbf{f}_t^{\otimes k} \quad \text{weakly in } L^1(\mathbb{R}^{6k}) ,
\end{align}
where $\mathbf{f}_t$ is given in Proposition
\ref{thm:vlasov_wellposedness}. 
\end{theorem}

A direct consequence of Theorem \ref{thm:magnetic_vlasov_limit} is the following convergence:

\begin{corollary}[Convergence in Wasserstein $\mathbb{W}_p$, $p \in [1,2)$]
\label{cor:wasserstein_p_convergence}
Let Assumptions \ref{assump:A1}--\ref{assump:A4} and \ref{assump:I0}--\ref{assump:I2} hold, and let
$T > 0$. Then for every $k \ge 1$, and $t \in [0,T]$ it holds
\begin{align}
\label{eq:wasserstein_p_convergence}
    \mathbb{W}_p\left( \mhus_{N,t}^{(k)}, \mathbf{f}_t^{\otimes k} \right)
    &\; \xrightarrow[\hbar = N^{-1/3}]{\hbar \to 0} \; 0 , \quad\text{for}\quad p \in [1,2).
\end{align}
\end{corollary}
The current density of $\Psi_{N,t}$ at the position $x$ and the current density of the weak solution are
\begin{align}
\label{eq:quantum_current_def}
    \mathbf{j}_{\Psi_{N,t}}(x) &\coloneqq \operatorname{Re} \langle a_x \Psi_{N,t}, \mpi_x a_x \Psi_{N,t} \rangle ,
    \qquad
    \mathbf{J}_{t}(q) \coloneqq \int_{\mathbb{R}^3} dv \, v \, \mathbf{f}_t(q,v) .
\end{align}
As a consequence of Corollary \ref{cor:wasserstein_p_convergence}, the current density of the $N$-body state converges:

\begin{corollary}[Convergence of the Current Density]
\label{cor:quantitative_transport_rate_noV}
Let Assumptions \ref{assump:A1}--\ref{assump:A4} and \ref{assump:I0}--\ref{assump:I2} hold, let
$T > 0$. Then for every $t \in [0,T]$ and
$\phi \in W^{1,\infty}(\mathbb{R}^3)$ it holds
\begin{align}
\label{eq:current_convergence_bound_cor_noV}
    \left| \int_{\mathbb{R}^3} dq \, \phi(q) \Big( \tfrac{1}{N} \mathbf{j}_{\Psi_{N,t}}(q) - \mathbf{J}_t(q) \Big) \right| &\le \C (1+t)^2 \, \Big( \mathbb{W}_1\left(\mhus_{N,t}^{(1)}, \mathbf{f}_t\right) + \sqrt{\hbar} \Big) \; \xrightarrow[\hbar = N^{-1/3}]{\hbar \to 0} \; 0 ,
\end{align}
with $\C$ independent of $N$, $\hbar$ and $t$.
\end{corollary}

Moreover, the convergence result can even be achieved in the Wasserstein distance $\mathbb{W}_2$. Golse, Mouhot and Paul obtained the combined limit with a quantum analogue of $\mathbb{W}_2$ \cite{golse2016mean, golse2017schrodinger}. The magnetic Husimi measure is a probability measure on the phase space, where $\mathbf{f}_t$ lives, so $\mathbb{W}_2$ applies in its classical form \cite[Def.~6.1]{villani2008optimal}. 

We consider the total energy of the magnetic Vlasov equation and give a further assumption on the initial data:
\begin{align}
\label{eq:classical_energy_definition}
    \mathcal{E}[\mathbf{f}] &\coloneqq \iint_{\mathbb{R}^6} dq \, dv \, \Big( \tfrac{1}{2} |v|^2 + \Vext(q) \Big) \mathbf{f}(q,v)
        + \frac{1}{2} \int_{\mathbb{R}^{12}} (dq \, dv)^{\otimes 2} \, W(q_1 - q_2) \, \mathbf{f}(q_1,v_1) \, \mathbf{f}(q_2,v_2) .
\end{align}

\begin{assumptionI}[Spatial Moment and Energy at $t = 0$]
\label{assump:I4}
With $\mathbf{f}_0$ of Assumption \ref{assump:I0}, it holds
\begin{align}
\label{eq:assumption_moment_convergence}
    & \lim\nolimits_{N \to \infty} \iint_{\mathbb{R}^6} dq \, dv \, |q|^2 \, \mhus_{N,0}^{(1)}(q,v)
    = \iint_{\mathbb{R}^6} dq \, dv \, |q|^2 \, \mathbf{f}_0(q,v) , \\
\label{eq:criterion_energy_noV}
    & \limsup\nolimits_{N \to \infty} \frac{1}{N} \big\langle \Psi_{N,0}, \mathcal{H} \Psi_{N,0} \big\rangle
    \le \mathcal{E}[\mathbf{f}_0] .
\end{align}
\end{assumptionI}
The reverse inequality $\liminf_{N \to \infty} \frac{1}{N} \langle \Psi_{N,0}, \mathcal{H} \Psi_{N,0} \rangle \ge \mathcal{E}[\mathbf{f}_0]$ holds under Assumptions \ref{assump:A1}--\ref{assump:A4}, \ref{assump:I1} and \ref{assump:I2}, so \eqref{eq:criterion_energy_noV} states the convergence of the initial energy per particle to $\mathcal{E}[\mathbf{f}_0]$.

\vskip3mm
The second main result of this paper is:
\begin{theorem}[Convergence in Wasserstein $\mathbb{W}_2$]
\label{thm:w2_energy_criterion_noV}
Let Assumptions \ref{assump:A1}--\ref{assump:A4} and \ref{assump:I0}--\ref{assump:I4} hold, and let $T > 0$. Then it
holds, for every $k \ge 1$ and $t \in [0,T]$,
\begin{align}
\label{eq:criterion_transport_noV}
    \mathbb{W}_2\Big( \mhus_{N,t}^{(k)}, \mathbf{f}_t^{\otimes k} \Big)
    \; \xrightarrow[\hbar = N^{-1/3}]{\hbar \to 0} \; 0 .
\end{align}
\end{theorem}

\paragraph{Discussion of the initial data.} Theorems \ref{thm:magnetic_vlasov_limit} and \ref{thm:w2_energy_criterion_noV}
rest on assumptions on the initial states $\Psi_{N,0}$ and on their magnetic Husimi measures. The next remark shows
that the set of initial states which satisfy these assumptions is not empty: all of them hold for the Slater determinants (defined as in \cite{elgart2004nonlinear}) of trapped ground states.

\begin{remark}[Slater Determinants as Initial Data]
\label{rem:initial_data_ground_states}
The hypotheses of Theorems \ref{thm:magnetic_vlasov_limit} and
\ref{thm:w2_energy_criterion_noV} hold for the Slater determinants of trapped ground states, a case of the initial states of
\cite[Remark~1.5]{ChenLeeLiew2021}. More precisely, let $\Psi_{N,0}$ be the Slater determinant of the $N$ lowest eigenfunctions of
$- \tfrac12 \hbar^2 \Delta + \Vext_0$, where $\Vext_0$ satisfies Assumption \ref{assump:A2} and $\Vext_0(x) \ge c_0 |x|^2 - C_0$,
$c_0 > 0$. 

Initially, the trap $\Vext_0$ prepares the state, at $t = 0$ it is switched off, and
$\Psi_{N,0}$ evolves by \eqref{eq:Schrodinger_scaled} with the given Hamiltonian $H_N$. In Appendix \ref{app:slater_determinants} we prove:
\begin{enumerate}
\renewcommand{\theenumi}{(\roman{enumi})}
\renewcommand{\labelenumi}{\theenumi}
    \item Assumption \ref{assump:I1} holds
    (the result is given in Proposition \ref{prop:slater_hypotheses}).
   \item Assumptions \ref{assump:I0} and \ref{assump:I2} hold, with the initial data
    \begin{align}
    \label{eq:intro_ground_state_f0}
        \mathbf{f}_0(q,v) &= (2\pi)^{-3} \, \One\Big( |v - A(q)|^2 + 2\Vext_0(q) \le s \Big) ,
    \end{align}
which is the minimizer of $\displaystyle \iint_{\mathbb{R}^6} dq \, dv \, \big( \tfrac12 |v - A(q)|^2 + \Vext_0(q) \big) \mathbf{f}$ under the
Pauli bound $0 \le \mathbf{f} \le (2\pi)^{-3}$, with $s$ determined by
    $\displaystyle \iint_{\mathbb{R}^6} dq \, dv \, \mathbf{f}_0 = 1$ (the result is also given in Proposition \ref{prop:slater_hypotheses}).
        \item \label{item:slater_I4} Assumption \ref{assump:I4} holds (the result is given in Proposition \ref{prop:slater_I4}).
    
\end{enumerate}
\end{remark}

\subsection{Proof Strategy} The main contribution of the present paper is Theorem \ref{thm:magnetic_vlasov_limit}. Its proof rests on three steps. The {\bf first step} is an identity extending the BBGKY hierarchy of \cite{ChenLeeLiew2021} to the magnetic case.

\begin{theorem}[$k$-Particle Magnetic BBGKY Hierarchy in Phase Space]
\label{thm:magnetic_hierarchy_k}
Let Assumptions \ref{assump:A1}--\ref{assump:A4} hold. For any $1 \le k \le N$, the $k$-particle magnetic Husimi measure $\mhus_{N,t}^{(k)}$ satisfies, at every $Z_k \in \mathbb{R}^{6k}$,
\begin{align}
\label{eq:magnetic_hierarchy_k}
    &\partial_t \mhus_{N,t}^{(k)} + \sum_{j=1}^k \Big( v_j \cdot \nabla_{q_j} \mhus_{N,t}^{(k)} - \big( v_j \times B(q_j) + \nabla \Vext(q_j) \big) \cdot \nabla_{v_j} \mhus_{N,t}^{(k)}
        - \nabla_{v_j} \cdot \mathcal{F}_{\text{int}, j}\big[\mhus_{N,t}^{(k+1)}\big] \Big) = \mathbf{R}^{(k)} ,
\end{align}
with the mean-field force on the $j$-th particle
\begin{align}
\label{eq:mean_field_force_def_k}
    \mathcal{F}_{\text{int}, j}[\mhus_{N,t}^{(k+1)}](Z_k) &\coloneqq \tfrac{N-k}{N} \iint_{\mathbb{R}^6} dq_{k+1} \, dv_{k+1} \, \nabla W(q_j - q_{k+1}) \, \mhus_{N,t}^{(k+1)}(Z_k, q_{k+1}, v_{k+1}) ,
\end{align}
and the quantum remainder $\mathbf{R}^{(k)}$ given explicitly in Section \ref{sec:main_results}.
\end{theorem}

The quantum remainder $\mathbf{R}^{(k)}$ collects all terms beyond the classical magnetic transport. Its closed form is a sum of a magnetic remainder, itself split into five parts, and of the external, mean-field and internal-correlation terms, each an explicit expectation of field operators in the state $\Psi_{N,t}$. The proof is given in Section \ref{sec:main_results}.

The {\bf second step} is a bound on that remainder, uniform in $N, t$. 
\begin{theorem}[Uniform Semiclassical Remainder Estimates]
\label{thm:unified_remainders_k}
Let Assumptions \ref{assump:A1}--\ref{assump:A4} hold. Let $\Phi \in C_c^\infty(\mathbb{R}^{6k})$. For every $\rate \in \left(0, \frac{1}{2}\right)$ and $1 \le k \le N$:
\begin{align}
\label{eq:unified_bound_R_k}
  \left| \int_{\mathbb{R}^{6k}} \dZk \, \Phi(Z_k) \mathbf{R}^{(k)}(Z_k) \right| &\le \C(k) \hbar^{\frac{1}{2} - \rate} \xrightarrow[\hbar = N^{-1/3}]{\hbar \to 0} 0,
\end{align}
where $\C(k) < \infty$ is independent of $N$, $\hbar$, and $t$.
\end{theorem}
The proof in Section \ref{sec:remainder_bounds} bounds $\mathbf{R}^{(k)}$ component by component. The exponent $\tfrac{1}{2}$ is the ceiling set by the dispersion $\sqrt{\hbar}$ of the coherent state, and the rate $\hbar^{\frac{1}{2} - \rate}$ stays below it. The obstacle in the second step stems from the divergence of $A$. Acting on a velocity-shifted coherent state, the commutator $[\mpi_{x_j}, A(x_j)] = -i\hbar \, (\nabla \cdot A)(x_j)$ leaves $\nabla \cdot A(q_j)$ at the centre $q_j$ of the coherent state, and the vector potential, linearized at $q_j$, produces in the scaled variable $z = (x_j - q_j)/\sqrt{\hbar}$ the term $(\nabla A(q_j) z) \cdot \nabla f(z)$. Both are of order one in $\hbar$, and neither vanishes as $\hbar \to 0$. The projection of $(\nabla A(q_j) z) \cdot \nabla f(z)$ onto $f$ cancels $\nabla \cdot A(q_j)$. The residual is orthogonal to $f$, and this orthogonality supplies the rate of Theorem \ref{thm:unified_remainders_k} (Section \ref{subsubsec:bound_I_div_j_k}).

The {\bf third step} in Section \ref{sec:combined_limit} provides the a priori bounds on the magnetic Husimi measures and shows that they are weakly compact in $L^1$. Furthermore, we prove that any convergent subsequence converges to an admissible solution of the corresponding infinite magnetic hierarchy, which is defined as follows.
\begin{definition}[Admissible Solutions]
\label{def:admissible_class}
Let $T > 0$. The class $\mathcal{A}_T$ consists of the families $\mathbf{f} = (\mathbf{f}^{(k)}_t)$,
indexed by $k \ge 1$ and $t \in [0,T]$, with the following properties.
\begin{enumerate}
\renewcommand{\theenumi}{(\alph{enumi})}
\renewcommand{\labelenumi}{\theenumi}
\item \label{item:admissible_a}
Each $\mathbf{f}^{(k)}_t$ is a symmetric probability density on $\mathbb{R}^{6k}$, and the family is consistent:
\begin{align}
\label{eq:admissible_consistency}
    \iint_{\mathbb{R}^6} dq_{k+1} \, dv_{k+1} \, \mathbf{f}^{(k+1)}_t &= \mathbf{f}^{(k)}_t .
\end{align}
\item \label{item:admissible_b}
For every $k \ge 1$ and every $\varphi \in C_b(\mathbb{R}^{6k})$ the map
\begin{align}
\label{eq:admissible_continuity}
    t \;\mapsto\; \int_{\mathbb{R}^{6k}} (dq \, dv)^{\otimes k} \, \varphi(Z_k) \, \mathbf{f}^{(k)}_t(Z_k)
\end{align}
is continuous on $[0,T]$.
\item \label{item:admissible_c}
The second phase-space moment is bounded,
\begin{align}
\label{eq:admissible_moments}
    \sup_{t \in [0,T]} \iint_{\mathbb{R}^6} dq \, dv \, \big( |q|^2 + |v|^2 \big) \, \mathbf{f}^{(1)}_t(q,v) &< \infty .
\end{align}
\item \label{item:admissible_d}
For every $k \ge 1$ and $\Phi \in C_c^\infty([0,T) \times \mathbb{R}^{6k})$,
\begin{align}
\label{eq:hierarchy_weak_form}
    &\int_0^T dt \int_{\mathbb{R}^{6k}} (dq \, dv)^{\otimes k} \, \mathbf{f}^{(k)}_t \bigg( \partial_t \Phi
        + \sum\nolimits_{j=1}^k \Big( v_j \cdot \nabla_{q_j} \Phi - \big( v_j \times B(q_j) + \nabla \Vext(q_j) \big) \cdot \nabla_{v_j} \Phi \Big) \bigg) \nonumber \\
    &\quad - \sum\nolimits_{j=1}^k \int_0^T dt \int_{\mathbb{R}^{6(k+1)}} (dq \, dv)^{\otimes (k+1)}
        \nabla_{v_j} \Phi(t, Z_k) \cdot \nabla W(q_j - q_{k+1}) \, \mathbf{f}^{(k+1)}_t(Z_{k+1}) \nonumber \\
    &\quad + \int_{\mathbb{R}^{6k}} (dq \, dv)^{\otimes k} \, \Phi(0, Z_k) \, \mathbf{f}^{(k)}_0(Z_k) = 0 .
\end{align}

\end{enumerate}
A family $\mathbf{f} = (\mathbf{f}^{(k)}_t) \in \mathcal{A}_T$ is an \emph{admissible solution} of the infinite
magnetic hierarchy \eqref{eq:hierarchy_weak_form} with initial data $(\mathbf{f}^{(k)}_0)_{k \ge 1}$.
\end{definition}
Actually, we can follow and adapt the idea of Spohn \cite{golse2013empirical, spohn1981vlasov} to show the uniqueness of the admissible solution. For completeness, we list the result:
\begin{proposition}[The Infinite Magnetic Hierarchy] \label{thm:vlasov_infinite_hierarchy}
Let Assumptions \ref{assump:A1}--\ref{assump:A3} and \ref{assump:I0} hold, and let $T > 0$. Then the infinite magnetic hierarchy \eqref{eq:hierarchy_weak_form} has exactly one solution in $\mathcal{A}_T$ with given admissible initial data
$(\mathbf{f}_0^{\otimes k})_{k \ge 1}$, namely the only solution is given by $(\mathbf{f}_t^{\otimes k})_{k \ge 1}$ with $\mathbf{f}_t$ from
Proposition \ref{thm:vlasov_wellposedness}.
\end{proposition}

\begin{remark}
    Assumption \ref{assump:I0} enters only through Proposition \ref{thm:vlasov_wellposedness}. Without it, Proposition \ref{prop:main_magnetic_limit} still gives convergence of subsequences: under Assumptions \ref{assump:A1}--\ref{assump:A4} and \ref{assump:I1} alone, the magnetic Husimi measures converge weakly in $L^1$, in every sector and for a subsequence, to an admissible solution of the infinite magnetic hierarchy \eqref{eq:hierarchy_weak_form}.
\end{remark}

Therefore, the result in Theorem \ref{thm:magnetic_vlasov_limit} is obtained. This weak convergence in $L^1$ yields narrow convergence, tested with bounded continuous functions. In order to prove the convergence in $\mathbb{W}_2$ in Theorem \ref{thm:w2_energy_criterion_noV}, it is enough to show the convergence of the second moment of $\mhus_{N,t}^{(1)}$ to that of $\mathbf{f}_t$ \cite[Thm.~6.9]{villani2008optimal}. For this, starting from the convergence of the initial spatial $|q|^2$ moment, we control the spatial $|q|^2$ moment by the velocity $|v|^2$ moment, which converges by energy conservation and the convergence of the initial energy, since the difference of the quantum and the classical energy vanishes as $\hbar \to 0$.

\subsection{Organization of the Paper}
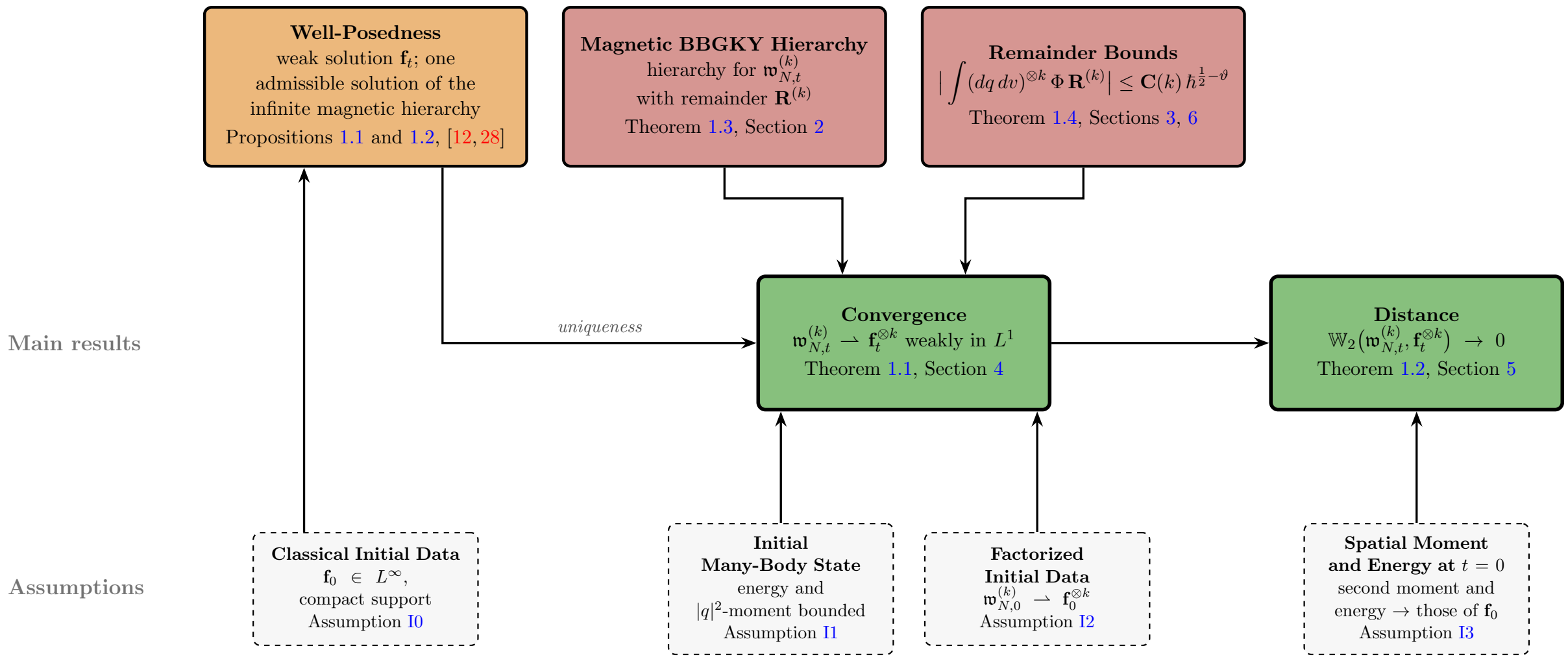
\begin{figure}[H]
    \centering
    \definecolor{coreGreen}{RGB}{135,192,125}
    \definecolor{mattRed}{RGB}{214,150,146}
    \definecolor{mattOrange}{RGB}{236,184,124}
    \resizebox{\textwidth}{!}{%
    \begin{tikzpicture}[>=Stealth,
        cbox/.style={draw=black, rounded corners, align=center, inner sep=1.6ex, minimum height=2.6cm,
                     text width=5.2cm, font=\fontsize{11}{14}\selectfont,
                     execute at begin node={\hyphenpenalty=10000\exhyphenpenalty=10000}},
                stepbox/.style={cbox, fill=mattRed, line width=1.6pt, text width=5.8cm, minimum height=3.1cm},
                litbox/.style={cbox, fill=mattOrange, line width=1.6pt, text width=5.8cm, minimum height=3.1cm},
        main/.style={cbox, fill=coreGreen, line width=2pt},
        assum/.style={cbox, dashed, thick, fill=black!3, text width=3.9cm, minimum height=2.2cm,
                      font=\fontsize{10}{12.5}\selectfont},
        lvl/.style={font=\bfseries\large, text=black!55, anchor=west},
        arr/.style={->, line width=1.2pt},
        lab/.style={font=\fontsize{10}{12}\selectfont\itshape, text=black!70, fill=white, inner sep=1.5pt},
        ]
        \node[lvl] at (-8.6, -5.0)  {Main results};
        \node[lvl] at (-8.6, -9.8)  {Assumptions};

        \node[litbox] (T1) at (-1.5, 0)
            {\textbf{Well-Posedness}\\ weak solution $\mathbf{f}_t$; one admissible solution of the infinite magnetic hierarchy\\[0.4ex]
             Propositions \ref{thm:vlasov_wellposedness} and \ref{thm:vlasov_infinite_hierarchy}, \cite{dobrushin1979vlasov, spohn1981vlasov}};

        \node[stepbox] (T4) at (5.5, 0)
            {\textbf{Magnetic BBGKY Hierarchy}\\ hierarchy for $\mhus^{(k)}_{N,t}$ with remainder $\mathbf{R}^{(k)}$\\[0.4ex]
             Theorem \ref{thm:magnetic_hierarchy_k}, Section \ref{sec:main_results}};
        \node[stepbox] (T5) at (12.5, 0)
            {\textbf{Remainder Bounds}\\ $\big| \displaystyle \int \dZk \, \Phi \, \mathbf{R}^{(k)} \big| \le \C(k) \, \hbar^{\frac{1}{2} - \rate}$\\[0.4ex]
             Theorem \ref{thm:unified_remainders_k}, Sections \ref{sec:remainder_bounds}, \ref{sec:unified_bounds}};

        \node[main] (T2) at (9.0, -5.0)
            {\textbf{Convergence}\\ $\mhus^{(k)}_{N,t} \rightharpoonup \mathbf{f}_t^{\otimes k}$ weakly in $L^1$\\[0.4ex]
             Theorem \ref{thm:magnetic_vlasov_limit}, Section \ref{sec:combined_limit}};
        \node[main] (T3) at (19.0, -5.0)
            {\textbf{Distance}\\ $\mathbb{W}_2\big(\mhus^{(k)}_{N,t}, \mathbf{f}_t^{\otimes k}\big) \to 0$\\[0.4ex]
             Theorem \ref{thm:w2_energy_criterion_noV}, Section \ref{sec:wasserstein_noV}};

        \node[assum] (I0) at (-1.5, -9.8)
            {\textbf{Classical Initial Data}\\ $\mathbf{f}_0 \in L^\infty$, compact support\\ Assumption \ref{assump:I0}};
        \node[assum] (I1) at (6.6, -9.8)
            {\textbf{Initial Many-Body State}\\ energy and $|q|^2$-moment bounded\\ Assumption \ref{assump:I1}};
        \node[assum] (I2) at (11.6, -9.8)
            {\textbf{Factorized Initial Data}\\ $\mhus^{(k)}_{N,0} \rightharpoonup \mathbf{f}_0^{\otimes k}$\\ Assumption \ref{assump:I2}};
        \node[assum] (I3) at (19.0, -9.8)
            {\textbf{Spatial Moment and Energy at $t = 0$}\\ second moment and energy $\to$ those of $\mathbf{f}_0$\\ Assumption \ref{assump:I4}};

        \draw[arr] (T4.south) -- ++(0,-0.6) -| ([xshift=-1.2cm]T2.north);
        \draw[arr] (T5.south) -- ++(0,-0.6) -| ([xshift=1.2cm]T2.north);
        \draw[arr] ([xshift=1.5cm]T1.south) |- node[lab, pos=0.75, above=2pt] {uniqueness} (T2.west);
        \draw[arr] (T2.east) -- (T3.west);
        \draw[arr] ([xshift=-1.2cm]I0.north) -- ([xshift=-1.2cm]I0.north |- T1.south);
        \draw[arr] (I1.north) -- (I1.north |- T2.south);
        \draw[arr] (I2.north) -- (I2.north |- T2.south);
        \draw[arr] (I3.north) -- (T3.south);
    \end{tikzpicture}%
    }
  \caption{Content of the paper: the two steps of the proof (red), the well-posedness from the literature (orange),
    the main results (green), and the initial assumptions they need in addition to
    \ref{assump:A1}--\ref{assump:A4} (dashed). Each red and green box names where it is proved. TikZ code drafted with the assistance of AI.}
    \label{fig:content_overview}
\end{figure}

The paper is organized as follows, and Figure \ref{fig:content_overview} shows its content. Section \ref{sec:main_results} derives the magnetic hierarchy of Theorem \ref{thm:magnetic_hierarchy_k}. Section \ref{sec:remainder_bounds} reduces the estimate for the remainder to one master estimate, which is verified component by component to finish the estimates. We postpone the details of reducing every remainder to the master form to Section \ref{sec:unified_bounds}, which completes the proof of Theorem \ref{thm:unified_remainders_k}. Section \ref{sec:combined_limit} provides the a priori bounds on the energy and the phase-space moments and combines these results to get the compactness argument for completing the proof of Theorem \ref{thm:magnetic_vlasov_limit}. In Section \ref{sec:wasserstein_noV} we prove Theorem \ref{thm:w2_energy_criterion_noV} by showing the convergence of second moments.

Appendix \ref{app:slater_determinants} proves that the assumptions on the initial data can be achieved by Slater determinants, and Appendix
\ref{app:corollary_proofs} proves Corollaries \ref{cor:wasserstein_p_convergence} and
\ref{cor:quantitative_transport_rate_noV}.

\section{The \texorpdfstring{$k$}{k}-Particle Magnetic Hierarchy for Husimi Measures}
\label{sec:main_results}
This section derives the evolution equation of the $k$-particle magnetic Husimi measure, the magnetic BBGKY
hierarchy in phase space with an explicit quantum remainder.

\begin{proof}[Proof of Theorem \ref{thm:magnetic_hierarchy_k}]
For $X = (x_1, \dots, x_k)$ and $Y = (y_1, \dots, y_k)$ we write
\begin{align}
\label{eq:spatial_operator_notation}
    a_X &\coloneqq a_{x_k} \cdots a_{x_1}, \qquad \dXkYk \coloneqq dx_1 \, dy_1 \cdots dx_k \, dy_k.
\end{align}
We first define the quantum remainder $\mathbf{R}^{(k)}$ of \eqref{eq:magnetic_hierarchy_k} and then show that $\mhus^{(k)}_{N,t}$ satisfies \eqref{eq:magnetic_hierarchy_k} with it. The remainder is
\begin{align}
\label{eq:total_remainder_definitions}
    \mathbf{R}^{(k)}(Z_k) &\coloneqq \sum_{j=1}^k \mathcal{R}_{\text{mag}, j}^{(k)}(Z_k) + \sum_{j=1}^k \mathcal{R}^{(k)}_{\text{ext}, j}(Z_k) + \sum_{j=1}^k \nabla_{v_j} \cdot \mathcal{R}_{\text{mf}, j}^{(k)}(Z_k)
        + \tfrac{1}{2\hbar N} \sum_{1 \le i \neq j \le k} \mathcal{R}_{\text{cor}, i, j}^{(k)}(Z_k),
\end{align}
with the magnetic remainder
\begin{align}
\label{eq:prop_R_mag_def}
    \mathcal{R}_{\text{mag}, j}^{(k)}(Z_k) &\coloneqq \mathcal{R}_{\text{dis}, j}^{(k)}(Z_k) + \mathcal{R}_{\text{for}, j}^{(k)}(Z_k) + \mathcal{R}_{\text{gau}, j}^{(k)}(Z_k) + \mathcal{R}_{\text{dia}, j}^{(k)}(Z_k) + \mathcal{R}_{\text{div}, j}^{(k)}(Z_k)
\end{align}
whose five components are
\begin{align}
\label{eq:prop_R_disp}
    \mathcal{R}_{\text{dis}, j}^{(k)}(Z_k) &\coloneqq \tfrac{c_{N,k}}{(2\pi)^{3k}} \tfrac{1}{\hbar} \operatorname{Im} \left\langle \Psi_{N,t}, \mathfrak{A}_k^\dagger[\mathfrak{f}] \, \mathfrak{A}_{k,j} \big[ (i\hbar \nabla_{q_j}|_{p_j})^2 \mathfrak{f}_j \big] \Psi_{N,t} \right\rangle, \\
\label{eq:prop_R_for}
    \mathcal{R}_{\text{for}, j}^{(k)}(Z_k) &\coloneqq \tfrac{c_{N,k}}{(2\pi)^{3k}} \tfrac{2}{\hbar} \operatorname{Im} \Big\langle \Psi_{N,t}, \mathfrak{A}_k^\dagger[\mathfrak{f}] \mathfrak{A}_{k,j} \Big[ \Big( ( \mathrm M_A(x_j,q_j))^T v_j \Big) \cdot (x_j-q_j) \mathfrak{f}_j \Big] \Psi_{N,t} \Big\rangle, \\
\label{eq:prop_R_gau}
    \mathcal{R}_{\text{gau}, j}^{(k)}(Z_k) &\coloneqq \tfrac{c_{N,k}}{(2\pi)^{3k}} \operatorname{Re} \Big\langle \Psi_{N,t}, \mathfrak{A}_k^\dagger[\mathfrak{f}] \mathfrak{A}_{k,j} \Big[ \big( (\nabla \cdot A)(x_j) - (\nabla \cdot A)(q_j) \big) \mathfrak{f}_j \Big] \Psi_{N,t} \Big\rangle, \\
\label{eq:prop_R_dia}
    \mathcal{R}_{\text{dia}, j}^{(k)}(Z_k) &\coloneqq \tfrac{c_{N,k}}{(2\pi)^{3k}} \tfrac{1}{\hbar} \operatorname{Im} \left\langle \Psi_{N,t}, \mathfrak{A}_k^\dagger[\mathfrak{f}] \, \mathfrak{A}_{k,j} \Big[ |A(x_j)-A(q_j)|^2 \mathfrak{f}_j \Big] \Psi_{N,t} \right\rangle, \\
\label{eq:prop_R_FS}
    \mathcal{R}_{\text{div}, j}^{(k)}(Z_k) &\coloneqq \tfrac{c_{N,k}}{(2\pi)^{3k}} \operatorname{Re} \Big\langle \Psi_{N,t}, \mathfrak{A}_k^\dagger[\mathfrak{f}] \mathfrak{A}_{k,j} \Big[ \big( \nabla \cdot A(q_j)- 2 \big(A(x_j) - A(q_j)\big) \cdot \nabla_{q_j}\vert_{p_j} \big) \mathfrak{f}_j \Big] \Psi_{N,t} \Big\rangle
\end{align}
with the averaged deviation of $\nabla A$ on the segment from $q_j$ to $x_j$,
\begin{align} \label{eq:MAaveragedgradient_k}
 \mathrm M_A(x_j,q_j) &\coloneqq \int_0^1 ds \, \nabla A\big(q_j + s(x_j-q_j)\big) - \nabla A(q_j).
\end{align}
The mean-field and internal-correlation remainders are
\begin{align}
\label{eq:R_mf_kin_def_prop}
    \mathcal{R}_{\text{mf}, j}^{(k)}(Z_k) &\coloneqq \tfrac{c_{N,k}}{(2\pi)^{3(k+1)}} \operatorname{Re} \int_0^1 ds \iint_{\mathbb{R}^6} dq_{k+1} \, dv_{k+1}
        \int_{\mathbb{R}^{6(k+1)}} \dXkYkplus \left( \prod\nolimits_{m \neq j}^k \mathfrak{f}_m(x_m) \overline{\mathfrak{f}_m(y_m)} \right) \nonumber \\
    &\qquad \times \mathfrak{f}_j(x_j) \overline{\mathfrak{f}_j(y_j)} \mathfrak{f}_{k+1}(x_{k+1}) \overline{\mathfrak{f}_{k+1}(y_{k+1})}
        \Big[ \nabla W\big(s y_j + (1-s) x_j - x_{k+1}\big) - \nabla W(q_j - q_{k+1}) \Big] \nonumber \\
    &\qquad \times \left\langle \Psi_{N,t}, a_X^\dagger a_{x_{k+1}}^\dagger a_{y_{k+1}} a_Y \Psi_{N,t} \right\rangle, \\
\label{eq:R_intrn_kin_def_prop}
    \mathcal{R}_{\text{cor}, i, j}^{(k)}(Z_k) &\coloneqq \tfrac{c_{N,k}}{(2\pi)^{3k}} \operatorname{Im} \int_{\mathbb{R}^{6k}} \dXkYk \left( \prod\nolimits_{m \neq i,j}^k \mathfrak{f}_m(x_m) \overline{\mathfrak{f}_m(y_m)} \right)
        \mathfrak{f}_i(x_i) \overline{\mathfrak{f}_i(y_i)} \mathfrak{f}_j(x_j) \overline{\mathfrak{f}_j(y_j)} \nonumber \\
    &\qquad \times \Big( W(y_i-y_j) - W(x_i-x_j) \Big) \left\langle \Psi_{N,t}, a_X^\dagger a_Y \Psi_{N,t} \right\rangle
\end{align}
and the external remainder is
\begin{align}
\label{eq:rkinextdef}
    \mathcal{R}^{(k)}_{\text{ext}, j}(Z_k) &\coloneqq \tfrac{c_{N,k}}{(2\pi)^{3k}} \tfrac{2}{\hbar} \operatorname{Im} \Big\langle \Psi_{N,t}, \mathfrak{A}_k^\dagger[\mathfrak{f}] \mathfrak{A}_{k, j}\Big[ \Big(\Vext(x_j) - (x_j-q_j)\cdot \nabla \Vext(q_j)\Big) \mathfrak{f}_j \Big] \Psi_{N,t} \Big\rangle.
\end{align}

\paragraph{Derivation of the hierarchy.} We show that $\mhus^{(k)}_{N,t}$ satisfies \eqref{eq:magnetic_hierarchy_k} with $\mathbf{R}^{(k)}$ of \eqref{eq:total_remainder_definitions}. Throughout this proof, the dependence of the operators on $Z_k$ is not displayed. Let $\Psi_{N,t}$ be the solution to \eqref{eq:Schrodinger_scaled}. Decomposing the Hamiltonian $\mathcal{H} = \mathcal{K}_A + \mathcal\Vext + \mathcal{W}$, the time evolution of the $k$-particle magnetic Husimi measure yields the commutator sum:
\begin{align}
\label{eq:hierarchy_commutator_split}
    \partial_t \mhus_{N,t}^{(k)} &= \tfrac{c_{N,k}}{(2\pi)^{3k}} \tfrac{1}{i\hbar} \left\langle \Psi_{N,t}, \left[\mathfrak{A}_k^\dagger[\mathfrak{f}] \, \mathfrak{A}_k[\mathfrak{f}], \mathcal{H}\right] \Psi_{N,t} \right\rangle
    \eqqcolon \mathcal{T}_{\text{kin}}^{(k)} + \mathcal{T}_{\text{ext}}^{(k)} + \mathcal{T}_{\text{int}}^{(k)}.
\end{align}
\paragraph{The magnetic kinetic term \texorpdfstring{$\mathcal{T}_{\text{kin}}^{(k)}$}{T\_kin\textasciicircum{(k)}}.}
A direct computation gives
\begin{align}
\label{eq:single_particle_commutator}
    [a(\mathfrak{f}_j), \mathcal{K}_A] &= \tfrac{1}{2} a\left( \mpi_{x_j}^2 \mathfrak{f}_j \right), \qquad [a^\dagger(\mathfrak{f}_j), \mathcal{K}_A] = -\tfrac{1}{2} a^\dagger\left( \mpi_{x_j}^2 \mathfrak{f}_j \right), \\
\label{eq:annihilation_string_commutator}
    \left[ \mathfrak{A}_k[\mathfrak{f}], \mathcal{K}_A \right] &= \sum_{j=1}^k \left( \prod\nolimits_{m=k}^{j+1} a(\mathfrak{f}_m) \right) \left[ a(\mathfrak{f}_j), \mathcal{K}_A \right] \left( \prod\nolimits_{m=j-1}^1 a(\mathfrak{f}_m) \right) = \tfrac{1}{2} \sum_{j=1}^k \mathfrak{A}_{k,j}\left[ \mpi_{x_j}^2 \mathfrak{f}_j \right].
\end{align}
The self-adjointness of $\mathcal{K}_A$ reduces the imaginary component within the expectation value:
\begin{align}
\label{eq:tilde_T_kin_exact_k}
    \mathcal{T}_{\text{kin}}^{(k)} &\coloneqq \tfrac{c_{N,k}}{(2\pi)^{3k}} \tfrac{1}{i\hbar} \left\langle \Psi_{N,t}, \left[\mathfrak{A}_k^\dagger[\mathfrak{f}] \, \mathfrak{A}_k[\mathfrak{f}], \mathcal{K}_A\right] \Psi_{N,t} \right\rangle = \tfrac{c_{N,k}}{(2\pi)^{3k}} \tfrac{1}{\hbar} \sum_{j=1}^k \operatorname{Im} \left\langle \Psi_{N,t}, \mathfrak{A}_k^\dagger[\mathfrak{f}] \, \mathfrak{A}_{k,j}\left[ \mpi_{x_j}^2 \mathfrak{f}_j \right] \Psi_{N,t} \right\rangle.
\end{align}

With the definition \eqref{eq:coherent_shifted} of $\mathfrak{f}_j(x_j)$ we have:
\begin{align}
\label{eq:grad_relation_derived}
    \nabla_{x_j} \mathfrak{f}_j(x_j) &= - \nabla_{q_j}\Big|_{p_j} \mathfrak{f}_j(x_j) + \tfrac{i (v_j-A(q_j))}{\hbar} \mathfrak{f}_j(x_j).
\end{align}

Multiplying by $-i\hbar$ and adding $A(x_j)\mathfrak{f}_j(x_j)$ gives
\begin{align}
\label{eq:pi_linear_action_derived}
    \mpi_{x_j} \mathfrak{f}_j(x_j) &= \left( v_j + i\hbar \nabla_{q_j}\Big|_{p_j} + A(x_j) - A(q_j) \right) \mathfrak{f}_j(x_j) \\
\label{eq:pi_squared_initial_derived}
    \mpi_{x_j}^2 \mathfrak{f}_j(x_j) &= \mpi_{x_j} \cdot \left[ \left( v_j + i\hbar \nabla_{q_j}\Big|_{p_j} + A(x_j) - A(q_j) \right) \mathfrak{f}_j(x_j) \right].
\end{align}

The commutator $[\mpi_{x_j}, A(x_j)-A(q_j)] = -i\hbar (\nabla \cdot A)(x_j)$ yields:
\begin{align*}
    \mpi_{x_j} \cdot \left( v_j \mathfrak{f}_j(x_j) \right) &= \left[ v_j^2 + v_j \cdot \left( i\hbar \nabla_{q_j}\Big|_{p_j} \right)
        + v_j \cdot (A(x_j) - A(q_j)) \right] \mathfrak{f}_j(x_j), \\
    \mpi_{x_j} \cdot \Big( i\hbar \nabla_{q_j}\Big|_{p_j} \mathfrak{f}_j(x_j) \Big) &= \left[ v_j \cdot \Big( i\hbar \nabla_{q_j}\Big|_{p_j} \Big) + \Big( i\hbar \nabla_{q_j}\Big|_{p_j} \Big)^2
        + (A(x_j) - A(q_j)) \cdot \Big( i\hbar \nabla_{q_j}\Big|_{p_j} \Big) \right] \mathfrak{f}_j(x_j), \\
    \mpi_{x_j} \cdot \left( (A(x_j) - A(q_j)) \mathfrak{f}_j(x_j) \right) &= \Bigg[ -i\hbar (\nabla \cdot A)(x_j) + (A(x_j) - A(q_j)) \cdot v_j
        + (A(x_j) - A(q_j)) \cdot \left( i\hbar \nabla_{q_j}\Big|_{p_j} \right)  \\
    &\qquad + |A(x_j) - A(q_j)|^2 \Bigg] \mathfrak{f}_j(x_j).
\end{align*}

We define the seven operator terms:
\begin{align}
\label{eq:operator_terms_derived}
    \mathfrak{s}_1 &\coloneqq v_j^2 \mathfrak{f}_j, & \mathfrak{s}_2 &\coloneqq 2 v_j \cdot (i\hbar \nabla_{q_j}|_{p_j}) \mathfrak{f}_j, \nonumber \\
    \mathfrak{s}_3 &\coloneqq (i\hbar \nabla_{q_j}|_{p_j})^2 \mathfrak{f}_j, & \mathfrak{s}_4 &\coloneqq 2 ((A(x_j)-A(q_j)) \cdot v_j) \mathfrak{f}_j, \nonumber \\
    \mathfrak{s}_5 &\coloneqq 2 (A(x_j)-A(q_j)) \cdot (i\hbar \nabla_{q_j}|_{p_j}) \mathfrak{f}_j, & \mathfrak{s}_6 &\coloneqq |A(x_j)-A(q_j)|^2 \mathfrak{f}_j, \nonumber \\
    \mathfrak{s}_7 &\coloneqq -i\hbar (\nabla \cdot A)(x_j) \mathfrak{f}_j.
\end{align}
and this splits the kinetic term \eqref{eq:tilde_T_kin_exact_k} into a sum:
\begin{align}
\label{eq:recombination_derived}
    \mpi_{x_j}^2 \mathfrak{f}_j &= \sum_{m=1}^7 \mathfrak{s}_m, \quad \implies \quad \mathcal{T}_{\text{kin}}^{(k)} = \sum_{m=1}^7 \mathfrak{S}_m, \nonumber \\
    \mathfrak{S}_m &\coloneqq \tfrac{c_{N,k}}{(2\pi)^{3k}} \tfrac{1}{\hbar} \sum_{j=1}^k \operatorname{Im} \big\langle \Psi_{N,t}, \mathfrak{A}_k^\dagger[\mathfrak{f}] \, \mathfrak{A}_{k,j}[ \mathfrak{s}_m ] \Psi_{N,t} \big\rangle.
\end{align}

Three of the seven terms follow immediately:
\begin{align}
\label{eq:S_expectations_new}
    \mathfrak{S}_1 =   0, \qquad \qquad
    \mathfrak{S}_3 =  \sum_{j=1}^k \mathcal{R}_{\text{dis}, j}^{(k)}, \qquad \qquad
    \mathfrak{S}_6 = \sum_{j=1}^k\mathcal{R}_{\text{dia}, j}^{(k)}.
\end{align}

Since $\mhus_{N,t}^{(k)}$ is real, we have:
\begin{align} \label{eq:husimifromre}
   \nabla_{v_j} \mhus_{N,t}^{(k)} &= \tfrac{c_{N,k}}{(2\pi)^{3k}} \, 2 \operatorname{Re} \langle \Psi_{N,t}, \mathfrak{A}_k^\dagger[\mathfrak f] \mathfrak{A}_{k,j}[ \nabla_{v_j} \mathfrak{f}_j ] \Psi_{N,t} \rangle.
\end{align}

By \eqref{eq:state_level_chain_rule}, $\mathfrak{S}_2$ reads:
\begin{align}
\label{eq:T2_def_combined_derived_new}
    \mathfrak{S}_2 &= \tfrac{c_{N,k}}{(2\pi)^{3k}} \tfrac{1}{\hbar} \sum_{j=1}^k \operatorname{Im} \big\langle \Psi_{N,t}, \mathfrak{A}_k^\dagger[\mathfrak{f}] \, \mathfrak{A}_{k,j}\left[ 2 v_j \cdot (i\hbar \nabla_{q_j}|_{p_j}) \mathfrak{f}_j \right] \Psi_{N,t} \big\rangle \nonumber \\
    &= - \tfrac{c_{N,k}}{(2\pi)^{3k}} \, 2 \sum_{j=1}^k \operatorname{Re} \left\langle \Psi_{N,t}, \mathfrak{A}_k^\dagger[\mathfrak{f}] \, \mathfrak{A}_{k,j} \left[ v_j \cdot \left( \nabla_{q_j} + (\nabla A(q_j))^T \nabla_{v_j} \right) \mathfrak{f}_j \right] \Psi_{N,t} \right\rangle \nonumber \\
    &= - \sum_{j=1}^k v_j \cdot \nabla_{q_j} \mhus_{N,t}^{(k)} - \sum_{j=1}^k \Big( (\nabla A(q_j)) v_j \Big) \cdot \nabla_{v_j} \mhus_{N,t}^{(k)}.
\end{align}

Definition \eqref{eq:MAaveragedgradient_k} reformulates the spatial difference:
\begin{align}
\label{eq:spatial_difference_reformulated}
     A(x_j) - A(q_j) &= \nabla A(q_j) (x_j-q_j) +  \mathrm M_A(x_j, q_j)(x_j - q_j).
\end{align}

Applying $(x_j-q_j) \mathfrak{f}_j = (-i\hbar \nabla_{v_j} - q_j) \mathfrak{f}_j$ (see \eqref{eq:state_level_chain_rule}) to the linear interaction expands the phase-space coupling:
\begin{align}
\label{eq:T4_expanded_and_split_k_derived_new}
    \mathfrak{S}_4 &= \tfrac{c_{N,k}}{(2\pi)^{3k}} \tfrac{2}{\hbar} \sum_{j=1}^k \operatorname{Im} \Big\langle \Psi_{N,t},
        \mathfrak{A}_k^\dagger[\mathfrak{f}] \, \mathfrak{A}_{k,j} \Big[ \Big( \nabla A(q_j)(x_j - q_j) + \mathrm M_A(x_j, q_j)(x_j - q_j) \Big) \cdot v_j \mathfrak{f}_j \Big] \Psi_{N,t} \Big\rangle \nonumber \\
    &= \tfrac{c_{N,k}}{(2\pi)^{3k}} \tfrac{2}{\hbar} \sum_{j=1}^k \operatorname{Im} \Big\langle \Psi_{N,t}, \mathfrak{A}_k^\dagger[\mathfrak{f}] \, \mathfrak{A}_{k,j} \Big[ \Big( (\nabla A(q_j))^T v_j \Big) \cdot (-i\hbar \nabla_{v_j} - q_j) \mathfrak{f}_j \Big] \Psi_{N,t} \Big\rangle \nonumber \\
    &\quad + \tfrac{c_{N,k}}{(2\pi)^{3k}} \tfrac{2}{\hbar} \sum_{j=1}^k \operatorname{Im} \Big\langle \Psi_{N,t}, \mathfrak{A}_k^\dagger[\mathfrak{f}] \, \mathfrak{A}_{k,j} \Big[ \Big( ( \mathrm M_A(x_j,q_j))^T v_j \Big) \cdot (x_j-q_j) \mathfrak{f}_j \Big] \Psi_{N,t} \Big\rangle \nonumber \\
    &= \sum_{j=1}^k \Big( (\nabla A(q_j))^T v_j \Big) \cdot \nabla_{v_j} \mhus_{N,t}^{(k)} + \sum_{j=1}^k \mathcal{R}_{\text{for}, j}^{(k)}.
\end{align}

Adding and subtracting $\nabla \cdot A (q_j)$ and by the definition of $\mhus_{N,t}^{(k)}$:
\begin{align}
\label{eq:T7_Rem_derived}
    \mathfrak{S}_7 &= \sum_{j=1}^k (\nabla \cdot A)(q_j) \mhus_{N,t}^{(k)} + \sum_{j=1}^k \mathcal{R}_{\text{gau}, j}^{(k)}.
\end{align}

A direct computation gives
\begin{align}\label{eq:mag_vector_identity}
(\nabla A(q_j))^T v_j - (\nabla A(q_j)) v_j & = v_j \times (\nabla \times A(q_j)) = v_j \times B(q_j).
\end{align}

Grouping the terms of $\mathcal{T}_{\text{kin}}^{(k)}$ and applying the vector identity \eqref{eq:mag_vector_identity} yields the classical magnetic advection. The remaining terms form the remainder $\mathcal{R}_{\text{mag}, j}^{(k)}$ defined in \eqref{eq:prop_R_mag_def}:
\begin{align}
\label{eq:T_kin_total_synthesis}
    \mathcal{T}_{\text{kin}}^{(k)} &= \sum_{m=1}^7 \mathfrak{S}_m = \mathfrak{S}_2 + \sum_{j=1}^k \mathcal{R}_{\text{dis}, j}^{(k)} + \mathfrak{S}_4 + \mathfrak{S}_5 + \sum_{j=1}^k \mathcal{R}_{\text{dia}, j}^{(k)} + \mathfrak{S}_7 \nonumber \\
    &= - \sum_{j=1}^k \left( v_j \cdot \nabla_{q_j} \mhus_{N,t}^{(k)} + \Big( (\nabla A(q_j)) v_j \Big) \cdot \nabla_{v_j} \mhus_{N,t}^{(k)} \right) + \sum_{j=1}^k \mathcal{R}_{\text{dis}, j}^{(k)}
        + \sum_{j=1}^k \Big( (\nabla A(q_j))^T v_j \Big) \cdot \nabla_{v_j} \mhus_{N,t}^{(k)} \nonumber \\
    &\quad + \sum_{j=1}^k \mathcal{R}_{\text{for}, j}^{(k)}
        - \tfrac{c_{N,k}}{(2\pi)^{3k}} \, 2 \sum_{j=1}^k \operatorname{Re} \left\langle \Psi_{N,t}, \mathfrak{A}_k^\dagger[\mathfrak{f}] \, \mathfrak{A}_{k,j} \left[ (A(x_j)- A(q_j)) \cdot \nabla_{q_j}|_{p_j} \mathfrak{f}_j \right] \Psi_{N,t} \right\rangle \nonumber \\
    &\quad + \sum_{j=1}^k \mathcal{R}_{\text{dia}, j}^{(k)} + \sum_{j=1}^k (\nabla \cdot A)(q_j) \mhus_{N,t}^{(k)} + \sum_{j=1}^k \mathcal{R}_{\text{gau}, j}^{(k)} \nonumber \\
    &= - \sum_{j=1}^k v_j \cdot \nabla_{q_j} \mhus_{N,t}^{(k)} + \sum_{j=1}^k \Big( (\nabla A(q_j))^T v_j - (\nabla A(q_j)) v_j \Big) \cdot \nabla_{v_j} \mhus_{N,t}^{(k)}
        + \sum_{j=1}^k (\nabla \cdot A)(q_j) \mhus_{N,t}^{(k)} \nonumber \\
    &\quad - \tfrac{c_{N,k}}{(2\pi)^{3k}} \, 2 \sum_{j=1}^k \operatorname{Re} \Big\langle \Psi_{N,t}, \mathfrak{A}_k^\dagger[\mathfrak{f}] \, \mathfrak{A}_{k,j} \Big[ (A(x_j)- A(q_j)) \cdot \nabla_{q_j}|_{p_j} \mathfrak{f}_j \Big] \Psi_{N,t} \Big\rangle \nonumber \\
    &\quad + \sum_{j=1}^k \Big( \mathcal{R}_{\text{dis}, j}^{(k)} + \mathcal{R}_{\text{for}, j}^{(k)} + \mathcal{R}_{\text{dia}, j}^{(k)} + \mathcal{R}_{\text{gau}, j}^{(k)} \Big) \nonumber \\
    &= - \sum_{j=1}^k \left( v_j \cdot \nabla_{q_j} - (v_j \times B(q_j)) \cdot \nabla_{v_j} \right) \mhus_{N,t}^{(k)} \nonumber \\
    &\quad + \sum_{j=1}^k \Big( \mathcal{R}_{\text{dis}, j}^{(k)} + \mathcal{R}_{\text{for}, j}^{(k)} + \mathcal{R}_{\text{div}, j}^{(k)} + \mathcal{R}_{\text{dia}, j}^{(k)} + \mathcal{R}_{\text{gau}, j}^{(k)} \Big).
\end{align}

\paragraph{The interaction term \texorpdfstring{$\mathcal{T}_{\text{int}}^{(k)}$}{(T\_\{int\}\textasciicircum\{(k)\})}.} \label{subsec:interaction_term}
The interaction term $\mathcal{T}_{\text{int}}^{(k)}$ is generated by the commutator with $\mathcal{W}$. Computation via the canonical anticommutation relations and the derivation property of the commutator yields:
\begin{align}
\label{eq:app_leibniz_expansion}
    [a_{y_k} \dots a_{y_1}, \mathcal{W}] &= \tfrac{1}{N} \sum_{j=1}^k \int_{\mathbb{R}^3} dy_{k+1} \, W(y_j-y_{k+1}) \Big( a_{y_k} \dots a_{y_{j+1}} a_{y_{k+1}}^\dagger a_{y_{k+1}} a_{y_j} \dots a_{y_1} \Big).
\end{align}

Normal ordering $a_{y_{k+1}}^\dagger a_{y_{k+1}}$ to the left of the $k-j$ annihilation operators separates the sequence into a mean-field background coupling and internal correlation pairing $\{a_{y_i}, a_{y_{k+1}}^\dagger\} = \delta(y_i - y_{k+1})$:
\begin{align*}
    a_{y_k} \dots a_{y_{j+1}} a_{y_{k+1}}^\dagger a_{y_{k+1}} &= \underbrace{a_{y_{k+1}}^\dagger a_{y_{k+1}} a_{y_k} \dots a_{y_{j+1}}}_{\text{Background}} + \underbrace{\sum_{i=j+1}^k (-1)^{i-j-1} \delta(y_i - y_{k+1}) a_{y_k} \dots \widehat{a_{y_i}} \dots a_{y_{j+1}} a_{y_{k+1}}}_{\text{Internal Correlation}},
\end{align*}
where $\widehat{a_{y_i}}$ denotes operator omission. The background component covers the $(k+1)$-th marginal by inserting the resolution of the identity \eqref{eq:resolution_identity}. For the internal term, the spatial resolution $y_{k+1} = y_i$:

\begin{align}
\label{eq:app_bg_string_eval}
    \Big( \mathfrak{A}_k^\dagger[\mathfrak{f}] \Big[\mathfrak{A}_k[\mathfrak{f}], \mathcal{W}\Big] \Big)_{\text{bg}, j} &= \tfrac{1}{N} \tfrac{1}{(2\pi\hbar)^3} \iint_{\mathbb{R}^6} dq_{k+1} \, dv_{k+1} \int_{\mathbb{R}^{6(k+1)}} \dXkYkplus \left( \prod\nolimits_{m=1}^{k+1} \mathfrak{f}_m(x_m) \overline{\mathfrak{f}_m(y_m)} \right) \\
    &\qquad \qquad \qquad \qquad \qquad \qquad \qquad \qquad \qquad \times W(x_j-x_{k+1}) a_X^\dagger a_{x_{k+1}}^\dagger a_{y_{k+1}} a_Y, \nonumber \\[8pt]
\label{eq:app_internal_string_eval}
    \Big( \mathfrak{A}_k^\dagger[\mathfrak{f}] \Big[\mathfrak{A}_k[\mathfrak{f}], \mathcal{W}\Big] \Big)_{\text{int}, j} &= \tfrac{1}{N} \sum_{i=j+1}^k \int_{\mathbb{R}^{6k}} \dXkYk \left( \prod\nolimits_{m=1}^k \mathfrak{f}_m(x_m) \overline{\mathfrak{f}_m(y_m)} \right) W(x_j-x_i) a_X^\dagger a_Y.
\end{align}
Swapping the integration variables $x \leftrightarrow y$ yields the conjugate terms. Symmetrizing under $i \leftrightarrow j$ turns the triangular sum to $\frac{1}{2}\sum_{i \neq j}$:
{
\begin{align}
\label{eq:Gentbgtintrn}
    & [\mathfrak{A}_k^\dagger[\mathfrak{f}] \mathfrak{A}_k[\mathfrak{f}], \mathcal{W}] = \mathfrak{A}_k^\dagger[\mathfrak{f}] [\mathfrak{A}_k[\mathfrak{f}], \mathcal{W}] - \big( \mathfrak{A}_k^\dagger[\mathfrak{f}] [\mathfrak{A}_k[\mathfrak{f}], \mathcal{W}] \big)^\dagger \nonumber \\
    &= \Bigg[ \tfrac{1}{N} \sum_{j=1}^k \tfrac{1}{(2\pi\hbar)^3} \iint_{\mathbb{R}^6} dq_{k+1} \, dv_{k+1} \int_{\mathbb{R}^{6(k+1)}} \dXkYkplus \nonumber \\
    &\qquad \qquad \times \left( \prod\nolimits_{m=1}^{k+1} \mathfrak{f}_m(x_m) \overline{\mathfrak{f}_m(y_m)} \right) W(x_j-x_{k+1}) a_X^\dagger a_{x_{k+1}}^\dagger a_{y_{k+1}} a_Y \nonumber \\
    &\qquad + \tfrac{1}{2N} \sum_{1 \le i \neq j \le k} \int_{\mathbb{R}^{6k}} \dXkYk \left( \prod\nolimits_{m=1}^k \mathfrak{f}_m(x_m) \overline{\mathfrak{f}_m(y_m)} \right) W(x_j-x_i) a_X^\dagger a_Y \Bigg] \nonumber \\
    &\quad - \Bigg[ \tfrac{1}{N} \sum_{j=1}^k \tfrac{1}{(2\pi\hbar)^3} \iint_{\mathbb{R}^6} dq_{k+1} \, dv_{k+1} \int_{\mathbb{R}^{6(k+1)}} \dXkYkplus \nonumber \\
    &\qquad \qquad \times \left( \prod\nolimits_{m=1}^{k+1} \mathfrak{f}_m(x_m) \overline{\mathfrak{f}_m(y_m)} \right) W(y_j-x_{k+1}) a_X^\dagger a_{x_{k+1}}^\dagger a_{y_{k+1}} a_Y \nonumber \\
    &\qquad + \tfrac{1}{2N} \sum_{1 \le i \neq j \le k} \int_{\mathbb{R}^{6k}} \dXkYk \left( \prod\nolimits_{m=1}^k \mathfrak{f}_m(x_m) \overline{\mathfrak{f}_m(y_m)} \right) W(y_j-y_i) a_X^\dagger a_Y \Bigg] \nonumber \\
    &= \tfrac{1}{N} \sum_{j=1}^k \tfrac{1}{(2\pi\hbar)^3} \iint_{\mathbb{R}^6} dq_{k+1} \, dv_{k+1}
        \int_{\mathbb{R}^{6(k+1)}} \dXkYkplus \left( \prod\nolimits_{m=1}^{k+1} \mathfrak{f}_m(x_m) \overline{\mathfrak{f}_m(y_m)} \right) \nonumber \\
    &\qquad \times \Big( W(x_j-x_{k+1}) - W(y_j-x_{k+1}) \Big) a_X^\dagger a_{x_{k+1}}^\dagger a_{y_{k+1}} a_Y \nonumber \\
    &\quad + \tfrac{1}{2N} \sum_{1 \le i \neq j \le k} \int_{\mathbb{R}^{6k}} \dXkYk \left( \prod\nolimits_{m=1}^k \mathfrak{f}_m(x_m) \overline{\mathfrak{f}_m(y_m)} \right)
        \Big( W(x_j-x_i) - W(y_j-y_i) \Big) a_X^\dagger a_Y.
\end{align}
}
Applying \eqref{eq:Gentbgtintrn}, the commutator distributes across the $k$ particles, separating into a coupling to the $(k+1)$-th marginal $\mathcal{T}_{\text{mf}}^{(k)}$ and an internal correlation $\mathcal{T}_{\text{cor}}^{(k)}$:
\begin{align}
\label{eq:interaction_comm_k}
    \mathcal{T}_{\text{int}}^{(k)} &\coloneqq \tfrac{c_{N,k}}{(2\pi)^{3k}} \tfrac{1}{i\hbar} \left\langle \Psi_{N,t}, \left[ \mathfrak{A}_k^\dagger[\mathfrak{f}] \mathfrak{A}_k[\mathfrak{f}], \mathcal{W} \right] \Psi_{N,t} \right\rangle = \mathcal{T}_{\text{mf}}^{(k)} + \mathcal{T}_{\text{cor}}^{(k)}.
\end{align}
Expanding the operators $a(\mathfrak{f}_j)$ into the continuous spatial basis and inserting the resolution of the identity \eqref{eq:resolution_identity} for the $(k+1)$-th particle, the respective terms take the form of spatial field strings $a_X^\dagger$ and $a_Y$:
\begin{align}
\label{eq:T_bg_exact}
    \mathcal{T}_{\text{mf}}^{(k)} &= - \tfrac{c_{N,k}}{(2\pi)^{3k}} \tfrac{1}{i\hbar N} \sum_{j=1}^k \tfrac{1}{(2\pi\hbar)^3} \iint_{\mathbb{R}^6} dq_{k+1} \, dv_{k+1}
        \int_{\mathbb{R}^{6(k+1)}} \dXkYkplus \left( \prod\nolimits_{m=1}^{k+1} \mathfrak{f}_m(x_m) \overline{\mathfrak{f}_m(y_m)} \right) \nonumber \\
    &\qquad \times \Big( W(x_j-x_{k+1}) - W(y_j-y_{k+1}) \Big) \left\langle \Psi_{N,t}, a_X^\dagger a_{x_{k+1}}^\dagger a_{y_{k+1}} a_Y \Psi_{N,t} \right\rangle, \\
\label{eq:T_intrn_exact}
    \mathcal{T}_{\text{cor}}^{(k)} &= - \tfrac{c_{N,k}}{(2\pi)^{3k}} \tfrac{1}{2 i\hbar N} \sum_{1 \le i \neq j \le k} \int_{\mathbb{R}^{6k}} \dXkYk \left( \prod\nolimits_{m=1}^k \mathfrak{f}_m(x_m) \overline{\mathfrak{f}_m(y_m)} \right) \nonumber \\
    &\qquad \times \Big( W(x_i-x_j) - W(y_i-y_j) \Big) \left\langle \Psi_{N,t}, a_X^\dagger a_Y \Psi_{N,t} \right\rangle.
\end{align}
To separate the classical target interaction term from the background coupling $\mathcal{T}_{\text{mf}}^{(k)}$, express the potential difference via the Fundamental Theorem of Calculus and apply \eqref{eq:state_level_chain_rule}:
\begin{align}
\label{eq:Fundamental Theorem of Calculus_background}
    W(w_j-y) - W(u_j-y) &= \int_0^1 ds \, \nabla W\big(su_j + (1-s)w_j - y\big) \cdot (w_j-u_j), \\
\label{eq:kinetic_grad_identity_k}
    (w_j-u_j) \mathfrak{f}_j(w_j) \overline{\mathfrak{f}_j(u_j)} &= -i\hbar \nabla_{v_j} \left( \mathfrak{f}_j(w_j) \overline{\mathfrak{f}_j(u_j)} \right).
\end{align}

Substituting \eqref{eq:kinetic_grad_identity_k} into \eqref{eq:T_bg_exact}, inserting the identity \eqref{eq:resolution_identity}, adding and subtracting $\nabla W(q_j - q_{k+1})$ reconstructs the $(k+1)$-particle magnetic Husimi measure $\mhus_{N,t}^{(k+1)}$ \cite{ChenLeeLiew2021}. The reconstruction carries the ratio $c_{N,k} / c_{N,k+1} = (N-k)/N$ of \eqref{eq:kinetic_husimi_M_def}, the prefactor of \eqref{eq:mean_field_force_def_k}. The Hermiticity of the interaction commutator, combined with the respective symmetry and antisymmetry of the potential differences under coordinate exchange $x \leftrightarrow y$, constrains the background coupling $\mathcal{T}_{\text{mf}}^{(k)}$ to the real axis and the internal correlation $\mathcal{T}_{\text{cor}}^{(k)}$ to the imaginary axis, yielding the $\operatorname{Re}(\cdot)$ and $\operatorname{Im}(\cdot)$ projections for $\mathcal{R}_{\text{mf}, j}^{(k)}$ \eqref{eq:R_mf_kin_def_prop} and $\mathcal{R}_{\text{cor}, i, j}^{(k)}$ \eqref{eq:R_intrn_kin_def_prop}.
The background and internal couplings admit the decompositions:
\begin{align}
\label{eq:T_intrn_sum_ij}
    \mathcal{T}_{\text{mf}}^{(k)} &= \sum_{j=1}^k \nabla_{v_j} \cdot \mathcal{F}_{\text{int}, j}[\mhus_{N,t}^{(k+1)}] + \sum_{j=1}^k \nabla_{v_j} \cdot \mathcal{R}_{\text{mf}, j}^{(k)}, \qquad
    \mathcal{T}_{\text{cor}}^{(k)} = \tfrac{1}{2\hbar N} \sum_{1 \le i \neq j \le k} \mathcal{R}_{\text{cor}, i, j}^{(k)} .
\end{align}
The interaction term decomposition yields:
\begin{align}
    \mathcal{T}_{\text{int}}^{(k)} &=\sum_{j=1}^k \nabla_{v_j} \cdot \mathcal{F}_{\text{int}, j}[\mhus_{N,t}^{(k+1)}] + \sum_{j=1}^k \nabla_{v_j} \cdot \mathcal{R}_{\text{mf}, j}^{(k)} + \tfrac{1}{2\hbar N} \sum_{1 \le i \neq j \le k} \mathcal{R}_{\text{cor}, i, j}^{(k)} .\label{eq:Tintfinal}
\end{align}

\paragraph{The external term \texorpdfstring{$\mathcal{T}_{\text{ext}}^{(k)}$}{T-ext(k)}.}

$\mathcal{T}_{\text{ext}}^{(k)}$ is governed by $\mathcal\Vext$ and its commutator with the $k$-particle annihilation string expands as:
\begin{align}
\label{eq:V_ext_commutator}
    [\mathfrak{A}_k[\mathfrak{f}], \mathcal\Vext]   &= \int_{\mathbb{R}^3} dx \, \Vext(x) [\mathfrak{A}_k[\mathfrak{f}],a_x^\dagger a_x ] = \sum_{j=1}^k \mathfrak{A}_{k,j}[ \Vext \mathfrak{f}_j ].
\end{align}

The derivation property of the commutator produces the imaginary component for each particle $j$:
\begin{align}
\label{eq:Text_exact_Im_k}
    \mathcal{T}_{\text{ext}}^{(k)} &\coloneqq \tfrac{c_{N,k}}{(2\pi)^{3k}} \tfrac{1}{i\hbar} \big\langle \Psi_{N,t}, [\mathfrak{A}_k^\dagger[\mathfrak{f}] \mathfrak{A}_k[\mathfrak{f}], \mathcal\Vext] \Psi_{N,t} \big\rangle = \sum_{j=1}^k \tfrac{c_{N,k}}{(2\pi)^{3k}} \tfrac{2}{\hbar} \operatorname{Im} \big\langle \Psi_{N,t}, \mathfrak{A}_k^\dagger[\mathfrak{f}] \mathfrak{A}_{k,j}\left[ \Vext \mathfrak{f}_j \right] \Psi_{N,t} \big\rangle.
\end{align}
Adding and subtracting the classical force term:
\begin{align}
\label{eq:Text_expansion}
    \Vext(x_j) &= \nabla \Vext(q_j) \cdot (x_j - q_j) + \left( \Vext(x_j) - \nabla \Vext(q_j) \cdot (x_j - q_j) \right).
\end{align}
Inserted into \eqref{eq:Text_exact_Im_k} the $q_j$ component vanishes in $\operatorname{Im}(\cdot)$. Applying $x_j \mathfrak{f}_j = -i\hbar \nabla_{v_j} \mathfrak{f}_j$ \eqref{eq:state_level_chain_rule} and \eqref{eq:husimifromre} yields:
\begin{align}
\label{eq:Text_advection_derivation}
    &\tfrac{c_{N,k}}{(2\pi)^{3k}} \tfrac{2}{\hbar} \operatorname{Im} \big\langle \Psi_{N,t}, \mathfrak{A}_k^\dagger [\mathfrak{f}] \mathfrak{A}_{k,j}\left[ \nabla \Vext(q_j) \cdot (x_j-q_j) \mathfrak{f}_j \right] \Psi_{N,t} \big\rangle \nonumber \\
    &\quad = \tfrac{c_{N,k}}{(2\pi)^{3k}} \, 2 \nabla \Vext(q_j) \cdot \operatorname{Re} \big\langle \Psi_{N,t}, \mathfrak{A}_k^\dagger[\mathfrak{f}]  \mathfrak{A}_{k,j}\left[ \nabla_{v_j} \mathfrak{f}_j \right] \Psi_{N,t} \big\rangle = \nabla \Vext(q_j) \cdot \nabla_{v_j} \mhus_{N,t}^{(k)}.
\end{align}
Rearranging \eqref{eq:Text_exact_Im_k} gives
\begin{align}
\label{eq:Text_final_decomposition}
    \mathcal{T}_{\text{ext}}^{(k)} &= \sum_{j=1}^k \nabla \Vext(q_j) \cdot \nabla_{v_j} \mhus_{N,t}^{(k)} + \sum_{j=1}^k \mathcal{R}_{\text{ext}, j}^{(k)}.
\end{align}

\paragraph{Conclusion: the full \texorpdfstring{$k$}{k}-particle magnetic hierarchy.}
\label{subsubsec:conclusion_hierarchy}
The splitting of $\partial_t \mhus_{N,t}^{(k)}$ into \eqref{eq:hierarchy_commutator_split} with each identity \eqref{eq:T_kin_total_synthesis}, \eqref{eq:Text_final_decomposition} and \eqref{eq:Tintfinal} results in \eqref{eq:magnetic_hierarchy_k}.
\end{proof}

\section{Remainder Bounding Framework}
\label{sec:remainder_bounds}

This section proves the uniform bound on the remainder $\mathbf{R}^{(k)}$ of the hierarchy
\eqref{eq:magnetic_hierarchy_k}. Every component is reduced to one master form and bounded by a non-stationary phase
argument. The reduction of the components is carried out in Section \ref{sec:unified_bounds}.

\begin{proof}[Proof of Theorem \ref{thm:unified_remainders_k}]
We establish uniform bounds on the magnetic, external and interaction remainders through a unified bounding framework.

\paragraph{Two-stage proof strategy.}
\label{sss:stratproofthmremainder}

The $k$-particle quantum remainders are tested against $\Phi \in C_c^\infty(\mathbb{R}^{6k})$, so we define the magnetic index set $\Lambda_{\text{mag}} \coloneqq \{\text{dis}, \text{for}, \text{gau}, \text{dia}, \text{div}\}$ of \eqref{eq:total_remainder_definitions} and the universal index set $\Lambda \coloneqq \Lambda_{\text{mag}} \cup \{\text{ext}, \text{mf}, \text{cor}\}$, such that:
\begin{align}\label{eq:unified_strong_weak_remainders}
\begin{split}
    \mathcal{I}_{\star, j}^{(k)} &\coloneqq \int_{\mathbb{R}^{6k}} \dZk \, \Phi(Z_k) \mathcal{R}_{\star, j}^{(k)}(Z_k), \quad \star \in \Lambda_{\text{mag}}, \\
    \mathcal{I}_{\text{ext}, j}^{(k)} &\coloneqq \int_{\mathbb{R}^{6k}} \dZk \, \Phi(Z_k) \mathcal{R}_{\text{ext}, j}^{(k)}(Z_k), \\
    \mathcal{I}_{\text{mf}, j}^{(k)} &\coloneqq - \int_{\mathbb{R}^{6k}} \dZk \, \nabla_{v_j} \Phi(Z_k) \cdot \mathcal{R}_{\text{mf}, j}^{(k)}(Z_k), \\
    \mathcal{I}_{\text{cor}, i, j}^{(k)} &\coloneqq \int_{\mathbb{R}^{6k}} \dZk \, \Phi(Z_k) \mathcal{R}_{\text{cor}, i, j}^{(k)}(Z_k).
\end{split}
\end{align}
Each component is estimated in two stages.

\paragraph{Stage one: remainder preparation.} \label{sss:remainderpreparation}

Taylor expansion in $x_j - q_j = \mathcal{O}(\sqrt{\hbar})$ writes each component as an integral against $\langle a_X \Psi_{N,t}, a_Y \Psi_{N,t} \rangle$. The Cauchy--Schwarz and Young inequalities decouple the trace: for functions $u$ of $X$ and $w$ of $Y$,
\begin{align*}
\label{eq:stage_one_young}
    \| a_X \Psi_{N,t} \| \, \| a_Y \Psi_{N,t} \| \, |u(X)| \, |w(Y)|
    &\le \tfrac{1}{2} \| a_X \Psi_{N,t} \|^2 |u(X)|^2 + \tfrac{1}{2} \| a_Y \Psi_{N,t} \|^2 |w(Y)|^2,
\end{align*}
and the exchange $X \leftrightarrow Y$ maps the second term onto the first and the phase
$e^{\frac{i}{\hbar} \sum_{n=1}^k (v_n - A(q_n)) \cdot (x_n - y_n)}$ in the $v$-integral onto
$e^{-\frac{i}{\hbar} \sum_{n=1}^k (v_n - A(q_n)) \cdot (x_n - y_n)}$, which gives $\frac{1}{2}\sum_\pm$.

\paragraph{Stage two: master bounding framework.} \label{sss:MasterBoundingFramework}

The following proposition provides the uniform bound and the order of decay for each remainder. The parameters $(\beta_\star, \mrate_\star, \mathcal{G}_{\star, j}, g_\star)$ are listed in Table \ref{tab:unified_remainder_parameters}.
\begin{proposition}[Master Bounding Framework]
\label{thm:master_remainder_bound}
Let $\star \in \Lambda$. Let $\beta_\star \in \mathbb{N}_0$ and let $\mathcal{I}_{\star, j}^{(k)}$ be a remainder bounded by
\begin{align}
\label{eq:master_assumed_integral}
    \left| \mathcal{I}_{\star, j}^{(k)} \right| &\le \tfrac{c_{N,k}}{(2\pi)^{3k}} \hbar^{\mrate_{\star} - 3k/2} \int_{\mathbb{R}^{3k}} (dq)^{\otimes k} \int_{\mathbb{R}^{6k}} \dXkYk
        \left\| a_X \Psi_{N,t} \right\|^2 \big| x_j - y_j \big|^{\beta_\star} g_{\star}\left(\tfrac{x_j-q_j}{\sqrt{\hbar}}\right) \nonumber \\
    &\quad \times \left(\prod\nolimits_{m \neq j}^k \left| f\left(\tfrac{x_m-q_m}{\sqrt{\hbar}}\right) \right|^2 \right)
        \tfrac{1}{2} \sum_{\pm} \left| \int_{\mathbb{R}^{3k}} (dv)^{\otimes k} \mathcal{G}_{\star, j}[\Phi](Z_k) \, e^{ \pm \frac{i}{\hbar} \sum_{n=1}^k (v_n - A(q_n)) \cdot (y_n - x_n) } \right|_\infty,
\end{align}
for a fixed index $j \in \{1, \dots, k\}$, a test observable $\Phi \in C_c^\infty(\mathbb{R}^{6k})$ with $\Phi = 0$ for $|v_n| > R_\Phi$, $1 \le n \le k$, a spatial profile function $g_\star \in L^1(\mathbb{R}^3; \mathbb{R}^+)$, and a function $\mathcal{G}_{\star, j}[\Phi]$ of $Z_k$, $x_j$ and $y_j$, linear in $\Phi$, smooth in $v$, zero for $|v_n| > R_\Phi$, with $v$-derivatives bounded uniformly in $(q, x_j, y_j)$. One of the two terms may carry $\overline{\mathcal{G}_{\star, j}[\Phi]}$ in place of $\mathcal{G}_{\star, j}[\Phi]$.


Then, for every $\epsilon > 0$,
\begin{align}
\label{eq:master_lemma_bound}
    \left| \mathcal{I}_{\star, j}^{(k)} \right| &\le \C(k) \hbar^{\mrate_\star + \beta_\star\alpha - \epsilon},
\end{align}
with $\alpha \coloneqq 1 - \frac{\epsilon}{3k}$ and $s_\star \coloneqq \max\big(4 + \beta_\star, \lceil \frac{3k}{\epsilon} \rceil\big)$.
\end{proposition}
By Stage one, every component $\star \in \Lambda$ is bounded by \eqref{eq:master_assumed_integral}: the $Y$-dependence
remains only in the weight $|x_j - y_j|^{\beta_\star}$ and in the phase of the $v$-integral. For
$\mathcal{I}_{\text{div}, j}^{(k)}$ this holds for its orthogonal and non-linear parts separately (Section
\ref{subsubsec:bound_I_div_j_k}).

The value $5 = 4 + \beta_{\text{dia}}$ in $s = \max\big(5, \lceil 3k/\epsilon \rceil\big)$ comes from the diamagnetic remainder. The rate of Theorem \ref{thm:unified_remainders_k} follows for
\begin{align}
\label{eq:gamma_unified_explicit_defs}
    \epsilon\left(1 + \tfrac{1}{3k}\right) &\le \rate < \tfrac{1}{2}
\end{align}
The sums over $j$, $i \neq j$ and $\star$ enter $\C(k)$:
\begin{align}
\label{eq:R_ultimate_limits_unified_C}
    \left| \int_{\mathbb{R}^{6k}} \dZk \, \Phi(Z_k) \mathbf{R}^{(k)}(Z_k) \right| &\le \sum_{j=1}^k \left( \sum_{\star \in \Lambda_{\text{mag}}} \big| \mathcal{I}_{\star, j}^{(k)} \big| + \big| \mathcal{I}_{\text{ext}, j}^{(k)} \big| + \big| \mathcal{I}_{\text{mf}, j}^{(k)} \big| \right) + \tfrac{\hbar^2}{2} \sum_{1 \le i \neq j \le k} \big| \mathcal{I}_{\text{cor}, i, j}^{(k)} \big| \\
    &\le \C(k) \hbar^{\frac{1}{2} - \rate}.
\end{align}
The correlation term is of order $\hbar^{3-\rate}$, by its prefactor $\frac{\hbar^2}{2}$ and \eqref{eq:master_lemma_bound}.

\begin{proof}[Proof of Proposition \ref{thm:master_remainder_bound}]

To separate the semiclassical scaling from highly oscillatory velocity phases, define the slow-oscillation domain $\Omega_{\hbar}$ and its exterior tail $\Omega_{\hbar}^c$ via the spatial scaling exponent $\alpha \in (0, 1)$ and the maximum norm $|x|_\infty \coloneqq \max_{i \le 3} |x_i|$:
\begin{align}
\label{eq:omega_hbar_def}
    \Omega_{\hbar} &\coloneqq \left\{ x \in \mathbb{R}^3 : |x|_\infty \le \hbar^\alpha \right\} .
\end{align}
The following lemma is \cite[Lemma~2.5]{ChenLeeLiew2021} before inserting $|x|_\infty > \hbar^\alpha$.
\begin{lemma}[Estimate of Oscillation]
\label{lemma:oscillation_estimate}
Let $\phi \in C_c^\infty(\mathbb{R}^3)$. For any integer $s \in \mathbb{N}$ and any spatial coordinate $x \in \Omega_{\hbar}^c$, the following non-stationary phase estimate holds:
\begin{align}
\label{eq:oscillation_bound_single}
    \left| \int_{\mathbb{R}^3} dv \, e^{\frac{i}{\hbar} v \cdot x} \phi(v) \right| &\le \tfrac{\hbar^s}{|x|_\infty^s} \big| \left\| \nabla_v^s \phi \right\|_{L^1(\mathbb{R}^3)} \big|_\infty.
\end{align}
\end{lemma}

\begin{proof} [Proof of Lemma \ref{lemma:oscillation_estimate}]
For $x \in \Omega_{\hbar}^c$, assume $|x|_\infty = |x_1| > \hbar^\alpha$. Over the compact support $\operatorname{supp}(\phi) \subset B_r(0)$, applying the identity $\frac{\hbar}{i x_1} \nabla_{v,1} \big( e^{\frac{i}{\hbar} v_1 x_1} \big) = e^{\frac{i}{\hbar} v_1 x_1}$ and integrating by parts $s$ times over $v_1$ yields:
\begin{align}
\label{eq:oscillation_ibp_bound}
    \left| \int_{\mathbb{R}^3} dv \, e^{\frac{i}{\hbar} v \cdot x} \phi(v) \right| &= \bigg| \int_{-r}^r dv_2 \, e^{\frac{i}{\hbar} v_2 x_2} \int_{-r}^r dv_3 \, e^{\frac{i}{\hbar} v_3 x_3} \left( \tfrac{-\hbar}{i x_1} \right)^{s} \int_{-r}^r dv_1 \, e^{\frac{i}{\hbar} v_1 x_1} \nabla_{v,1}^{s} \phi(v) \bigg| \\
    &\le \tfrac{\hbar^{s}}{|x_1|^{s}} \int_{\mathbb{R}^3} dv \left| \nabla_{v,1}^{s} \phi(v) \right| \le \tfrac{\hbar^{s}}{|x|_\infty^{s}} \big| \left\| \nabla_v^{s} \phi \right\|_{L^1(\mathbb{R}^3)} \big|_\infty. \nonumber
\end{align}
\end{proof}

Starting from \eqref{eq:master_assumed_integral}, we expand the integration over the phase-space partition $\mathbb{R}^3 = \Omega_{\hbar} \cup \Omega_{\hbar}^c$ \eqref{eq:omega_hbar_def}, which decomposes the $k$-fold indicator product into configurations $S \subseteq \{1, \dots, k\}$. Grouping by the exterior particle count $l \coloneqq |S^c|$ yields the binomial multiplicity $\binom{k}{l}$ and bounds the sum by the maximal configuration for each fixed $l$:
\begin{align}
\label{eq:master_post_partition}
    \left| \mathcal{I}_{\star, j}^{(k)} \right| &\le \tfrac{c_{N,k}}{(2\pi)^{3k}} \hbar^{\mrate_\star - 3k/2} \int_{\mathbb{R}^{3k}} (dq)^{\otimes k} \int_{\mathbb{R}^{6k}} \dXkYk
        \prod\nolimits_{n=1}^k \Big( \chi_{\Omega_{\hbar}}(x_n - y_n) + \chi_{\Omega_{\hbar}^c}(x_n - y_n) \Big) \, \big| x_j - y_j \big|^{\beta_\star} \nonumber \\
    &\quad \times g_\star\left(\tfrac{x_j-q_j}{\sqrt{\hbar}}\right) \left( \prod\nolimits_{m \neq j}^k \left| f\left(\tfrac{x_m-q_m}{\sqrt{\hbar}}\right) \right|^2 \right) \left\| a_X \Psi_{N,t} \right\|^2 \nonumber \\
    &\quad \times \tfrac{1}{2} \sum_{\pm} \left| \int_{\mathbb{R}^{3k}} (dv)^{\otimes k} \left[\mathcal{G}_{\star, j}[\Phi](Z_k)\right] \, e^{ \pm \frac{i}{\hbar} \sum_{n=1}^k (v_n - A(q_n)) \cdot (y_n - x_n) } \right|_{\infty} \nonumber \\
    &\le \tfrac{c_{N,k}}{(2\pi)^{3k}} \hbar^{\mrate_\star - 3k/2} \sum_{l=0}^k \binom{k}{l} \max_{S : |S^c|=l} \int_{\mathbb{R}^{3k}} (dq)^{\otimes k} \int_{\mathbb{R}^{6k}} \dXkYk
        \left\| a_X \Psi_{N,t} \right\|^2 \, \big| x_j - y_j \big|^{\beta_\star} \nonumber \\
    &\quad \times \left( \prod\nolimits_{n=1}^{k-l} \chi_{\Omega_{\hbar}}(x_n - y_n) \right) \left( \prod\nolimits_{r=k-l+1}^k \chi_{\Omega_{\hbar}^c}(x_r - y_r) \right)
        g_\star\left(\tfrac{x_j-q_j}{\sqrt{\hbar}}\right) \left( \prod\nolimits_{m \neq j}^k \left| f\left(\tfrac{x_m-q_m}{\sqrt{\hbar}}\right) \right|^2 \right) \nonumber \\
    &\quad \times \tfrac{1}{2} \sum_{\pm} \left| \int_{\mathbb{R}^{3k}} (dv)^{\otimes k} \left[\mathcal{G}_{\star, j}[\Phi](Z_k)\right] \, e^{ \pm \frac{i}{\hbar} \sum_{n=1}^{k} (v_n - A(q_n)) \cdot (y_n - x_n) } \right|_{\infty}.
\end{align}
We integrate in the order $v \to y \to q \to x$.

\emph{Integration in $v$.}
After the supremum over $q_1, \dots, q_k$, $v_1, \dots, v_{k-l}$, $x_j$ and $y_j$, the non-stationary phase estimate applies to both phases. Applying Lemma \ref{lemma:oscillation_estimate} to the $l$ exterior particles with $|x_r - y_r|_\infty > \hbar^\alpha$ decouples the integrals and bounds the interior velocities over the compact support $B_{R_\Phi}(0)$:
\begin{align}
\label{eq:master_vel_int_outside_bound}
    & \tfrac{1}{2} \sum_{\pm} \left| \int_{\mathbb{R}^{3k}} (dv)^{\otimes k} \mathcal{G}_{\star, j}[\Phi](Z_k) \, e^{ \pm\frac{i}{\hbar} \sum_{n=1}^k (v_n - A(q_n)) \cdot (y_n - x_n) } \right|_\infty \nonumber \\
    &\le \tfrac{1}{2} \sum_{\pm} \left( \int_{B_{R_\Phi}(0)^{k-l}} \prod\nolimits_{n=1}^{k-l} dv_n \right) \sup_{\substack{q_1, \dots, q_k, \, x_j, \, y_j \\ v_1, \dots, v_{k-l}}} \Bigg| \int_{\mathbb{R}^{3l}} \left(\prod\nolimits_{r=k-l+1}^k dv_r\right) \nonumber \\
    &\quad \times \mathcal{G}_{\star, j}[\Phi](Z_k) \, e^{ \pm \frac{i}{\hbar} \sum_{r=k-l+1}^k v_r \cdot (y_r - x_r) } \Bigg|_\infty \nonumber \\
    &\le \tfrac{1}{2} \sum_{\pm} |B_{R_\Phi}(0)|^{k-l} \max_{0 \le l \le k} \left\| \, \left| \left( \prod\nolimits_{r=k-l+1}^k \nabla_{v_r}^{s_\star} \right) \mathcal{G}_{\star, j}[\Phi] \right|_\infty \right\|_{L^\infty(\mathbb{R}^{6k-3l}; L^1(\mathbb{R}^{3l}))}
        \left( \prod\nolimits_{r=k-l+1}^k \tfrac{\hbar^{s_\star}}{|x_r - y_r|_\infty^{s_\star}} \right) \nonumber \\
    &\le \left( \tfrac{2 \pi e R_\Phi^2}{3} \right)^{3(k-l)/2} \max_{0 \le l \le k} \left\| \, \left| \left( \prod\nolimits_{r=k-l+1}^k \nabla_{v_r}^{s_\star} \right) \mathcal{G}_{\star, j}[\Phi] \right|_\infty \right\|_{L^\infty(\mathbb{R}^{6k-3l}; L^1(\mathbb{R}^{3l}))} \left( \prod\nolimits_{r=k-l+1}^k \tfrac{\hbar^{s_\star}}{|x_r - y_r|_\infty^{s_\star}} \right) .
\end{align}
The $L^\infty$ norm includes the supremum over $x_j$ and $y_j$.

\emph{Integration in $y$.}
After substituting \eqref{eq:master_vel_int_outside_bound} into \eqref{eq:master_post_partition}, we apply the translation $z = x_r - y_r$ . The exterior domain $\Omega_{\hbar}^c$ forces $|z|_\infty > \hbar^\alpha$:
\begin{align}\label{eq:master_Y_integrals_rigorous}
    \int_{\mathbb{R}^3} dy_r \, \chi_{\Omega_{\hbar}^c}(x_r-y_r) \tfrac{\hbar^{s_\star}}{|x_r - y_r|_\infty^{s_\star}} &= \hbar^{s_\star} \int_{\hbar^\alpha}^\infty dr \, \tfrac{24 r^2}{r^{s_\star}} = \tfrac{24}{s_\star-3} \hbar^{3\alpha + s_\star(1-\alpha)}, \nonumber \\
    \int_{\mathbb{R}^3} dy_n \, \chi_{\Omega_{\hbar}}(x_n-y_n) &= \int_{|z|_\infty \le \hbar^\alpha} dz = 8\hbar^{3\alpha}.
\end{align}
The index $j$ carries the weight $|x_j - y_j|^{\beta_\star}$, and with $s_\star > 3 + \beta_\star$
from Proposition \ref{thm:master_remainder_bound} its integrals are bounded by:
\begin{align}
\label{eq:master_Y_integrals_weighted}
    \int_{\mathbb{R}^3} dy_j \, \chi_{\Omega_\hbar}(x_j-y_j) \, \big| x_j - y_j \big|^{\beta_\star} &\le \big( \sqrt{3} \hbar^{\alpha} \big)^{\beta_\star} \int_{|z|_\infty \le \hbar^\alpha} dz = \Big( 8 \hbar^{3\alpha} \Big) \, 3^{\beta_\star/2} \hbar^{\beta_\star\alpha}, \nonumber \\
    \int_{\mathbb{R}^3} dy_j \, \chi_{\Omega_\hbar^c}(x_j-y_j) \, \big| x_j - y_j \big|^{\beta_\star} \tfrac{\hbar^{s_\star}}{|x_j - y_j|_\infty^{s_\star}} &\le 3^{\beta_\star/2} \hbar^{s_\star} \int_{\hbar^\alpha}^\infty dr \, \tfrac{24 r^{2+\beta_\star}}{r^{s_\star}} \nonumber \\
    &\quad = \left( \tfrac{24}{s_\star-3} \hbar^{3\alpha + s_\star(1-\alpha)} \right) 3^{\beta_\star/2} \tfrac{s_\star-3}{s_\star-3-\beta_\star} \, \hbar^{\beta_\star\alpha}.
\end{align}
The factor $3^{\beta_\star/2} \frac{s_\star-3}{s_\star-3-\beta_\star} \hbar^{\beta_\star\alpha}$ holds whether $j$ is in the interior set $S$ or the exterior set $S^c$, so it comes out of the partition sum $\sum_{l=0}^k \binom{k}{l} \max_{S : |S^c|=l}$.

\emph{Integration in $q$ and $x$.} The $q$-integrals yield $\hbar^{3k/2} \|g_\star\|_{L^1(\mathbb{R}^3)}$. The trace over $X$ evaluates via the global particle number operator $\mathcal{N}$, and the normalization $c_{N,k}$ \eqref{eq:kinetic_husimi_M_def} cancels the falling factorial:
\begin{align}
\label{eq:master_spatial_trace_evaluation}
    &\tfrac{c_{N,k}}{(2\pi)^{3k}} \int_{\mathbb{R}^{3k}} (dx)^{\otimes k} \left\| a_{x_k} \dots a_{x_1} \Psi_{N,t} \right\|^2
        \int_{\mathbb{R}^{3k}} (dq)^{\otimes k} \left[ g_\star\left(\tfrac{x_j-q_j}{\sqrt{\hbar}}\right) \prod\nolimits_{m \neq j}^k \left| f\left(\tfrac{x_m-q_m}{\sqrt{\hbar}}\right) \right|^2 \right] \nonumber \\
    &= \tfrac{c_{N,k}}{(2\pi)^{3k}} \|g_\star\|_{L^1(\mathbb{R}^3)} \hbar^{3k/2} \left\langle \Psi_{N,t}, \mathcal{N}(\mathcal{N}-1)\dots(\mathcal{N}-k+1) \Psi_{N,t} \right\rangle \nonumber \\
    &= \tfrac{\|g_\star\|_{L^1(\mathbb{R}^3)}}{(2\pi)^{3k}} \hbar^{3k/2} N^k.
\end{align}
Note that on the $N$-particle sector
\begin{align*}
      \left\langle \Psi_{N,t}, \mathcal{N}(\mathcal{N}-1)\dots(\mathcal{N}-k+1) \Psi_{N,t} \right\rangle = N(N-1)\dots(N-k+1) = N^k / c_{N,k}. \nonumber
\end{align*}
Applying the bound \eqref{eq:master_vel_int_outside_bound}, the $L^\infty$ volumes \eqref{eq:master_Y_integrals_rigorous}--\eqref{eq:master_Y_integrals_weighted}, and trace bound \eqref{eq:master_spatial_trace_evaluation} into \eqref{eq:master_post_partition}, the combinatorial sum resolves via the binomial theorem:
\begin{align}
\label{eq:master_final_resolved}
    \left| \mathcal{I}_{\star, j}^{(k)} \right| &\le \tfrac{\|g_\star\|_{L^1(\mathbb{R}^3)}}{(2\pi)^{3k}} \max_{0 \le l \le k} \left\| \left( \prod\nolimits_{r=k-l+1}^k \nabla_{v_r}^{s_\star} \right) \mathcal{G}_{\star, j}[\Phi] \right\|_{L^\infty(\mathbb{R}^{6k-3l}; L^1(\mathbb{R}^{3l}))}
        \times 3^{\beta_\star/2} \tfrac{s_\star-3}{s_\star-3-\beta_\star} \, \hbar^{\beta_\star\alpha} \hbar^{\mrate_\star - 3k/2} \nonumber \\
    &\quad \times \left[ \sum_{l=0}^k \binom{k}{l} \left( 8\left(\tfrac{2 \pi e R_\Phi^2}{3}\right)^{3/2}\hbar^{3\alpha} \right)^{k-l} \left( \tfrac{24}{s_\star-3} \hbar^{3\alpha + s_\star(1-\alpha)} \right)^l \right]
        \left( \hbar^{3k/2} \hbar^{-3k} \right) \nonumber \\
    &= \tfrac{\|g_\star\|_{L^1(\mathbb{R}^3)}}{(2\pi)^{3k}} \max_{0 \le l \le k} \left\| \left( \prod\nolimits_{r=k-l+1}^k \nabla_{v_r}^{s_\star} \right) \mathcal{G}_{\star, j}[\Phi] \right\|_{L^\infty(\mathbb{R}^{6k-3l}; L^1(\mathbb{R}^{3l}))}
        \times 3^{\beta_\star/2} \tfrac{s_\star-3}{s_\star-3-\beta_\star} \nonumber \\
    &\quad \times \hbar^{\mrate_\star + \beta_\star\alpha} \left( 8\left(\tfrac{2 \pi e R_\Phi^2}{3}\right)^{3/2}\hbar^{3\alpha - 3} + \tfrac{24}{s_\star-3} \hbar^{3\alpha - 3 + s_\star(1-\alpha)} \right)^k.
\end{align}

We introduce a parameter $\epsilon > 0$ and set $\alpha \coloneqq 1 - \frac{\epsilon}{3k}$, this fixes $3\alpha - 3 = -\frac{\epsilon}{k}$. Choosing the derivative count $s_\star \coloneqq \max\left(4 + \beta_\star, \left\lceil \frac{3k}{\epsilon} \right\rceil\right)$ ensures $s_\star(1-\alpha) = s_\star\left(\frac{\epsilon}{3k}\right) \ge 1$, so that the exterior tail decays at least as fast as the core:
\begin{align}\label{eq:master_final_exponent_annihilation}
    \left| \mathcal{I}_{\star, j}^{(k)} \right| &\le \tfrac{\|g_\star\|_{L^1(\mathbb{R}^3)}}{(2\pi)^{3k}} \max\nolimits_{0 \le l \le k} \left\| \left( \prod\nolimits_{r=k-l+1}^k \nabla_{v_r}^{s_\star} \right) \mathcal{G}_{\star, j}[\Phi] \right\|_{L^\infty(\mathbb{R}^{6k-3l}; L^1(\mathbb{R}^{3l}))} \nonumber \\
    &\quad \times \hbar^{\mrate_\star + \beta_\star\alpha} \left( \left[ 8\left(\tfrac{2 \pi e R_\Phi^2}{3}\right)^{3/2} + \tfrac{24}{s_\star-3} \right] \hbar^{-\tfrac{\epsilon}{k}} \right)^k 3^{\beta_\star/2} \tfrac{s_\star-3}{s_\star-3-\beta_\star} \nonumber \\
    &\le \C(k) \, \hbar^{\mrate_\star + \beta_\star\alpha - \epsilon}.
\end{align}
\end{proof}

The parameters $(\beta_\star,\mrate_\star,g_\star,\mathcal{G}_{\star,j})$ for each component of the Master Bounding Framework, Proposition \ref{thm:master_remainder_bound}, are:
\begin{table}[!htbp]
\renewcommand{\arraystretch}{1.15}
\resizebox{\textwidth}{!}{%
\begin{tabular}{l|c|c|c|c|l|l}
\hline
\textbf{Component ($\star$)} & \textbf{Source} & \textbf{Section} & $\bm{\beta_\star}$ & $\bm{\mrate_\star}$ & $\mathbf{g_\star(x)}$ & $\bm{\mathcal{G}_{\star, j}[\Phi](Z_k)}$ \\
\hline
Dispersion (dis) & \eqref{eq:I_dis_fully_expanded} & \ref{subsubsec:bound_I_dis_j_k} & $0$ & $1/2$ & $\frac{1}{2} \left( |f(x)|^2 + |\nabla f(x)|^2 \right)$ & $\left( \nabla_{q_j} + (\nabla A(q_j))^T \nabla_{v_j} \right) \Phi(Z_k)$ \\
Force (for) & \eqref{eq:I_for_abs_partitioned} & \ref{subsubsec:bound_I_for_j_k} & $0$ & $1/2$ & $\frac{3\sqrt{3}}{2} \|\nabla^2 A\|_{L^\infty(\mathbb{R}^3)} \left( 1 + |x|^2 \right) |f(x)|^2$ & $\nabla_{v_j} \Phi(Z_k) v_j^T + \Phi(Z_k) \mathbf{I}_{3\times 3}$ \\
Gauge (gau) & \eqref{eq:I_gau_abs_bound} & \ref{subsubsec:bound_I_gau_j_k} & $0$ & $1/2$ & $\frac{1}{2} \|\nabla^2 A\|_{L^\infty(\mathbb{R}^3)} \left( 1 + |x|^2 \right) |f(x)|^2$ & $\Phi(Z_k)$ \\
Diamagnetic (dia) & \eqref{eq:I_dia_abs_bound} & \ref{subsubsec:bound_I_dia_j_k} & $1$ & $-1/2$ & $\|\nabla A\|_{L^\infty(\mathbb{R}^3)}^2 \left( 1 + |x|^2 \right) |f(x)|^2$ & $\Phi(Z_k)$ \\
Orthogonal (Orth) & \eqref{eq:I_Orth_master_form} & \ref{subsubsec:bound_I_div_j_k} & $0$ & $1/2$ & $|f(x)|^2$ & $\tilde{\Phi}_{\text{Orth}}(Z_k)$, see \eqref{eq:tilde_Phi_Orth_def} \\
Non-linear (NL) & \eqref{eq:I_NL_resolved} & \ref{subsubsec:bound_I_div_j_k} & $0$ & $1/2$ & $\frac{1}{2} \|\nabla^2 A\|_{L^\infty(\mathbb{R}^3)} \left( |x|^4 \left| \nabla f(x) \right|^2 + |f(x)|^2 \right)$ & $\Phi(Z_k)$ \\
External (ext) & \eqref{eq:I_ext_abs_bound} & \ref{subsubsec:bound_I_ext_j_k} & $0$ & $1/2$ & $\frac{1}{2} \|\nabla^2 \Vext\|_{L^\infty(\mathbb{R}^3)} \left( 1 + |x|^2 \right) |f(x)|^2$ & $\nabla_{v_j} \Phi(Z_k)$ \\
Mean-Field (mf) & \eqref{eq:I_mf_final_master_bound} & \ref{subsubsec:bound_I_mf_j_k} & $0$ & $1/2$ & $\frac{1}{2} \|\nabla^2 W\|_{L^\infty(\mathbb{R}^3)} \big( 1 + \| |\cdot| f \|_{L^2(\mathbb{R}^3)}^2 \big) \left( 1 + |x|^2 \right) |f(x)|^2$ & $\nabla_{v_j} \Phi(Z_k)$ \\
Correlation (cor) & \eqref{eq:I_cor_abs_bound} & \ref{subsubsec:bound_I_cor_j_k} & $0$ & $1$ & $\|\nabla W\|_{L^\infty(\mathbb{R}^3)} |f(x)|^2$ & $(\nabla_{v_i} - \nabla_{v_j}) \Phi(Z_k)$ \\
\hline
\end{tabular}%
}
\caption{The column Source references the absolute bound of the component, the column Section the subsection with
its proof. Table drafted with the assistance of AI.}
\label{tab:unified_remainder_parameters}
\end{table}

\end{proof}

\section{Convergence to the Magnetic Vlasov Equation}
\label{sec:combined_limit}

This section proves the convergence to the magnetic Vlasov equation, as in \cite{ChenLeeLiew2021} for $A = 0$. The moment
bound of Section \ref{sec:apriori_bounds} gives weak compactness in $L^1$. The remainder estimates identify every
limit point as an admissible solution of the infinite magnetic hierarchy, and its uniqueness selects the tensor
family of the Vlasov solution. The section ends with the proof of Theorem \ref{thm:magnetic_vlasov_limit}, the main
result.

\subsection{Second Phase-Space Moment Bound}
\label{sec:apriori_bounds}
This subsection bounds the second phase-space moment of the magnetic Husimi measures, uniformly in $N$ with growth
$(1 + t)^2$. The bound rests on the conservation of the energy per particle, bounded by \eqref{eq:assumption_initial_data},
\begin{align}
\label{eq:energy_conservation}
    \tfrac{1}{N} \langle \Psi_{N,t}, \mathcal{H} \Psi_{N,t} \rangle
    &= \tfrac{1}{N} \langle \Psi_{N,0}, \mathcal{H} \Psi_{N,0} \rangle \le \C .
\end{align}

Since $\Vext \ge -c_1$ and $W \ge -\|W\|_{L^\infty(\mathbb{R}^3)}$, the kinetic energy per particle is bounded uniformly in
$N$ and $t \ge 0$,
\begin{align}
\label{eq:kinetic_energy_uniform}
    \tfrac{1}{N} \langle \Psi_{N,t}, \mathcal{K}_A \Psi_{N,t} \rangle
    &= \tfrac{1}{N} \langle \Psi_{N,t}, \mathcal{H} \Psi_{N,t} \rangle
    - \tfrac{1}{N} \langle \Psi_{N,t}, \mathcal\Vext \Psi_{N,t} \rangle
    - \tfrac{1}{N} \langle \Psi_{N,t}, \mathcal{W} \Psi_{N,t} \rangle \nonumber \\
    &\le \C + \tfrac{c_1}{N} \langle \Psi_{N,t}, \mathcal{N} \Psi_{N,t} \rangle
    + \tfrac{\|W\|_{L^\infty(\mathbb{R}^3)}}{2N^2} \langle \Psi_{N,t}, \mathcal{N}(\mathcal{N}-1) \Psi_{N,t} \rangle \nonumber\\
    &\le \C.
\end{align}
The moment estimate for the magnetic Husimi measures follows.

\begin{proposition}[Uniform $k$-Particle Phase-Space Moment Bound]
\label{prop:finite_moments_k}
Let Assumptions \ref{assump:A1}--\ref{assump:A4} and \ref{assump:I1} hold. For
$k \in \{1, \dots, N\}$ and $t \ge 0$ it holds
\begin{align}
\label{eq:prop_k_moment_bound}
    \int_{\mathbb{R}^{6k}} \dZk \, \sum_{j=1}^k \Big( |q_j|^2 + |v_j|^2 \Big)
    \mhus_{N,t}^{(k)}(Z_k)
    &\le \C k \big( 1 + t \big)^2 .
\end{align}
\end{proposition}

\begin{proof}[Proof of Proposition \ref{prop:finite_moments_k}]
For $k=1$, using
\begin{align}
\label{eq:kinetic_velocity_expansion}
    v^2 e^{\frac{i}{\hbar} v \cdot (x-y)} &= (-i\hbar \nabla_x) \cdot (i\hbar \nabla_y) e^{\frac{i}{\hbar} v \cdot (x-y)},
\end{align}
we can integrate by $v$ and produce $\delta(x-y)$:
\begin{align}
\label{eq:proof_kin_husimi_derivation1}
    & \iint_{\mathbb{R}^6} dq \, dv \, v^2 \mhus_{N,t}^{(1)}(q,v) \nonumber \\
    &= \tfrac{\hbar^{-3/2}}{(2\pi)^3} \int_{\mathbb{R}^3} dq \iint_{\mathbb{R}^6} dx \, dy \int_{\mathbb{R}^3} dv \, f\left(\tfrac{x-q}{\sqrt{\hbar}}\right) f\left(\tfrac{y-q}{\sqrt{\hbar}}\right) \nonumber \\
    &\qquad \times e^{-\frac{i}{\hbar} A(q) \cdot (x-y)} (-i\hbar \nabla_x) \cdot (i\hbar \nabla_y) e^{\frac{i}{\hbar} v \cdot (x-y)} \langle \Psi_{N,t}, a_x^\dagger a_y \Psi_{N,t} \rangle \nonumber \\
    &= \hbar^{3/2} \int_{\mathbb{R}^3} dq \int_{\mathbb{R}^3} dx \, \left\| \left( -i\hbar \nabla_x + A(q) \right) \left[ f\left(\tfrac{x-q}{\sqrt{\hbar}}\right) a_x \Psi_{N,t} \right] \right\|^2 \nonumber \\
    &= \hbar^{3/2} \int_{\mathbb{R}^3} dq \int_{\mathbb{R}^3} dx \, \Big\| f\left(\tfrac{x-q}{\sqrt{\hbar}}\right) \mpi_x a_x \Psi_{N,t} - \big(A(x) - A(q)\big) f\left(\tfrac{x-q}{\sqrt{\hbar}}\right) a_x \Psi_{N,t}
        - i \sqrt{\hbar} \, (\nabla f)\left(\tfrac{x-q}{\sqrt{\hbar}}\right) a_x \Psi_{N,t} \Big\|^2 \nonumber \\
    &= \hbar^{3/2} \int_{\mathbb{R}^3} dx \int_{\mathbb{R}^3} dq \, \Bigg( \Big\| f\left(\tfrac{x-q}{\sqrt{\hbar}}\right) \mpi_x a_x \Psi_{N,t} \nonumber \\
    &\qquad \qquad - \big(A(x) - A(q)\big) f\left(\tfrac{x-q}{\sqrt{\hbar}}\right) a_x \Psi_{N,t} \Big\|^2 + \hbar \left| (\nabla f)\left(\tfrac{x-q}{\sqrt{\hbar}}\right) \right|^2 \| a_x \Psi_{N,t} \|^2 \Bigg),
\end{align}
one cross term, a gradient in $q$, integrates to zero, and the other one with $A(x) - A(q)$ is purely imaginary. Using the kinetic energy bound \eqref{eq:kinetic_energy_uniform} in the last line we have
\begin{align}
\label{eq:proof_kin_husimi_derivation}
    \iint_{\mathbb{R}^6} dq \, dv \, v^2 \mhus_{N,t}^{(1)}(q,v)
    &\le 2 \hbar^3 \int_{\mathbb{R}^3} dx \, \| \mpi_x a_x \Psi_{N,t} \|^2 +\, \hbar^4 \int_{\mathbb{R}^3} dx \, \| a_x \Psi_{N,t} \|^2 \nonumber \\
    &\qquad \times \Big( \|\nabla f\|_{L^2(\mathbb{R}^3)}^2 + 2 \|\nabla A\|_{L^\infty(\mathbb{R}^3)}^2 \big\| |\cdot| f \big\|_{L^2(\mathbb{R}^3)}^2 \Big)
     \nonumber \\
    &\le \tfrac{4}{N} \langle \Psi_{N,t}, \mathcal{K}_A \Psi_{N,t} \rangle
    + \C \hbar \le \C \big( 1  + \hbar \big) .
\end{align}
The integration over $v$ produces
$(2\pi\hbar)^3 \delta(x-y)$ by \eqref{eq:resolution_identity}:
\begin{align}
\label{eq:husimi_spatial_density}
    \iint_{\mathbb{R}^6} dq \, dv \, |q|^2 \, \mhus_{N,t}^{(1)}(q,v)
    &= \hbar^{\frac{3}{2}} \int_{\mathbb{R}^3} dq \, |q|^2 \int_{\mathbb{R}^3} dx \, f^2\big( \tfrac{x-q}{\sqrt{\hbar}} \big) 
    \langle \Psi_{N,t}, a_x^\dagger a_x \Psi_{N,t} \rangle .
\end{align}
By \eqref{eq:Schrodinger_scaled} only $\mathcal{K}_A$ of \eqref{eq:H_fock_def} contributes to the time derivative
of \eqref{eq:husimi_spatial_density}. With $[\mpi_x, g(x)] = -i\hbar \nabla g(x)$ and the symmetry of $\mpi_x$,
\begin{align}
\label{eq:spatial_moment_operator_derivative}
    \tfrac{d}{dt} \iint_{\mathbb{R}^6} dq \, dv \, |q|^2 \, \mhus_{N,t}^{(1)}(q,v) &= i \, \hbar^{1/2} \int_{\mathbb{R}^3} dq \, |q|^2 \, \Big\langle \Psi_{N,t}, \Big[ \mathcal{K}_A , \int_{\mathbb{R}^3} dx \, f^2\big( \tfrac{x-q}{\sqrt{\hbar}} \big) \, a_x^\dagger a_x \Big] \Psi_{N,t} \Big\rangle \nonumber \\
    &= i \, \hbar^{1/2} \int_{\mathbb{R}^3} dq \, |q|^2 \int_{\mathbb{R}^3} dx \, \Big\langle a_x \Psi_{N,t}, \Big[ \mpi_x^2 , f^2\big( \tfrac{x-q}{\sqrt{\hbar}} \big) \Big] a_x \Psi_{N,t} \Big\rangle \nonumber \\
    &= \hbar^{3/2} \int_{\mathbb{R}^3} dq \, |q|^2 \int_{\mathbb{R}^3} dx \, \nabla_x f^2\big( \tfrac{x-q}{\sqrt{\hbar}} \big) \cdot \operatorname{Re} \big\langle a_x \Psi_{N,t}, \mpi_x a_x \Psi_{N,t} \big\rangle \nonumber \\
    &= - \hbar^{3/2} \int_{\mathbb{R}^3} dx \int_{\mathbb{R}^3} dq \, |q|^2 \, \nabla_q f^2\big( \tfrac{x-q}{\sqrt{\hbar}} \big) \cdot \operatorname{Re} \big\langle a_x \Psi_{N,t}, \mpi_x a_x \Psi_{N,t} \big\rangle \nonumber \\
    &= 2 \hbar^{3/2} \int_{\mathbb{R}^3} dx \int_{\mathbb{R}^3} dq \, q \, f^2\big( \tfrac{x-q}{\sqrt{\hbar}} \big) \cdot \operatorname{Re} \big\langle a_x \Psi_{N,t}, \mpi_x a_x \Psi_{N,t} \big\rangle \nonumber \\
    &\le 2 \hbar^{3/2} \int_{\mathbb{R}^3} dx \int_{\mathbb{R}^3} dq \, |q| \, f^2\big( \tfrac{x-q}{\sqrt{\hbar}} \big) \, \| a_x \Psi_{N,t} \| \, \| \mpi_x a_x \Psi_{N,t} \| \nonumber \\
    &\le 2 \Big( \hbar^{3/2} \int_{\mathbb{R}^3} dx \int_{\mathbb{R}^3} dq \, f^2\big( \tfrac{x-q}{\sqrt{\hbar}} \big) \, \| \mpi_x a_x \Psi_{N,t} \|^2 \Big)^{1/2} \nonumber \\
    &\qquad \times \Big( \hbar^{3/2} \int_{\mathbb{R}^3} dq \, |q|^2 \int_{\mathbb{R}^3} dx \, f^2\big( \tfrac{x-q}{\sqrt{\hbar}} \big) \, \| a_x \Psi_{N,t} \|^2 \Big)^{1/2} \nonumber \\
    &= 2 \Big( \tfrac{2}{N} \langle \Psi_{N,t}, \mathcal{K}_A \Psi_{N,t} \rangle \Big)^{1/2} \Big( \iint_{\mathbb{R}^6} dq \, dv \, |q|^2 \, \mhus_{N,t}^{(1)}(q,v) \Big)^{1/2} \nonumber \\
    &\le \C \Big( 1 + \iint_{\mathbb{R}^6} dq \, dv \, |q|^2 \, \mhus_{N,t}^{(1)}(q,v) \Big)^{1/2} .
\end{align}
Hence
\begin{align*}
\label{eq:spatial_moment_sqrt_inequality}
    &\tfrac{d}{dt} \Big( 1 + \iint_{\mathbb{R}^6} dq \, dv \, |q|^2 \, \mhus_{N,t}^{(1)}(q,v) \Big)^{1/2} = \frac{ \tfrac{d}{dt} \iint_{\mathbb{R}^6} dq \, dv \, |q|^2 \, \mhus_{N,t}^{(1)}(q,v) }
    { 2 \big( 1 + \iint_{\mathbb{R}^6} dq \, dv \, |q|^2 \, \mhus_{N,t}^{(1)}(q,v) \big)^{1/2} } \le \C.
\end{align*}
Integration over $[0,t]$, \eqref{eq:assumption_initial_data} give
\begin{align}
\label{eq:spatial_moment_gronwall_evaluated}
    \iint_{\mathbb{R}^6} dq \, dv \, |q|^2 \, \mhus_{N,t}^{(1)}(q,v)
    \le \C \big( 1 + t \big)^2 .
\end{align}
By the symmetry of $\mhus_{N,t}^{(k)}$ under $(q_i,v_i) \leftrightarrow (q_j,v_j)$ and $k-1$ applications of
\eqref{eq:normalized_husimi_measure}, the left side of \eqref{eq:prop_k_moment_bound} is $k$ times the case $k=1$. With \eqref{eq:proof_kin_husimi_derivation} and \eqref{eq:spatial_moment_gronwall_evaluated}, this gives \eqref{eq:prop_k_moment_bound}.
\end{proof}

\subsection{Weak Compactness} This subsection gives the compactness argument of magnetic Husimi measures at every time $t\in [0,T]$. The proof applies the Dunford--Pettis theorem at fixed time, Cantor's diagonal procedure over rational $t$ and $k$, and the Lipschitz continuity in time in the weak topology.

\begin{proposition}[Weak Compactness]
\label{prop:weakcompactmagnetic}
Let Assumptions \ref{assump:A1}--\ref{assump:A4} and \ref{assump:I1} hold and let $T > 0$. There exist a
subsequence $(N_l)_{l \in \mathbb{N}}$, independent of $t \in [0,T]$ and of $k \ge 1$, and functions
$\mathbf{f}_t^{(k)} \in L^1(\mathbb{R}^{6k})$ such that for every $t \in [0,T]$, every $k \ge 1$ and every
$\varphi \in L^\infty(\mathbb{R}^{6k})$
\begin{align}
\label{eq:weak_compactness_limit}
    \lim\nolimits_{l \to \infty} \int_{\mathbb{R}^{6k}} \dZk \, \mhus_{N_l,t}^{(k)}(Z_k) \, \varphi(Z_k)
    &= \int_{\mathbb{R}^{6k}} \dZk \, \mathbf{f}_t^{(k)}(Z_k) \, \varphi(Z_k) .
\end{align}
\end{proposition}

\begin{proof}[Proof of Proposition \ref{prop:weakcompactmagnetic}]
\emph{Fixed time.} By \eqref{eq:normalized_husimi_measure}, $\|\mhus_{N,t}^{(k)}\|_{L^1(\mathbb{R}^{6k})} = 1$
and, $\|\mhus_{N,t}^{(k)}\|_{L^\infty(\mathbb{R}^{6k})} \le c_{N,k} (2\pi)^{-3k}$, so
$\int_E \dZk \, \mhus_{N,t}^{(k)} \le 2 (2\pi)^{-3k} |E|$ for every measurable $E \subset \mathbb{R}^{6k}$, and
$(\mhus_{N,t}^{(k)})_N$ is uniformly integrable and
\begin{align}
\label{eq:markov_tightness}
    & \int_{\{ \sum_j (|q_j|^2 + |v_j|^2) \ge r \}} \dZk \, \mhus_{N,t}^{(k)}(Z_k)
    \nonumber \\
    &\le \tfrac{1}{r} \int_{\mathbb{R}^{6k}} \dZk \, \sum_{j=1}^k \Big( |q_j|^2 + |v_j|^2 \Big) \mhus_{N,t}^{(k)}(Z_k)
    \le \tfrac{k}{r} \, \C \, ( 1 + T )^2
    \xrightarrow[r \to \infty]{} 0
\end{align}
uniformly in $N$ and $t \in [0,T]$. Hence $(\mhus_{N,t}^{(k)})_N$ is tight, and the Dunford--Pettis theorem
\cite[Thm.~4.30]{brezis2011functional} renders it relatively weakly compact in $L^1(\mathbb{R}^{6k})$ at every
fixed $t$. Therefore every subsequence of $(\mhus_{N,t}^{(k)})_N$ has a further subsequence converging weakly in $L^1(\mathbb{R}^{6k})$.

\emph{Lipschitz bound in $t$.} Let $\Phi \in C_c^\infty(\mathbb{R}^{6k})$ and $0 \le s \le t \le T$. Tested
against $\Phi$ and integrated over $[s,t]$, \eqref{eq:magnetic_hierarchy_k} gets
\begin{align}
\label{eq:husimi_integrated_hierarchy}
    & \Big| \int_{\mathbb{R}^{6k}} \dZk \, \Phi \, \Big( \mhus_{N,t}^{(k)} - \mhus_{N,s}^{(k)} \Big) \Big|
    \nonumber \\
    &= \Bigg| \int_s^t d\tau \sum_{j=1}^k \int_{\mathbb{R}^{6k}} \dZk \, \mhus_{N,\tau}^{(k)}
    \Big( v_j \cdot \nabla_{q_j} \Phi - \big( v_j \times B(q_j) + \nabla \Vext(q_j) \big) \cdot \nabla_{v_j} \Phi \Big)
    \nonumber \\
    &\quad \hphantom{{}= \Bigg|} - \int_s^t d\tau \sum_{j=1}^k \int_{\mathbb{R}^{6k}} \dZk \, \nabla_{v_j} \Phi \cdot
    \mathcal{F}_{\text{int}, j}\big[ \mhus_{N,\tau}^{(k+1)} \big]
 + \int_s^t d\tau \int_{\mathbb{R}^{6k}} \dZk \, \Phi \, \mathbf{R}^{(k)} \Bigg|
    \le \C(k) \, |t - s| .
\end{align}
The first term is bounded by mass one \eqref{eq:normalized_husimi_measure}, since $\operatorname{supp} \Phi$ is
compact and $B$, $\nabla \Vext$ are continuous. The second term is bounded by $\|\nabla W\|_{L^\infty(\mathbb{R}^3)}$
and mass one of $\mhus_{N,\tau}^{(k+1)}$, the third by Theorem \ref{thm:unified_remainders_k}. The constant
$\C(k, \Phi)$ does not depend on $N$, $\hbar$, $s$ and $t$.

\emph{Rational times.} By the first part, Cantor's diagonal procedure over the countable set
$\big( [0,T] \cap \mathbb{Q} \big) \times \mathbb{N}$ gives a sequence $(N_l)_{l \in \mathbb{N}}$ such that, for every
$s \in [0,T] \cap \mathbb{Q}$ and $k \ge 1$, there is $\mathbf{f}_s^{(k)} \in L^1(\mathbb{R}^{6k})$ with
\begin{align}
\label{eq:weak_limit_rational}
    \mhus_{N_l,s}^{(k)} &\rightharpoonup \mathbf{f}_s^{(k)} \quad \text{weakly in } L^1(\mathbb{R}^{6k}) .
\end{align}

\emph{All times.} Let $t \in [0,T]$, $s \in [0,T] \cap \mathbb{Q}$ and $\Phi \in C_c^\infty(\mathbb{R}^{6k})$.
By \eqref{eq:husimi_integrated_hierarchy},
\begin{align}
\label{eq:husimi_cauchy_all_times}
    \Big| \int_{\mathbb{R}^{6k}} \dZk \, \Phi \, \Big( \mhus_{N_l,t}^{(k)} - \mhus_{N_m,t}^{(k)} \Big) \Big| &\le 2 \, \C(k, \Phi) \, |t - s| + \Big| \int_{\mathbb{R}^{6k}} \dZk \, \Phi \, \Big( \mhus_{N_l,s}^{(k)} - \mhus_{N_m,s}^{(k)} \Big) \Big| .
\end{align}
The first term is small for $s$ close to $t$, the second for large $l, m$ by \eqref{eq:weak_limit_rational}.
Hence $\big( \int_{\mathbb{R}^{6k}} \dZk \, \Phi \, \mhus_{N_l,t}^{(k)} \big)_l$ converges for every
$t \in [0,T]$ and $\Phi \in C_c^\infty(\mathbb{R}^{6k})$.

\emph{From $C_c^\infty$ to $L^\infty$.} Since $C_c^\infty(\mathbb{R}^{6k}) \subset L^\infty(\mathbb{R}^{6k})$, the limits coincide.

\end{proof}

\subsection{Convergence to the Infinite Magnetic Hierarchy}

Based on the compactness of the magnetic Husimi measures and the estimates of the remainders, the following
proposition identifies the limit points as admissible solutions of the infinite magnetic hierarchy.

\begin{proposition}[Convergence to the Infinite Magnetic Hierarchy]
\label{prop:main_magnetic_limit}
Let Assumptions \ref{assump:A1}--\ref{assump:A4} and \ref{assump:I1} hold and let $T > 0$. The limit functions
$(\mathbf{f}^{(k)}_t)_{k \ge 1}$ of Proposition \ref{prop:weakcompactmagnetic} form a family in $\mathcal{A}_T$
(Definition \ref{def:admissible_class}) with initial data $(\mathbf{f}^{(k)}_0)_{k \ge 1}$. The same holds for every
subsequence of $(N)$: it has a further subsequence converging weakly in $L^1(\mathbb{R}^{6k})$ for every $k \ge 1$
and $t \in [0,T]$, and the limits form a family in $\mathcal{A}_T$.
\end{proposition}

\begin{proof}[Proof of Proposition \ref{prop:main_magnetic_limit}]
Let $(N_l)$ and $\mathbf{f}^{(k)}_t$ be as in Proposition \ref{prop:weakcompactmagnetic},
$\Phi \in C_c^\infty(\mathbb{R}^{6k})$ and $0 \le s \le t \le T$. In \eqref{eq:husimi_integrated_hierarchy} we can
take the limit $l \to \infty$ due to the boundedness of the integrands and the remainder term vanishes by Theorem \ref{thm:unified_remainders_k},
since $\hbar^{\frac{1}{2} - \rate} \to 0$:
\begin{align}
\label{eq:limit_hierarchy_identity}
    &\int_{\mathbb{R}^{6k}} \dZk \, \Phi \, \Big( \mathbf{f}^{(k)}_t - \mathbf{f}^{(k)}_s \Big) \nonumber \\
    &= \int_s^t d\tau \sum_{j=1}^k \int_{\mathbb{R}^{6k}} \dZk \, \mathbf{f}^{(k)}_\tau
    \Big( v_j \cdot \nabla_{q_j} \Phi - \big( v_j \times B(q_j) + \nabla \Vext(q_j) \big) \cdot \nabla_{v_j} \Phi \Big)
    \nonumber \\
    &\quad - \int_s^t d\tau \sum_{j=1}^k \int_{\mathbb{R}^{6(k+1)}} (dq \, dv)^{\otimes (k+1)} \,
    \nabla_{v_j} \Phi(Z_k) \cdot \nabla W(q_j - q_{k+1}) \, \mathbf{f}^{(k+1)}_\tau(Z_{k+1}) .
\end{align}
The bounds hold for every $N$, so the same applies for every subsequence of $(N)$.

\emph{Admissibility.} Condition \ref{item:admissible_a}: testing \eqref{eq:weak_compactness_limit} with
$\varphi = \One_E$, $E \subset \mathbb{R}^{6k}$, $|E| < \infty$, with $\varphi = 1$ and with $\varphi \circ \pi$ for a
permutation $\pi$ of the $k$ phase-space variables confirms that $\mathbf{f}^{(k)}_t$ is a symmetric probability
density. Testing with $\varphi \otimes 1$ in the $(k+1)$-particle sector and \eqref{eq:normalized_husimi_measure}
gives the consistency $\displaystyle \iint_{\mathbb{R}^6} dq_{k+1} \, dv_{k+1} \, \mathbf{f}^{(k+1)}_t = \mathbf{f}^{(k)}_t$.

Condition \ref{item:admissible_c}: testing \eqref{eq:weak_compactness_limit} at $k = 1$ with $\varphi=\min(|q|^2 + |v|^2, R)$, $R > 0$,  with the bound \eqref{eq:prop_k_moment_bound}, uniform in
$\hbar$ and $t \in [0,T]$ gives \eqref{eq:admissible_moments} as $R \to \infty$.

Condition \ref{item:admissible_b}: for $\Phi \in C_c^\infty(\mathbb{R}^{6k})$, the integrands on the right side of
\eqref{eq:limit_hierarchy_identity} are bounded, so $t \mapsto \displaystyle \int \dZk \, \Phi \, \mathbf{f}^{(k)}_t$
is Lipschitz on $[0,T]$. The tightness from Condition \ref{item:admissible_c}, lifted to the $k$-particle sector by the
symmetry and consistency of Condition \ref{item:admissible_a}, extends the continuity of
$t \mapsto \int \dZk \, \varphi \, \mathbf{f}^{(k)}_t$ from $\varphi \in C_c^\infty(\mathbb{R}^{6k})$ to
$\varphi \in C_b(\mathbb{R}^{6k})$.

Condition \ref{item:admissible_d}: for $\Phi \in C_c^\infty([0,T) \times \mathbb{R}^{6k})$, integrating
\eqref{eq:limit_hierarchy_identity} on $[0,t]$ with the test function $\partial_t \Phi(t, \cdot)$ over $t \in [0,T]$ and
exchanging the time integrals gives \eqref{eq:hierarchy_weak_form}, since $\Phi(T, \cdot) = 0$.

\end{proof}

\subsection{Identification of the Vlasov Solution}
Now we can prove the main result. The uniqueness in the class $\mathcal{A}_T$ identifies the limit
(Proposition \ref{thm:vlasov_infinite_hierarchy}).

\begin{proof}[Proof of Theorem \ref{thm:magnetic_vlasov_limit}]
The bound \eqref{eq:normalized_husimi_measure} passes to the weak limit \eqref{eq:initial_condition_tensor}, so we have
$\mathbf{f}_0 \le (2\pi)^{-3}$. The weak solution $\mathbf{f}_t$ exists (Proposition \ref{thm:vlasov_wellposedness})
and $(\mathbf{f}_t^{\otimes k})_{k \ge 1}$ is the only family in $\mathcal{A}_T$ with initial data
$(\mathbf{f}_0^{\otimes k})_{k \ge 1}$ (Proposition \ref{thm:vlasov_infinite_hierarchy}). By Proposition \ref{prop:main_magnetic_limit}, every subsequence of $(N)$ has a
further subsequence for which $\mhus_{N,t}^{(k)} \rightharpoonup \mathbf{g}_t^{(k)}$ weakly in $L^1(\mathbb{R}^{6k})$
for all $k \ge 1$ and $t \in [0,T]$, with $(\mathbf{g}_t^{(k)})_{k \ge 1} \in \mathcal{A}_T$. By
\eqref{eq:initial_condition_tensor} we have $\mathbf{g}_0^{(k)} = \mathbf{f}_0^{\otimes k}$, so
$\mathbf{g}_t^{(k)} = \mathbf{f}_t^{\otimes k}$ by uniqueness in $\mathcal A_T$. The limit is the same for every subsequence, hence the sequence converges and finally \eqref{eq:full_sequence_convergence} holds.
\end{proof}

\section{Convergence in Wasserstein \texorpdfstring{$\mathbb{W}_2$}{W2}}
\label{sec:wasserstein_noV}

Theorem \ref{thm:magnetic_vlasov_limit} gives weak convergence in $L^1$. This section strengthens it to convergence in $\mathbb{W}_2$ and proves Theorem \ref{thm:w2_energy_criterion_noV}. For $p \in [1,2]$, we use the notation $\mathcal{P}_p(\mathbb{R}^{6k})$
for the probability measures on $\mathbb{R}^{6k}$ with finite $p$-th moment and $\mathbb{W}_p(\mu, \nu)$ for the
$p$-Wasserstein distance of $\mu, \nu \in \mathcal{P}_p(\mathbb{R}^{6k})$ \cite[Def.~6.1, Def.~6.4]{villani2008optimal}.

Based on the result obtained in Section \ref{sec:combined_limit}, together with the fact that  the weak convergence in
$L^1$ of Theorem \ref{thm:magnetic_vlasov_limit} is narrow convergence, the convergence in $\mathbb{W}_2$ is reduced to prove the convergence of the second moment of $\mhus_{N,t}^{(1)}$ to that of $\mathbf{f}_t$, \cite[Thm.~6.9]{villani2008optimal}.
The convergence of the velocity $|v|^2$ moment comes from energy conservation and the energy inequality of
Assumption \ref{assump:I4}. Namely, we relate the quantum and the classical energy (Lemma
\ref{lemma:energy_decomposition_noV}) and show that their difference vanishes as $\hbar \to 0$. For the convergence of the spatial $|q|^2$
moment, we consider the propagation of the spatial tail in Lemma \ref{lem:spatial_tail_propagation} and show that the tail of the $|q|^2$ moment can be bounded by the $|v|^2$
moment on $\{ |q| \ge R/2 \}$, which is controlled by the bound on the energy difference of Lemma
\ref{lemma:energy_decomposition_noV}.

We decompose the energy into the one-body and the interaction energy of the magnetic Husimi measures,
\begin{align}
\label{eq:energy_husimi_decomposition_noV}
    \mathcal{V}_N(t) &\coloneqq \iint_{\mathbb{R}^6} dq \, dv \,
    \Big( \tfrac{1}{2}|v|^2 + \Vext(q) \Big) \mhus_{N,t}^{(1)}(q,v),
    \notag \\
    \mathcal{W}_N(t) &\coloneqq \tfrac{1}{2} \int_{\mathbb{R}^{12}} (dq \, dv)^{\otimes 2} \,
    W(q_1 - q_2) \, \mhus_{N,t}^{(2)}(q_1,v_1,q_2,v_2).
\end{align}
Their sum is the classical energy \eqref{eq:classical_energy_definition} of the magnetic Husimi measures, with $\mhus_{N,t}^{(2)}$ in place of $\mathbf{f}^{\otimes 2}$. The following lemma compares it with the energy per particle
$\frac{1}{N} \langle \Psi_{N,t}, \mathcal{H} \Psi_{N,t} \rangle$, with $\mathcal{H}$ of \eqref{eq:H_fock_def}.

\begin{lemma}[Quantum and Classical Energy Difference]
\label{lemma:energy_decomposition_noV}
Let Assumptions \ref{assump:A1}--\ref{assump:A4} and \ref{assump:I1} hold. For $t \ge 0$ it holds
\begin{align}
\label{eq:energy_decomposition_bound_noV}
    \Big| \tfrac{1}{N} \big\langle \Psi_{N,t}, \mathcal{H} \Psi_{N,t} \big\rangle - \mathcal{V}_N(t) - \mathcal{W}_N(t) \Big|
    &\le \C (1 + t) \sqrt{\hbar} .
\end{align}
\end{lemma}

\begin{proof}[Proof of Lemma \ref{lemma:energy_decomposition_noV}]
As in \eqref{eq:proof_kin_husimi_derivation1}, with the kinetic energy bound \eqref{eq:kinetic_energy_uniform} in the last step:
\begin{align}
\label{eq:energy_dec_kinetic_noV}
    & \left| \tfrac{1}{2} \iint_{\mathbb{R}^6} dq \, dv \, |v|^2 \mhus_{N,t}^{(1)}(q,v) - \tfrac{1}{N} \big\langle \Psi_{N,t}, \mathcal{K}_A \Psi_{N,t} \big\rangle \right| \nonumber \\
    &= \Bigg| \tfrac{\hbar^{3/2}}{2} \int_{\mathbb{R}^3} dx \int_{\mathbb{R}^3} dq \, \Bigg( \big|A(x) - A(q)\big|^2 f^2\left(\tfrac{x-q}{\sqrt{\hbar}}\right) \| a_x \Psi_{N,t} \|^2 \nonumber \\
    &\qquad - 2 \operatorname{Re} \left\langle f\left(\tfrac{x-q}{\sqrt{\hbar}}\right) \mpi_x a_x \Psi_{N,t}, \, \big(A(x) - A(q)\big) f\left(\tfrac{x-q}{\sqrt{\hbar}}\right) a_x \Psi_{N,t} \right\rangle
        + \hbar \left| (\nabla f)\left(\tfrac{x-q}{\sqrt{\hbar}}\right) \right|^2 \| a_x \Psi_{N,t} \|^2 \Bigg) \Bigg| \nonumber \\
    &\le \tfrac{\hbar^3}{2} \int_{\mathbb{R}^3} dx \, \Bigg[ 2 \|\nabla A\|_{L^\infty(\mathbb{R}^3)} \big\| \sqrt{\left| \cdot \right|} f \big\|_{L^2(\mathbb{R}^3)}^2 \sqrt{\hbar} \, \big\| \mpi_x a_x \Psi_{N,t} \big\| \big\| a_x \Psi_{N,t} \big\| \nonumber \\
    &\qquad + \Big( \|\nabla f\|_{L^2(\mathbb{R}^3)}^2 + \|\nabla A\|_{L^\infty(\mathbb{R}^3)}^2 \big\| |\cdot| f \big\|_{L^2(\mathbb{R}^3)}^2 \Big) \hbar \, \big\| a_x \Psi_{N,t} \big\|^2 \Bigg] \nonumber \\
    &\le \C \sqrt{\hbar} \, \Big( \tfrac{2}{N} \big\langle \Psi_{N,t}, \mathcal{K}_A \Psi_{N,t} \big\rangle \Big)^{1/2} + \C \hbar \le \C \sqrt{\hbar} .
\end{align}
Taylor's formula at $x$ with $x - s \sqrt{\hbar} u$, $s \in (0,1)$:
\begin{align}
\label{eq:energy_dec_external_expansion_noV}
    &\iint_{\mathbb{R}^6} dq \, dv \, \Vext(q) \, \mhus_{N,t}^{(1)}(q,v) \nonumber \\
    &= \hbar^{3/2} \int_{\mathbb{R}^3} dq \, \Vext(q)
    \int_{\mathbb{R}^3} dx \, f^2\big( \tfrac{x-q}{\sqrt{\hbar}} \big) \big\langle \Psi_{N,t}, a_x^\dagger a_x \Psi_{N,t} \big\rangle \nonumber \\
    &= \hbar^3 \int_{\mathbb{R}^3} dx \, \big\langle \Psi_{N,t}, a_x^\dagger a_x \Psi_{N,t} \big\rangle
    \int_{\mathbb{R}^3} du \, f^2(u) \, \Vext\big( x - \sqrt{\hbar} u \big) \nonumber \\
    &= \tfrac{1}{N} \big\langle \Psi_{N,t}, \mathcal{V} \Psi_{N,t} \big\rangle
    + \hbar^3 \int_{\mathbb{R}^3} dx \, \big\langle \Psi_{N,t}, a_x^\dagger a_x \Psi_{N,t} \big\rangle
    \Big( - \sqrt{\hbar} \, \nabla \Vext(x) \cdot \int_{\mathbb{R}^3} du \, f^2(u) \, u \nonumber \\
    &\qquad\qquad\qquad\qquad\qquad + \hbar \int_0^1 ds \, (1-s) \int_{\mathbb{R}^3} du \, f^2(u) \, u^{T} \nabla^2 \Vext(x - s \sqrt{\hbar} u) u \Big) .
\end{align}
This bounds the difference:
\begin{align}
\label{eq:energy_dec_external_bound_noV}
    &\left| \iint_{\mathbb{R}^6} dq \, dv \, \Vext(q) \, \mhus_{N,t}^{(1)}(q,v) - \tfrac{1}{N} \big\langle \Psi_{N,t}, \mathcal{V} \Psi_{N,t} \big\rangle \right| \nonumber \\
    &\le \hbar^3 \int_{\mathbb{R}^3} dx \, \big\langle \Psi_{N,t}, a_x^\dagger a_x \Psi_{N,t} \big\rangle
 \Big( \sqrt{\hbar} \, \big\| \sqrt{| \cdot |} f \big\|_{L^2(\mathbb{R}^3)}^2 \big( |\nabla \Vext(0)| + \|\nabla^2 \Vext\|_{L^\infty(\mathbb{R}^3)} |x| \big) \nonumber \\
    &\qquad + \tfrac{\hbar}{2} \|\nabla^2 \Vext\|_{L^\infty(\mathbb{R}^3)} \big\| |\cdot| f \big\|_{L^2(\mathbb{R}^3)}^2 \Big) \nonumber \\
    &\le \C \sqrt{\hbar} \, \Big( 1 + \hbar^3 \int_{\mathbb{R}^3} dx \, |x| \, \big\langle \Psi_{N,t}, a_x^\dagger a_x \Psi_{N,t} \big\rangle \Big) \nonumber \\
    &= \C \sqrt{\hbar} \, \Big( 1 + \hbar^{3/2} \iint_{\mathbb{R}^6} dq \, dx \, |x| \, f^2\big(\tfrac{x-q}{\sqrt{\hbar}}\big)
    \big\langle \Psi_{N,t}, a_x^\dagger a_x \Psi_{N,t} \big\rangle \Big) \nonumber \\
    &\le \C \sqrt{\hbar} \, \Big( 1 + \iint_{\mathbb{R}^6} dq \, dv \, |q| \, \mhus_{N,t}^{(1)}(q, v)
    + \sqrt{\hbar} \, \big\| \sqrt{| \cdot |} f \big\|_{L^2(\mathbb{R}^3)}^2 \Big) \nonumber \\
    &\le \C \sqrt{\hbar} \, \Big( 1 + \Big( \iint_{\mathbb{R}^6} dq \, dv \, |q|^2 \, \mhus_{N,t}^{(1)}(q,v) \Big)^{1/2} \Big)
    \le \C (1 + t) \sqrt{\hbar} .
\end{align}
The last line used Proposition \ref{prop:finite_moments_k} at $k = 1$. Substituting $q_j = x_j - \sqrt{\hbar} u_j$, i.e. $q_1 - q_2 = x_1 - x_2 - \sqrt{\hbar}(u_1 - u_2)$ gives:
\begin{align}
\label{eq:energy_dec_interaction_bound_noV}
    \left| \mathcal{W}_N(t) - \tfrac{1}{N} \big\langle \Psi_{N,t}, \mathcal{W} \Psi_{N,t} \big\rangle \right| &= \Big| \tfrac{1}{2 N (N-1)} \iint_{\mathbb{R}^6} dx_1 \, dx_2 \, \big\langle \Psi_{N,t}, a_{x_1}^\dagger a_{x_2}^\dagger a_{x_2} a_{x_1} \Psi_{N,t} \big\rangle \nonumber \\
    &\qquad \times \iint_{\mathbb{R}^6} du_1 \, du_2 \, f^2(u_1) f^2(u_2) \, W\big( x_1 - x_2 - \sqrt{\hbar}(u_1 - u_2) \big) \nonumber \\
    &\qquad - \tfrac{1}{2N^2} \iint_{\mathbb{R}^6} dx_1 \, dx_2 \, W(x_1 - x_2) \big\langle \Psi_{N,t}, a_{x_1}^\dagger a_{x_2}^\dagger a_{x_2} a_{x_1} \Psi_{N,t} \big\rangle \Big| \nonumber \\
    &\le \tfrac{1}{2 N (N-1)} \iint_{\mathbb{R}^6} dx_1 \, dx_2 \, \big\langle \Psi_{N,t}, a_{x_1}^\dagger a_{x_2}^\dagger a_{x_2} a_{x_1} \Psi_{N,t} \big\rangle \nonumber \\
    &\qquad \times \iint_{\mathbb{R}^6} du_1 \, du_2 \, f^2(u_1) f^2(u_2) \, \sqrt{\hbar} \, \|\nabla W\|_{L^\infty(\mathbb{R}^3)} \big( |u_1| + |u_2| \big) \nonumber \\
    &\quad + \tfrac{1}{N-1} \left| \tfrac{1}{N} \big\langle \Psi_{N,t}, \mathcal{W} \Psi_{N,t} \big\rangle \right| \nonumber \\
    &\le \sqrt{\hbar} \, \|\nabla W\|_{L^\infty(\mathbb{R}^3)} \big\| \sqrt{|\cdot|} f \big\|_{L^2(\mathbb{R}^3)}^2 + \hbar^3 \|W\|_{L^\infty(\mathbb{R}^3)} \le \C \sqrt{\hbar} .
\end{align}
Combining \eqref{eq:energy_dec_kinetic_noV}, \eqref{eq:energy_dec_external_bound_noV}, \eqref{eq:energy_dec_interaction_bound_noV} yields \eqref{eq:energy_decomposition_bound_noV}.

\end{proof}

The following lemma bounds the spatial $|q|^2$ tail at time $t$ by its value at time $0$ and the velocity $|v|^2$
moment on $\{ |q| \ge R/2 \}$.

\begin{lemma}[Propagation of the Spatial Tail]
\label{lem:spatial_tail_propagation}
Let Assumptions \ref{assump:A1}--\ref{assump:A4} and \ref{assump:I1} hold. For $R > 0$ and $t \ge 0$ it holds
\begin{align}
\label{eq:spatial_tail_bound}
    &\iint_{\{ |q| \ge R \}} dq \, dv \, |q|^2 \, \mhus_{N,t}^{(1)}(q,v)
    \le \iint_{\{ |q| \ge R/2 \}} dq \, dv \, |q|^2 \, \mhus_{N,0}^{(1)}(q,v) \nonumber \\
    &\quad + \C (1 + t) \int_0^t d\tau \, \Big( \iint_{\{ |q| \ge R/2 \}} dq \, dv \, |v|^2 \,
    \mhus_{N,\tau}^{(1)}(q,v) \Big)^{1/2} + \C (1 + t) \sqrt{\hbar} \, t .
\end{align}
\end{lemma}

\begin{proof}[Proof of Lemma \ref{lem:spatial_tail_propagation}]
Fix $\psi \in C^\infty(\mathbb{R}; [0,1])$ with $\psi \equiv 0$ on $(-\infty, \tfrac12]$ and $\psi \equiv 1$ on
$[1, \infty)$. Since $\psi^2 \le \psi$ and $|q| \le R$ on the support of $\psi'( |q| / R )$, $\big| \nabla_q \big( |q|^2 \, \psi( \tfrac{|q|}{R} )^2 \big) \big| \le \C \, |q| \, \psi( \tfrac{|q|}{R} )$, and the gradient vanishes on $\{ |q| \le R/2 \}$. We follow the steps in \eqref{eq:husimi_spatial_density} and \eqref{eq:spatial_moment_operator_derivative} with $|q|^2 \, \psi( \tfrac{|q|}{R} )^2$ in place of $|q|^2$ up to the point where $\nabla_q$ acts on it:
\begin{align}
\label{eq:spatial_tail_operator_derivative}
    & \tfrac{d}{dt} \iint_{\mathbb{R}^6} dq \, dv \, |q|^2 \, \psi\big( \tfrac{|q|}{R} \big)^2 \, \mhus_{N,t}^{(1)}(q,v) \nonumber \\
    &= \hbar^{3/2} \int_{\mathbb{R}^3} dx \int_{\mathbb{R}^3} dq \, \nabla_q \Big( |q|^2 \, \psi\big( \tfrac{|q|}{R} \big)^2 \Big) f^2\big( \tfrac{x-q}{\sqrt{\hbar}} \big) \cdot \operatorname{Re} \big\langle a_x \Psi_{N,t}, \mpi_x a_x \Psi_{N,t} \big\rangle \nonumber \\
    &\le \C \Big( \hbar^{3/2} \int_{\{ |q| \ge R/2 \}} dq \int_{\mathbb{R}^3} dx \, f^2\big( \tfrac{x-q}{\sqrt{\hbar}} \big) \, \| \mpi_x a_x \Psi_{N,t} \|^2 \Big)^{1/2} \nonumber \\
    &\qquad \times \Big( \hbar^{3/2} \int_{\mathbb{R}^3} dq \, |q|^2 \, \psi\big( \tfrac{|q|}{R} \big)^2 \int_{\mathbb{R}^3} dx \, f^2\big( \tfrac{x-q}{\sqrt{\hbar}} \big) \, \| a_x \Psi_{N,t} \|^2 \Big)^{1/2} \nonumber \\
    &\le \C \Big( \hbar^{3/2} \int_{\{ |q| \ge R/2 \}} dq \int_{\mathbb{R}^3} dx \, f^2\big( \tfrac{x-q}{\sqrt{\hbar}} \big) \, \| \mpi_x a_x \Psi_{N,t} \|^2 \Big)^{1/2}
        \Big( \iint_{\mathbb{R}^6} dq \, dv \, |q|^2 \, \mhus_{N,t}^{(1)}(q,v) \Big)^{1/2} .
\end{align}
The first three steps of \eqref{eq:proof_kin_husimi_derivation1} and
$\|a - b - c\|^2 \ge \tfrac12 \|a\|^2 - 2 \|b\|^2 - 2 \|c\|^2$ give:
\begin{align}
\label{eq:spatial_tail_kinetic_localized}
    &\iint_{\{ |q| \ge R/2 \}} dq \, dv \, |v|^2 \, \mhus_{N,t}^{(1)}(q,v) \nonumber \\
    &= \hbar^{3/2} \int_{\{ |q| \ge R/2 \}} dq \int_{\mathbb{R}^3} dx \, \Big\| f\big(\tfrac{x-q}{\sqrt{\hbar}}\big) \mpi_x a_x \Psi_{N,t}
    - \big( A(x) - A(q) \big) f\big(\tfrac{x-q}{\sqrt{\hbar}}\big) a_x \Psi_{N,t} \nonumber \\
    &\qquad \qquad \qquad \qquad \qquad \qquad \qquad - i \sqrt{\hbar} \, (\nabla f)\big(\tfrac{x-q}{\sqrt{\hbar}}\big) a_x \Psi_{N,t} \Big\|^2 \nonumber \\
    &\ge \hbar^{3/2} \int_{\{ |q| \ge R/2 \}} dq \int_{\mathbb{R}^3} dx \, \bigg( \tfrac12 f^2\big(\tfrac{x-q}{\sqrt{\hbar}}\big)
    \| \mpi_x a_x \Psi_{N,t} \|^2 \nonumber \\
    &\qquad\qquad \qquad \qquad \qquad \qquad \qquad  - 2 \Big( \big| A(x) - A(q) \big|^2 f^2\big(\tfrac{x-q}{\sqrt{\hbar}}\big)
    + \hbar \, \big| (\nabla f)\big(\tfrac{x-q}{\sqrt{\hbar}}\big) \big|^2 \Big) \| a_x \Psi_{N,t} \|^2 \bigg) \nonumber \\
    &\ge \tfrac12 \hbar^{3/2} \int_{\{ |q| \ge R/2 \}} dq \int_{\mathbb{R}^3} dx \,
    f^2\big(\tfrac{x-q}{\sqrt{\hbar}}\big) \| \mpi_x a_x \Psi_{N,t} \|^2 - 2 \hbar \Big( \|\nabla A\|_{L^\infty(\mathbb{R}^3)}^2 \big\| |\cdot| f \big\|_{L^2(\mathbb{R}^3)}^2
    + \|\nabla f\|_{L^2(\mathbb{R}^3)}^2 \Big) .
\end{align}
\eqref{eq:spatial_tail_operator_derivative} and
\eqref{eq:spatial_tail_kinetic_localized} give
\begin{align}
\label{eq:spatial_tail_derivative}
    \tfrac{d}{dt} \iint_{\mathbb{R}^6} dq \, dv \, |q|^2 \, \psi\big( \tfrac{|q|}{R} \big)^2 \, \mhus_{N,t}^{(1)}(q,v) &\le \C \Big( \Big( \iint_{\{ |q| \ge R/2 \}} dq \, dv \, |v|^2 \, \mhus_{N,t}^{(1)}(q,v) \Big)^{1/2} + \sqrt{\hbar} \Big) \nonumber \\
    &\qquad \times \Big( \iint_{\mathbb{R}^6} dq \, dv \, |q|^2 \, \mhus_{N,t}^{(1)}(q,v) \Big)^{1/2} \nonumber \\
    &\le \C (1 + t) \Big( \iint_{\{ |q| \ge R/2 \}} dq \, dv \, |v|^2 \, \mhus_{N,t}^{(1)}(q,v) \Big)^{1/2} + \C (1 + t) \sqrt{\hbar} ,
\end{align}
where we used \eqref{eq:prop_k_moment_bound}. Integration over $[0,t]$ gives \eqref{eq:spatial_tail_bound}.
\end{proof}

The proof of Theorem \ref{thm:w2_energy_criterion_noV} combines Theorem \ref{thm:magnetic_vlasov_limit} with Lemmas
\ref{lemma:energy_decomposition_noV} and \ref{lem:spatial_tail_propagation}.

\begin{proof}[Proof of Theorem \ref{thm:w2_energy_criterion_noV}]
From the results in
Proposition \ref{prop:finite_moments_k} and \eqref{eq:vlasov_support}, we obtain that
$\mhus_{N,t}^{(k)}$, $\mathbf{f}_t^{\otimes k} \in \mathcal{P}_2(\mathbb{R}^{6k})$. By \cite[Thm.~6.9]{villani2008optimal}, convergence in $\mathbb{W}_2$ on $\mathcal{P}_2(\mathbb{R}^{6k})$ is narrow convergence
together with the convergence of the second moment, property (i) of \cite[Def.~6.8]{villani2008optimal}.
Theorem \ref{thm:magnetic_vlasov_limit}, tested with $\varphi \in C_b \subset L^\infty$, gives the narrow convergence.
By the symmetry and consistency \eqref{eq:normalized_husimi_measure} of $\mhus_{N,t}^{(k)}$ and of
$\mathbf{f}_t^{\otimes k}$, Definition \ref{def:admissible_class} \ref{item:admissible_a}, it suffices to prove the
case $k = 1$:
\begin{align}
\label{eq:criterion_one_particle_moment_noV}
    \lim_{N \to \infty} \iint_{\mathbb{R}^6} dq \, dv \, \Big( |q|^2 + |v|^2 \Big)
    \mhus_{N,t}^{(1)}(q,v)
    &= \iint_{\mathbb{R}^6} dq \, dv \, \Big( |q|^2 + |v|^2 \Big) \mathbf{f}_t(q,v) .
\end{align}

\emph{Truncation and split.} For $t \in [0,T]$ let $R \ge \max\{ R(T), R(T)^2 \}$, with $R(T)$ of \eqref{eq:vlasov_support},
so that $|q| \le R$ and $|v|^2 \le R$ on $\operatorname{supp}\mathbf{f}_t$. Then
\begin{align}
\label{eq:criterion_joint_split_noV}
    & \Big| \iint_{\mathbb{R}^6} dq \, dv \, \big( |q|^2 + |v|^2 \big) \big( \mhus_{N,t}^{(1)} - \mathbf{f}_t \big) \Big| \nonumber \\
    &\le \Big| \iint_{\mathbb{R}^6} dq \, dv \, \big( \min\{ |q|^2, R^2 \} + \min\{ |v|^2, R \} \big) \big( \mhus_{N,t}^{(1)} - \mathbf{f}_t \big) \Big|
        + \iint_{\{ |q| \ge R \}} dq \, dv \, |q|^2 \, \mhus_{N,t}^{(1)} \nonumber \\
    &\quad + 2 \iint_{\mathbb{R}^6} dq \, dv \, \Big( \tfrac12 |v|^2 + \Vext + c_1 - \min\big\{ \tfrac12 |v|^2 + \Vext + c_1, \, \tfrac{R}{2} \big\} \Big) \, \mhus_{N,t}^{(1)}.
\end{align}
The bound follows from
$|q|^2 - \min\{ |q|^2, R^2 \} \le |q|^2 \One_{\{ |q| \ge R \}}$ and the monotonicity of $x \mapsto x - \min\{ x, R \}$ at
$|v|^2 \le |v|^2 + 2 ( \Vext + c_1 )$, Assumption \ref{assump:A2}. The narrow convergence $\mhus_{N,t}^{(1)} \rightharpoonup \mathbf{f}_t$ sends the first term on the right of \eqref{eq:criterion_joint_split_noV} to $0$ as $N \to \infty$. It remains to show that the $\limsup_{N \to \infty}$ of the spatial $|q|^2$ tail and of the energy tail tends to $0$ as $R \to \infty$.

\emph{The energy tail.} Let $\tau \in [0,T]$. By \eqref{eq:energy_husimi_decomposition_noV}, Lemma
\ref{lemma:energy_decomposition_noV} at $t = \tau$ and the energy conservation \eqref{eq:energy_conservation},
\begin{align}
\label{eq:criterion_energy_remainder_noV}
    & \iint_{\mathbb{R}^6} dq \, dv \, \Big( \tfrac12 |v|^2 + \Vext + c_1 - \min\big\{ \tfrac12 |v|^2 + \Vext + c_1, \, \tfrac{R}{2} \big\} \Big) \, \mhus_{N,\tau}^{(1)} \nonumber \\
    &= \mathcal{V}_N(\tau) + c_1 - \iint_{\mathbb{R}^6} dq \, dv \, \min\big\{ \tfrac12 |v|^2 + \Vext + c_1, \, \tfrac{R}{2} \big\} \, \mhus_{N,\tau}^{(1)} \nonumber \\
    &\le \tfrac{1}{N} \big\langle \Psi_{N,0}, \mathcal{H} \Psi_{N,0} \big\rangle + c_1 - \mathcal{W}_N(\tau) + \C (1 + T) \sqrt{\hbar}
        - \iint_{\mathbb{R}^6} dq \, dv \, \min\big\{ \tfrac12 |v|^2 + \Vext + c_1, \, \tfrac{R}{2} \big\} \, \mhus_{N,\tau}^{(1)} .
\end{align}
As $\hbar = N^{-1/3} \to 0$, the energy inequality \eqref{eq:criterion_energy_noV} of Assumption
\ref{assump:I4} bounds the initial energy, and the integral of the bounded continuous truncation converges by narrow
convergence. Theorem \ref{thm:magnetic_vlasov_limit} at $k = 2$ gives
\begin{align}
\label{eq:criterion_interaction_limit_noV}
    \lim_{N \to \infty} \mathcal{W}_N(\tau)
    &= \tfrac{1}{2} \int_{\mathbb{R}^{12}} (dq \, dv)^{\otimes 2} \, W(q_1 - q_2) \,
    \mathbf{f}_\tau(q_1,v_1) \mathbf{f}_\tau(q_2,v_2) \; \eqqcolon \; \mathcal{W}[\mathbf{f}_\tau] .
\end{align}
The classical energy of \eqref{eq:magnetic_vlasov_intro} is conserved, $\mathcal{E}[\mathbf{f}_\tau] = \mathcal{E}[\mathbf{f}_0]$.
With \eqref{eq:classical_energy_definition}, \eqref{eq:criterion_energy_remainder_noV} gives, for every $R > 0$, the bound on the energy tail:
\begin{align}
\label{eq:criterion_energy_remainder_limit_noV}
    &\limsup_{N \to \infty} \iint_{\mathbb{R}^6} dq \, dv \,
    \Big( \tfrac12 |v|^2 + \Vext + c_1 - \min\big\{ \tfrac12 |v|^2 + \Vext + c_1, \, \tfrac{R}{2} \big\} \Big)
    \, \mhus_{N,\tau}^{(1)} \nonumber \\
    &\le \mathcal{E}[\mathbf{f}_0] + c_1 - \mathcal{W}[\mathbf{f}_\tau]
    - \iint_{\mathbb{R}^6} dq \, dv \,
    \min\big\{ \tfrac12 |v|^2 + \Vext + c_1, \, \tfrac{R}{2} \big\} \, \mathbf{f}_\tau \nonumber \\
    &= \iint_{\mathbb{R}^6} dq \, dv \,
    \Big( \tfrac12 |v|^2 + \Vext + c_1 - \min\big\{ \tfrac12 |v|^2 + \Vext + c_1, \, \tfrac{R}{2} \big\} \Big)
    \, \mathbf{f}_\tau .
\end{align}

\emph{The spatial tail.} We apply Lemma \ref{lem:spatial_tail_propagation} to the spatial $|q|^2$ term of
\eqref{eq:criterion_joint_split_noV} and bound the right side of \eqref{eq:spatial_tail_bound}, starting with the first
term, the $|q|^2$ term at time $0$.
On $\{ |q| \ge \tfrac{R}{2} \}$ we have $|q|^2 \le 2 \big( |q|^2 - \min\{ |q|^2, \tfrac{R^2}{8} \} \big)$, and
Assumptions \ref{assump:I2} and \ref{assump:I4} give
\begin{align}
\label{eq:criterion_tail_initial_noV}
    \limsup_{N \to \infty} \iint_{\{ |q| \ge \frac{R}{2} \}} dq \, dv \, |q|^2 \, \mhus_{N,0}^{(1)}
    &\le 2 \iint_{\mathbb{R}^6} dq \, dv \, \big( |q|^2 - \min\{ |q|^2, \tfrac{R^2}{8} \} \big) \, \mathbf{f}_0 .
\end{align}
For the second term of the right side of \eqref{eq:spatial_tail_bound}, the time integral of the $|v|^2$ term, we split $|v|^2$ at $R$ for $\tau \in [0,T]$ as
in \eqref{eq:criterion_joint_split_noV}:
\begin{align}
\label{eq:criterion_velocity_tail_noV}
    & \iint_{\{ |q| \ge R/2 \}} dq \, dv \, |v|^2 \, \mhus_{N,\tau}^{(1)} \nonumber \\
    &\le \iint_{\mathbb{R}^6} dq \, dv \, \big( |v|^2 - \min\{ |v|^2, R \} \big) \, \mhus_{N,\tau}^{(1)} + R \iint_{\{ |q| \ge R/2 \}} dq \, dv \, \mhus_{N,\tau}^{(1)} \nonumber \\
    &\le 2 \iint_{\mathbb{R}^6} dq \, dv \, \Big( \tfrac12 |v|^2 + \Vext + c_1 - \min\big\{ \tfrac12 |v|^2 + \Vext + c_1, \, \tfrac{R}{2} \big\} \Big) \mhus_{N,\tau}^{(1)}
        + \tfrac{4}{R} \iint_{\mathbb{R}^6} dq \, dv \, |q|^2 \mhus_{N,\tau}^{(1)} \nonumber \\
    &\le 2 \iint_{\mathbb{R}^6} dq \, dv \, \Big( \tfrac12 |v|^2 + \Vext + c_1 - \min\big\{ \tfrac12 |v|^2 + \Vext + c_1, \, \tfrac{R}{2} \big\} \Big) \, \mhus_{N,\tau}^{(1)} + \tfrac{4 \C (1 + T)^2}{R} .
\end{align}
By Proposition \ref{prop:finite_moments_k} at $k = 1$, the square root of
$\iint_{\{ |q| \ge R/2 \}} dq \, dv \, |v|^2 \, \mhus_{N,\tau}^{(1)}$ is bounded by $\C^{\frac12} (1 + T)$, so the reverse
Fatou lemma and \eqref{eq:criterion_velocity_tail_noV}, \eqref{eq:criterion_energy_remainder_limit_noV} give
\begin{align}
\label{eq:criterion_velocity_tail_integrated_noV}
    & \limsup_{N \to \infty} \int_0^t d\tau \, \bigg( \iint_{\{ |q| \ge R/2 \}} dq \, dv \, |v|^2 \, \mhus_{N,\tau}^{(1)} \bigg)^{\frac12} \nonumber \\
    &\le \int_0^t d\tau \, \limsup_{N \to \infty} \bigg( \iint_{\{ |q| \ge R/2 \}} dq \, dv \, |v|^2 \, \mhus_{N,\tau}^{(1)} \bigg)^{\frac12} \nonumber \\
    &\le \int_0^t d\tau \, \bigg[ 2 \iint_{\mathbb{R}^6} dq \, dv \, \Big( \tfrac12 |v|^2 + \Vext + c_1 - \min\big\{ \tfrac12 |v|^2 + \Vext + c_1, \, \tfrac{R}{2} \big\} \Big) \, \mathbf{f}_\tau
        + \tfrac{4 \C (1 + T)^2}{R} \bigg]^{\frac12} .
\end{align}
With \eqref{eq:spatial_tail_bound}, \eqref{eq:criterion_tail_initial_noV},
\eqref{eq:criterion_velocity_tail_integrated_noV} and $\hbar \to 0$,
\begin{align}
\label{eq:criterion_tail_limsup_noV}
    &\limsup_{N \to \infty} \iint_{\{ |q| \ge R \}} dq \, dv \, |q|^2 \, \mhus_{N,t}^{(1)} \nonumber \\
    &\le 2 \iint_{\mathbb{R}^6} dq \, dv \, \big( |q|^2 - \min\{ |q|^2, \tfrac{R^2}{8} \} \big) \, \mathbf{f}_0 \nonumber\\
       &\quad + \C (1 + t) \int_0^t d\tau \, \bigg[ 2 \iint_{\mathbb{R}^6} dq \, dv \, \Big( \tfrac12 |v|^2 + \Vext + c_1 - \min\big\{ \tfrac12 |v|^2 + \Vext + c_1, \, \tfrac{R}{2} \big\} \Big) \, \mathbf{f}_\tau  + \tfrac{4 \C (1 + T)^2}{R} \bigg]^{\frac12} .
\end{align}
Since $\operatorname{supp}\mathbf{f}_\tau \subset \{ |q|^2 + |v|^2 \le R(T)^2 \}$ for $\tau \in [0,T]$, the integrals
against $\mathbf{f}_0$ and $\mathbf{f}_\tau$ vanish for $R$ large, and the right side tends to $0$ as $R \to \infty$.

\emph{Conclusion.} By the first term of \eqref{eq:criterion_joint_split_noV}, \eqref{eq:criterion_energy_remainder_limit_noV} at $\tau = t$ and
\eqref{eq:criterion_tail_limsup_noV}, the $\limsup_{N \to \infty}$ of the left side of
\eqref{eq:criterion_joint_split_noV} is at most a sum of terms that tend to $0$ as $R \to \infty$. The left side
does not depend on $R$, and $R \to \infty$ gives \eqref{eq:criterion_one_particle_moment_noV}. This is
\eqref{eq:criterion_transport_noV}.
\end{proof}

\paragraph{Energy convergence.} The narrow convergence \eqref{eq:initial_condition_tensor} at $k = 1$, tested with
\begin{align}
\label{eq:criterion_energy_test_function}
    \min\{ \tfrac12 |v|^2 + \Vext + c_1, \alpha \} &\in C_b(\mathbb{R}^6) \quad \text{and} \quad \alpha \to \infty,
\end{align}
gives
\begin{align}
\label{eq:criterion_energy_liminf_noV}
    \liminf_{N \to \infty} \iint_{\mathbb{R}^6} dq \, dv \, \big( \tfrac12 |v|^2 + \Vext + c_1 \big) \, \mhus_{N,0}^{(1)}
    &\ge \iint_{\mathbb{R}^6} dq \, dv \, \big( \tfrac12 |v|^2 + \Vext + c_1 \big) \, \mathbf{f}_0 .
\end{align}
With Lemma \ref{lemma:energy_decomposition_noV}, \eqref{eq:criterion_interaction_limit_noV} at $\tau = 0$ and
\eqref{eq:criterion_energy_liminf_noV} we obtain the first inequality below, Assumption \ref{assump:I4} the second:
\begin{align}
\label{eq:criterion_energy_lower_bound_noV}
    \liminf_{N \to \infty} \tfrac{1}{N} \big\langle \Psi_{N,0}, \mathcal{H} \Psi_{N,0} \big\rangle
    &\ge \iint_{\mathbb{R}^6} dq \, dv \, \Big( \tfrac{1}{2}|v|^2 + \Vext \Big) \mathbf{f}_0 + \mathcal{W}[\mathbf{f}_0]
    = \mathcal{E}[\mathbf{f}_0] \ge \limsup_{N \to \infty} \tfrac{1}{N} \big\langle \Psi_{N,0}, \mathcal{H} \Psi_{N,0} \big\rangle ,
\end{align}
so the energy inequality \eqref{eq:criterion_energy_noV} is equivalent to
$\lim_{N \to \infty} \frac{1}{N} \langle \Psi_{N,0}, \mathcal{H} \Psi_{N,0} \rangle = \mathcal{E}[\mathbf{f}_0]$.

\section{Unified Bounds of the Remainder Components} \label{sec:unified_bounds}

This section bounds each remainder component $\mathcal{I}_{\star, j}^{(k)}$ of \eqref{eq:unified_strong_weak_remainders},
one subsection per component. Throughout this section, for $1 \le m \le k+1$, we abbreviate
\begin{align}
\label{eq:section8_notation}
    f^x_m &\coloneqq f\Big( \frac{x_m - q_m}{\sqrt{\hbar}} \Big) , \qquad
    f^y_m \coloneqq f\Big( \frac{y_m - q_m}{\sqrt{\hbar}} \Big) , \qquad
    e^{\pm}_\hbar \coloneqq e^{\pm \frac{i}{\hbar} \sum_{n=1}^k (v_n - A(q_n)) \cdot (x_n - y_n)} ,
\end{align}
with $(\nabla f)^x_m$ and its $y$-counterpart defined in the same way. Since the phases $e^{\pm \frac{i}{\hbar} \sum_n (v_n - A(q_n)) \cdot (y_n - x_n)}$
are $e^{\mp}_\hbar$, the symmetric sums $\frac{1}{2} \sum_\pm$ below carry $e^{\pm}_\hbar$.

\subsection{Bounding the Dispersion Remainder Integral \texorpdfstring{$\mathcal{I}_{\text{dis}, j}^{(k)}$}{I-dis-j-k}}
\label{subsubsec:bound_I_dis_j_k}

The remainder $\mathcal{I}_{\text{dis}, j}^{(k)}$ \eqref{eq:unified_strong_weak_remainders} is defined with $\mathcal{R}_{\text{dis}, j}^{(k)}$ \eqref{eq:prop_R_disp}. Applying \eqref{eq:state_level_chain_rule}, expanding the imaginary projection via exchange of the spatial variables $x \leftrightarrow y$, and integrating by parts transfers the divergence-free vector field $\nabla_{q_j} + (\nabla A(q_j))^T \nabla_{v_j}$ onto the test observable $\Phi$:
\begin{align}
\label{eq:I_dis_derivation}
    \mathcal{I}_{\text{dis}, j}^{(k)} &= - \tfrac{c_{N,k}}{(2\pi)^{3k}} \hbar \int_{\mathbb{R}^{6k}} \dZk \int_{\mathbb{R}^{6k}} \dXkYk \, \Phi(Z_k)
        \operatorname{Im} \Bigg[ \Big( \big( \nabla_{q_j} + (\nabla A(q_j))^T \nabla_{v_j} \big)^2 \overline{\mathfrak{f}_j(y_j)} \Big) \mathfrak{f}_j(x_j) \nonumber \\
    &\quad\qquad\qquad\qquad\qquad\qquad\qquad\qquad\qquad\qquad\qquad \qquad\times \prod_{m \neq j}^k \mathfrak{f}_m(x_m) \overline{\mathfrak{f}_m(y_m)} \, \dmat{k}(X;Y) \Bigg] \nonumber \\
    &= \tfrac{c_{N,k}}{(2\pi)^{3k}} \frac{i\hbar}{2} \int_{\mathbb{R}^{6k}} \dZk \int_{\mathbb{R}^{6k}} \dXkYk \prod_{m \neq j}^k \mathfrak{f}_m(x_m) \overline{\mathfrak{f}_m(y_m)} \, \dmat{k}(X;Y) \nonumber \\
    &\quad \times \Phi(Z_k) \big( \nabla_{q_j} + (\nabla A(q_j))^T \nabla_{v_j} \big) \cdot \Big[ \mathfrak{f}_j(x_j) \big( \nabla_{q_j} + (\nabla A(q_j))^T \nabla_{v_j} \big) \overline{\mathfrak{f}_j(y_j)} \nonumber \\
    &\qquad \qquad\qquad\qquad\qquad\qquad\qquad\qquad\qquad\qquad\qquad \qquad- \overline{\mathfrak{f}_j(y_j)} \big( \nabla_{q_j} + (\nabla A(q_j))^T \nabla_{v_j} \big) \mathfrak{f}_j(x_j) \Big] \nonumber \\
    &= \tfrac{c_{N,k}}{(2\pi)^{3k}} \frac{i\hbar}{2} \int_{\mathbb{R}^{6k}} \dZk \int_{\mathbb{R}^{6k}} \dXkYk \prod_{m \neq j}^k \mathfrak{f}_m(x_m) \overline{\mathfrak{f}_m(y_m)} \, \dmat{k}(X;Y) \nonumber \\
    &\quad \times \Big[ \big( \nabla_{q_j} + (\nabla A(q_j))^T \nabla_{v_j} \big) \Phi(Z_k) \Big] \cdot \Big[ \overline{\mathfrak{f}_j(y_j)} \big( \nabla_{q_j} + (\nabla A(q_j))^T \nabla_{v_j} \big) \mathfrak{f}_j(x_j) \nonumber \\
    &\qquad \qquad\qquad\qquad\qquad\qquad\qquad\qquad\qquad\qquad\qquad \qquad- \mathfrak{f}_j(x_j) \big( \nabla_{q_j} + (\nabla A(q_j))^T \nabla_{v_j} \big) \overline{\mathfrak{f}_j(y_j)} \Big] \nonumber \\
    &= \tfrac{c_{N,k}}{(2\pi)^{3k}} \frac{i}{2} \hbar^{1/2 - 3k/2} \int_{\mathbb{R}^{6k}} \dZk \int_{\mathbb{R}^{6k}} \dXkYk \prod_{m \neq j}^k f^x_m f^y_m \, \dmat{k}(X;Y) \nonumber \\
    &\quad \times e^{+}_\hbar \Big[ \big( \nabla_{q_j} + (\nabla A(q_j))^T \nabla_{v_j} \big) \Phi(Z_k) \Big] \cdot \big( f^x_j (\nabla f)^y_j - (\nabla f)^x_j f^y_j \big) .
\end{align}
The Cauchy--Schwarz and Young inequalities decouple the trace. For the first component, this yields:
\begin{align}
\label{eq:am_gm_imbedded_disp_1}
    \Big| \Big\langle f^x_j \prod_{m \neq j}^k f^x_m \, a_X \Psi_{N,t} , \, (\nabla f)^y_j \prod_{m \neq j}^k f^y_m \, a_Y \Psi_{N,t} \Big\rangle \Big| &\le \frac{1}{2} |f^x_j|^2 \prod_{m \neq j}^k |f^x_m|^2 \, \| a_X \Psi_{N,t} \|^2 \nonumber \\
    &\quad + \frac{1}{2} |(\nabla f)^y_j|^2 \prod_{m \neq j}^k |f^y_m|^2 \, \| a_Y \Psi_{N,t} \|^2 .
\end{align}
The second component yields the same bound with the profile gradients reversed, and the map $X \leftrightarrow Y$ on the $Y$-dependent halves generates the symmetric phase addition $\frac{1}{2}\sum_\pm$:
\begin{align}
\label{eq:I_dis_fully_expanded}
    \big| \mathcal{I}_{\text{dis}, j}^{(k)} \big| &\le \tfrac{c_{N,k}}{(2\pi)^{3k}} \frac{1}{2} \hbar^{1/2 - 3k/2} \int_{\mathbb{R}^{3k}} (dq)^{\otimes k} \int_{\mathbb{R}^{6k}} \dXkYk \, \| a_X \Psi_{N,t} \|^2
        \Big( |f^x_j|^2 + |(\nabla f)^x_j|^2 \Big) \prod_{m \neq j}^k |f^x_m|^2 \nonumber \\
    &\quad \times \frac{1}{2} \sum_{\pm} \Big| \int_{\mathbb{R}^{3k}} (dv)^{\otimes k} \big( \nabla_{q_j} + (\nabla A(q_j))^T \nabla_{v_j} \big) \Phi(Z_k) \, e^{\pm}_\hbar \Big|_\infty .
\end{align}

\subsection{Bounding the Force Deviation Remainder Integral \texorpdfstring{$\mathcal{I}_{\text{for}, j}^{(k)}$}{I-for-j-k}}
\label{subsubsec:bound_I_for_j_k}

The remainder $\mathcal{I}_{\text{for}, j}^{(k)}$ \eqref{eq:unified_strong_weak_remainders} is defined with $\mathcal{R}_{\text{for}, j}^{(k)}$ \eqref{eq:prop_R_for} and $ \mathrm M_A$ \eqref{eq:MAaveragedgradient_k}. Applying the Fundamental Theorem of Calculus to extract $(y_j - x_j)$, we substitute $(y_j - x_j) e^{\frac{i}{\hbar} v_j \cdot (x_j - y_j)} = i\hbar \nabla_{v_j} e^{\frac{i}{\hbar} v_j \cdot (x_j - y_j)}$. Velocity integration by parts transfers the derivative:
\begin{align}
\label{eq:I_for_exact_FTC}
    \mathcal{I}_{\text{for}, j}^{(k)} &= \tfrac{c_{N,k}}{(2\pi)^{3k}} \frac{2}{\hbar} \int_{\mathbb{R}^{6k}} \dZk \int_{\mathbb{R}^{6k}} \dXkYk \, \Phi(Z_k)
        \operatorname{Im} \Bigg[ v_j \cdot \Big( \mathrm M_A(y_j, q_j)(y_j - q_j) \Big) \mathfrak{f}_j(x_j) \overline{\mathfrak{f}_j(y_j)} \nonumber \\
    &\quad \times \prod_{m \neq j}^k \mathfrak{f}_m(x_m) \overline{\mathfrak{f}_m(y_m)} \, \dmat{k}(X;Y) \Bigg] \nonumber \\
    &= \tfrac{c_{N,k}}{(2\pi)^{3k}} \frac{1}{i\hbar} \int_{\mathbb{R}^{6k}} \dZk \int_{\mathbb{R}^{6k}} \dXkYk \, \Phi(Z_k) \, \dmat{k}(X;Y) \nonumber \\
    &\quad \times \prod_{m = 1}^k \mathfrak{f}_m(x_m) \overline{\mathfrak{f}_m(y_m)} \, v_j \cdot \Big( \mathrm M_A(y_j, q_j)(y_j - q_j) - \mathrm M_A(x_j, q_j)(x_j - q_j) \Big) \nonumber \\
    &= \tfrac{c_{N,k}}{(2\pi)^{3k}} \hbar^{-3k/2} \int_{\mathbb{R}^{6k}} \dZk \int_{\mathbb{R}^{6k}} \dXkYk \prod_{m = 1}^k f^x_m f^y_m \, \dmat{k}(X;Y) \nonumber \\
    &\quad \times \Phi(Z_k) \Big[ \Big( \int_0^1 du \, \nabla_z \big( \mathrm M_A(z, q_j)(z - q_j) \big) \Big|_{z = (1-u) x_j + u y_j} \Big)^T v_j \Big] \cdot \nabla_{v_j} e^{+}_\hbar \nonumber \\
    &= - \tfrac{c_{N,k}}{(2\pi)^{3k}} \hbar^{-3k/2} \int_{\mathbb{R}^{6k}} \dZk \int_{\mathbb{R}^{6k}} \dXkYk \prod_{m = 1}^k f^x_m f^y_m \, e^{+}_\hbar \, \dmat{k}(X;Y) \nonumber \\
    &\quad \times \operatorname{Tr}_{\mathbb{C}^{3\times 3}} \Bigg( \Big( \int_0^1 du \, \nabla_z \big( \mathrm M_A(z, q_j)(z - q_j) \big) \Big|_{z = (1-u) x_j + u y_j} \Big)
        \big( \nabla_{v_j} \Phi(Z_k) \, v_j^T + \Phi(Z_k) \mathbf{I}_{3\times 3} \big) \Bigg) .
\end{align}
With $\nabla^2 A \in L^\infty(\mathbb{R}^3)$, \eqref{eq:MAaveragedgradient_k} yields:
\begin{align}
\label{eq:M_components_bound}
\| \mathrm M_A(z, q_j)\|_{\text{op}} &\le \tfrac{1}{2} \|\nabla^2 A\|_{L^\infty(\mathbb{R}^3)} |z - q_j|, \qquad
    \|\nabla_z \mathrm M_A(z, q_j)\|_{\text{op}} \le \tfrac{1}{2} \|\nabla^2 A\|_{L^\infty(\mathbb{R}^3)}.
\end{align}
This bounds the integral:
\begin{align}
\label{eq:FTC_Jacobian_bound}
    &\Big\| \int_0^1 du \, \nabla_z \big( \mathrm M_A(z, q_j)(z - q_j) \big) \Big|_{z = (1-u) x_j + u y_j} \Big\|_{\text{op}}
    \nonumber \\
    &\le \int_0^1 du \, \Big( \big\| \nabla_z \mathrm M_A(z, q_j) \big\|_{\text{op}} |z - q_j|
    + \big\| \mathrm M_A(z, q_j) \big\|_{\text{op}} \Big) \Big|_{z = (1-u) x_j + u y_j}
    \nonumber \\
    &\le \|\nabla^2 A\|_{L^\infty(\mathbb{R}^3)} \int_0^1 du \, \Big( (1-u)|x_j - q_j| + u|y_j - q_j| \Big)
    \nonumber \\
    &= \frac{1}{2} \|\nabla^2 A\|_{L^\infty(\mathbb{R}^3)} \Big( |x_j - q_j| + |y_j - q_j| \Big) .
\end{align}
The Cauchy--Schwarz and Young inequalities decouple the trace, and the map $X \leftrightarrow Y$ on the second half generates the symmetric phase addition $\frac{1}{2}\sum_\pm$:
\begin{align}
\label{eq:I_for_abs_partitioned}
    \big| \mathcal{I}_{\text{for}, j}^{(k)} \big| &\le \tfrac{c_{N,k}}{(2\pi)^{3k}} \frac{3\sqrt{3}}{2} \|\nabla^2 A\|_{L^\infty(\mathbb{R}^3)} \hbar^{1/2 - 3k/2} \int_{\mathbb{R}^{3k}} (dq)^{\otimes k} \int_{\mathbb{R}^{6k}} \dXkYk
        \| a_X \Psi_{N,t} \| \| a_Y \Psi_{N,t} \| \nonumber \\
    &\quad \times \Big( \Big| \frac{x_j - q_j}{\sqrt{\hbar}} \Big| + \Big| \frac{y_j - q_j}{\sqrt{\hbar}} \Big| \Big) |f^x_j| |f^y_j| \prod_{m \neq j}^k |f^x_m| |f^y_m|
        \Big| \int_{\mathbb{R}^{3k}} (dv)^{\otimes k} \big( \nabla_{v_j} \Phi(Z_k) \, v_j^T + \Phi(Z_k) \mathbf{I}_{3\times 3} \big) e^{+}_\hbar \Big|_\infty \nonumber \\
    &\le \tfrac{c_{N,k}}{(2\pi)^{3k}} \frac{3\sqrt{3}}{2} \|\nabla^2 A\|_{L^\infty(\mathbb{R}^3)} \hbar^{1/2 - 3k/2} \int_{\mathbb{R}^{3k}} (dq)^{\otimes k} \int_{\mathbb{R}^{6k}} \dXkYk
        \| a_X \Psi_{N,t} \|^2 \Big( 1 + \Big| \frac{x_j - q_j}{\sqrt{\hbar}} \Big|^2 \Big) |f^x_j|^2 \nonumber \\
    &\quad \times \prod_{m \neq j}^k |f^x_m|^2
        \frac{1}{2} \sum_{\pm} \Big| \int_{\mathbb{R}^{3k}} (dv)^{\otimes k} \big( \nabla_{v_j} \Phi(Z_k) \, v_j^T + \Phi(Z_k) \mathbf{I}_{3\times 3} \big) e^{\pm}_\hbar \Big|_\infty .
\end{align}

\subsection{Bounding the Gauge Remainder Integral \texorpdfstring{$\mathcal{I}_{\text{gau}, j}^{(k)}$}{I-gau-j-k}}
\label{subsubsec:bound_I_gau_j_k}

The remainder $\mathcal{I}_{\text{gau}, j}^{(k)}$ \eqref{eq:unified_strong_weak_remainders} is defined with $\mathcal{R}_{\text{gau}, j}^{(k)}$ \eqref{eq:prop_R_gau}. We resolve the real projection via the exchange $X \leftrightarrow Y$ to symmetrize the divergence term. Applying the global Lipschitz bound yields the semiclassical scaling. The Cauchy--Schwarz and Young inequalities decouple the trace, and the map $X \leftrightarrow Y$ on the second half generates the symmetric phase addition $\frac{1}{2}\sum_\pm$:
\begin{align}
\label{eq:I_gau_abs_bound}
     \big| \mathcal{I}_{\text{gau}, j}^{(k)} \big|    &\le \tfrac{c_{N,k}}{(2\pi)^{3k}} \frac{1}{2} \hbar^{-3k/2} \int_{\mathbb{R}^{3k}} (dq)^{\otimes k}
    \int_{\mathbb{R}^{6k}} \dXkYk \, \| a_X \Psi_{N,t} \| \| a_Y \Psi_{N,t} \|
    \nonumber \\
    &\quad \times \Big| (\nabla \cdot A)(y_j) - (\nabla \cdot A)(q_j) + (\nabla \cdot A)(x_j) - (\nabla \cdot A)(q_j) \Big|
    \nonumber \\
    &\quad \times |f^x_j| |f^y_j| \prod_{m \neq j}^k |f^x_m| |f^y_m|
    \, \Big| \int_{\mathbb{R}^{3k}} (dv)^{\otimes k} \Phi(Z_k) \, e^{+}_\hbar \Big| \nonumber
    \\
    &\le \tfrac{c_{N,k}}{(2\pi)^{3k}} \frac{1}{2} \|\nabla^2 A\|_{L^\infty(\mathbb{R}^3)} \hbar^{1/2 - 3k/2}
    \int_{\mathbb{R}^{3k}} (dq)^{\otimes k} \int_{\mathbb{R}^{6k}} \dXkYk \, \| a_X \Psi_{N,t} \|^2
    \nonumber \\
    &\quad \times \Big( 1 + \Big| \frac{x_j - q_j}{\sqrt{\hbar}} \Big|^2 \Big)
    |f^x_j|^2 \prod_{m \neq j}^k |f^x_m|^2 \, \frac{1}{2} \sum_{\pm}
    \Big| \int_{\mathbb{R}^{3k}} (dv)^{\otimes k} \Phi(Z_k) \, e^{\pm}_\hbar \Big| . 
\end{align}
\subsection{Bounding the Diamagnetic Remainder Integral \texorpdfstring{$\mathcal{I}_{\text{dia}, j}^{(k)}$}{I-dia-j-k}}
\label{subsubsec:bound_I_dia_j_k}

The remainder $\mathcal{I}_{\text{dia}, j}^{(k)}$ \eqref{eq:unified_strong_weak_remainders} is defined with $\mathcal{R}_{\text{dia}, j}^{(k)}$ \eqref{eq:prop_R_dia}. We resolve the imaginary projection via the exchange $X \leftrightarrow Y$ to symmetrize the quadratic difference. Applying the Lipschitz bound yields:
\begin{align}
\label{eq:dia_lipschitz_defect_abs}
    \Big| |A(y_j) - A(q_j)|^2 - |A(x_j) - A(q_j)|^2 \Big| &\le \sqrt{\hbar} \, \|\nabla A\|_{L^\infty(\mathbb{R}^3)}^2 |x_j - y_j| \Big( \Big| \frac{y_j - q_j}{\sqrt{\hbar}} \Big| + \Big| \frac{x_j - q_j}{\sqrt{\hbar}} \Big| \Big) .
\end{align}

The Cauchy--Schwarz and Young inequalities decouple the trace, with each amplitude paired against $\mathbf{1}$, and the map $X \leftrightarrow Y$ on the second half converts $\left| \frac{y_j-q_j}{\sqrt{\hbar}} \right|$ into $\left| \frac{x_j-q_j}{\sqrt{\hbar}} \right|^2$ and generates the symmetric phase addition $\frac{1}{2}\sum_\pm$:
\begin{align}
\label{eq:I_dia_abs_bound}
    \big| \mathcal{I}_{\text{dia}, j}^{(k)} \big| &\le \tfrac{c_{N,k}}{(2\pi)^{3k}} \frac{1}{2\hbar} \hbar^{-3k/2} \int_{\mathbb{R}^{3k}} (dq)^{\otimes k} \int_{\mathbb{R}^{6k}} \dXkYk \, \| a_X \Psi_{N,t} \| \| a_Y \Psi_{N,t} \| \nonumber \\
    &\quad \times \Big| |A(y_j) - A(q_j)|^2 - |A(x_j) - A(q_j)|^2 \Big|
        |f^x_j| |f^y_j| \prod_{m \neq j}^k |f^x_m| |f^y_m| \, \Big| \int_{\mathbb{R}^{3k}} (dv)^{\otimes k} \Phi(Z_k) \, e^{+}_\hbar \Big| \nonumber \\
    &\le \tfrac{c_{N,k}}{(2\pi)^{3k}} \frac{1}{2} \|\nabla A\|_{L^\infty(\mathbb{R}^3)}^2 \hbar^{-\frac{1}{2} - \frac{3k}{2}} \int_{\mathbb{R}^{3k}} (dq)^{\otimes k} \int_{\mathbb{R}^{6k}} \dXkYk
        \| a_X \Psi_{N,t} \| \| a_Y \Psi_{N,t} \| \, |x_j - y_j| \nonumber \\
    &\quad \times \Big( \Big| \frac{x_j - q_j}{\sqrt{\hbar}} \Big| + \Big| \frac{y_j - q_j}{\sqrt{\hbar}} \Big| \Big)
        |f^x_j| |f^y_j| \prod_{m \neq j}^k |f^x_m| |f^y_m| \, \Big| \int_{\mathbb{R}^{3k}} (dv)^{\otimes k} \Phi(Z_k) \, e^{+}_\hbar \Big| \nonumber \\
    &\le \tfrac{c_{N,k}}{(2\pi)^{3k}} \frac{1}{2} \|\nabla A\|_{L^\infty(\mathbb{R}^3)}^2 \hbar^{-\frac{1}{2} - \frac{3k}{2}} \int_{\mathbb{R}^{3k}} (dq)^{\otimes k} \int_{\mathbb{R}^{6k}} \dXkYk \nonumber \\
    &\quad \times \| a_X \Psi_{N,t} \|^2 \, |x_j - y_j| \Big( 1 + \Big| \frac{x_j - q_j}{\sqrt{\hbar}} \Big|^2 \Big) \prod_{m = 1}^k |f^x_m|^2
        \frac{1}{2} \sum_{\pm} \Big| \int_{\mathbb{R}^{3k}} (dv)^{\otimes k} \Phi(Z_k) \, e^{\pm}_\hbar \Big| .
\end{align}
The spatial difference $|x_j - y_j|$ supplies the active weight $\beta_{\text{dia}} = 1$.

\subsection{Bounding the Divergence Remainder Integral \texorpdfstring{$\mathcal{I}_{\text{div}, j}^{(k)}$}{I-div-j-k}}
\label{subsubsec:bound_I_div_j_k}
The remainder $\mathcal{I}_{\text{div}, j}^{(k)}$ \eqref{eq:unified_strong_weak_remainders} is defined with $\mathcal{R}_{\text{div}, j}^{(k)}$ \eqref{eq:prop_R_FS}. Expanding the phase-space operators \eqref{eq:phase_space_number_operators}, resolving the real projection and interchanging $x_n \leftrightarrow y_n$ leads to:
\begin{align}
\label{eq:I_div_j_raw_to_symmetrized}
    \mathcal{I}_{\text{div}, j}^{(k)} &= \tfrac{c_{N,k}}{(2\pi)^{3k}} \int_{\mathbb{R}^{6k}} \dZk \, \Phi(Z_k)
        \Bigg( - 2 \operatorname{Re} \Big\langle \Psi_{N,t}, \mathfrak{A}_k^\dagger[\mathfrak{f}] \, \mathfrak{A}_{k,j} \Big[ \Big( \big(A(\cdot) - A(q_j)\big) \cdot \nabla_{q_j}\big|_{p_j} \Big) \mathfrak{f}_j \Big] \Psi_{N,t} \Big\rangle \nonumber \\
    &\qquad + \big(\nabla \cdot A(q_j)\big) \operatorname{Re} \big\langle \Psi_{N,t}, \mathfrak{A}_k^\dagger[\mathfrak{f}] \, \mathfrak{A}_{k,j} [ \mathfrak{f}_j ] \Psi_{N,t} \big\rangle \Bigg) \nonumber \\
    &= \tfrac{c_{N,k}}{(2\pi)^{3k}} \hbar^{-3k/2} \int_{\mathbb{R}^{6k}} \dZk \int_{\mathbb{R}^{6k}} \dXkYk \, \Phi(Z_k) \prod_{m \neq j}^k f^x_m f^y_m \nonumber \\
    &\quad \times \Bigg( \frac{2}{\sqrt{\hbar}} \operatorname{Re} \Big[ e^{+}_\hbar \, \dmat{k}(X;Y) \, f^x_j \big(A(y_j) - A(q_j)\big) \cdot (\nabla f)^y_j \Big] \nonumber \\
    &\qquad + \big(\nabla \cdot A(q_j)\big) \operatorname{Re} \Big[ e^{+}_\hbar \, \dmat{k}(X;Y) \, f^x_j f^y_j \Big] \Bigg) \nonumber \\
    &= \tfrac{c_{N,k}}{(2\pi)^{3k}} \hbar^{-3k/2 - 1/2} \int_{\mathbb{R}^{6k}} \dZk \int_{\mathbb{R}^{6k}} \dXkYk \, \Phi(Z_k) \prod_{m \neq j}^k f^x_m f^y_m
        e^{+}_\hbar \, \dmat{k}(X;Y) \nonumber \\
    &\quad \times \Big[ f^x_j \big(A(y_j) - A(q_j)\big) \cdot (\nabla f)^y_j
        + \big(A(x_j) - A(q_j)\big) \cdot (\nabla f)^x_j \, f^y_j + \sqrt{\hbar} \big(\nabla \cdot A(q_j)\big) f^x_j f^y_j \Big] .
\end{align}
Expanding the magnetic difference via the averaged gradient matrix \eqref{eq:MAaveragedgradient_k} yields:
\begin{align}
\label{eq:magnetic_difference_expansion}
A(z) - A(q_j) &= \nabla A(q_j)(z-q_j) +  \mathrm M_A(z, q_j)(z-q_j).
\end{align}

\paragraph{Cancellation of the divergence term.}

The term $\nabla \cdot A(q_j)$ in the last row of \eqref{eq:I_div_j_raw_to_symmetrized} is of order one and cancels against the projection onto $f$ of the linear part of \eqref{eq:magnetic_difference_expansion}. For fixed
$q_j \in \mathbb{R}^3$ and the scaled variable $z \in \mathbb{R}^3$, the residual is
\begin{align}
\label{eq:orthprojection_refined}
    \Theta_j(z) &= \big( \nabla A(q_j) z \big) \cdot \nabla f(z) + \frac{1}{2}\big(\nabla \cdot A(q_j)\big) f(z) .
\end{align}
Substituting the algebraic identity $f \, \nabla f = \tfrac{1}{2} \nabla \big( f^2 \big)$ and applying integration by parts, with $z f(z)^2 \in W^{1,1}(\mathbb{R}^3; \mathbb{R}^3)$ by Assumption \ref{assump:A4}, yields the orthogonality condition:
\begin{align}
\label{eq:correct_orthogonality_v2}
\langle f, \Theta_j \rangle_{L^2(\mathbb{R}^3)} &= \int_{\mathbb{R}^3} dz \, f(z) \big( \nabla A(q_j) z \big) \cdot \nabla f(z) + \frac{1}{2}\big(\nabla \cdot A(q_j)\big) \int_{\mathbb{R}^3} dz \, f(z)^2 \nonumber \\
    &= \frac{1}{2} \int_{\mathbb{R}^3} dz \, \big( \nabla A(q_j) z \big) \cdot \nabla \big( f(z)^2 \big) + \frac{1}{2} \big(\nabla \cdot A(q_j)\big) \nonumber \\
    &= -\frac{1}{2} \int_{\mathbb{R}^3} dz \, \big( \nabla_z \cdot (\nabla A(q_j) z) \big) f(z)^2 + \frac{1}{2} \big(\nabla \cdot A(q_j)\big) \nonumber \\
    &= -\frac{1}{2} \big(\nabla \cdot A(q_j)\big) \int_{\mathbb{R}^3} dz \, f(z)^2 + \frac{1}{2} \big(\nabla \cdot A(q_j)\big) \nonumber \\
    &= 0.
\end{align}
\vspace{1ex}
\noindent
\begin{minipage}[c]{0.62\textwidth}
Isolating the linear gradient profile from definition \eqref{eq:orthprojection_refined} yields the orthogonal splitting:
\begin{align}
\label{eq:hj_decomposition_lemma_v2}
    \big( \nabla A(q_j) z \big) \cdot \nabla f(z) &= -\tfrac{1}{2}\big(\nabla \cdot A(q_j)\big) f(z) + \Theta_j(z) ,
\end{align}
with $-\tfrac{1}{2}\big(\nabla \cdot A(q_j)\big) f$ a multiple of $f$ and $\Theta_j \in \{f\}^\perp$.
\end{minipage}
\hfill
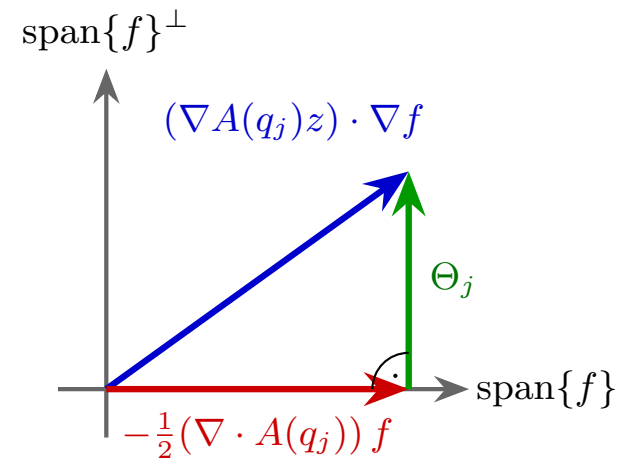
\begin{minipage}[c]{0.34\textwidth}
\centering
\begin{tikzpicture}[>=Stealth, scale=0.8, transform shape]
    \useasboundingbox (-2.2,-1.75) rectangle (2.7,2.0);
    \draw[thick, ->, black!60] (-2.0,-1.15) -- (1.4,-1.15)
        node[right, text=black, xshift=-2pt] {\small$\operatorname{span}\{f\}$};
    \draw[thick, ->, black!60] (-1.6,-1.55) -- (-1.6,1.5)
        node[above, text=black] {\small$\operatorname{span}\{f\}^{\perp}$};
    \draw[very thick, blue!80!black, ->] (-1.6,-1.15) -- (0.9,0.65);
    \node[blue!80!black, above] at (-0.05,0.78) {\small$(\nabla A(q_j)z)\cdot\nabla f$};
    \draw[very thick, red!80!black, ->] (-1.6,-1.15) -- (0.9,-1.15);
    \node[red!80!black, below] at (-0.35,-1.22) {\small$-\tfrac{1}{2}(\nabla\cdot A(q_j))\,f$};
    \draw[very thick, green!60!black, ->] (0.9,-1.15) -- (0.9,0.65);
    \node[green!45!black, right] at (0.95,-0.25) {\small$\Theta_j$};
    \draw (0.6,-1.15) arc (180:90:0.3);
    \fill (0.79,-1.04) circle (0.5pt);
\end{tikzpicture}

\vspace{0.5ex}
{\small Orthogonal projection onto $f$. TikZ code drafted with the assistance of AI.}
\end{minipage}
\vspace{1ex}

Substituting this splitting into the expanded magnetic difference cross-terms produces:
\begin{align}
\label{eq:cross_term_evaluation_v2}
   & f^x_j \Big( \nabla A(q_j) \frac{y_j - q_j}{\sqrt{\hbar}} \Big) \cdot (\nabla f)^y_j + \Big( \nabla A(q_j) \frac{x_j - q_j}{\sqrt{\hbar}} \Big) \cdot (\nabla f)^x_j \, f^y_j\nonumber \\
    &= - \big(\nabla \cdot A(q_j)\big) f^x_j f^y_j + \Theta_j\Big( \frac{x_j - q_j}{\sqrt{\hbar}} \Big) f^y_j \quad + f^x_j \Theta_j\Big( \frac{y_j - q_j}{\sqrt{\hbar}} \Big) .
\end{align}
Multiplied by $\sqrt{\hbar}$, its first term cancels the last term of \eqref{eq:I_div_j_raw_to_symmetrized}, which splits the remainder into the orthogonal and non-linear components:

\begin{align}
\label{eq:I_div_bifurcated_unified}
    \mathcal{I}_{\text{div}, j}^{(k)} &= \tfrac{c_{N,k}}{(2\pi)^{3k}} \hbar^{-3k/2 - 1/2} \int_{\mathbb{R}^{6k}} \dZk \int_{\mathbb{R}^{6k}} \dXkYk \, \Phi(Z_k) \prod_{m \neq j}^k f^x_m f^y_m
        e^{+}_\hbar \, \dmat{k}(X;Y) \nonumber \\
    &\quad \times \Bigg(\underbrace{ \cancel{- \sqrt{\hbar} \big(\nabla \cdot A(q_j)\big) f^x_j f^y_j} }_{\text{Parallel component}} + \underbrace{ \sqrt{\hbar} \Big[ \Theta_j\Big( \frac{x_j - q_j}{\sqrt{\hbar}} \Big) f^y_j + f^x_j \Theta_j\Big( \frac{y_j - q_j}{\sqrt{\hbar}} \Big) \Big] }_{\text{Orthogonal part } \mathcal{I}_{\text{Orth}, j}^{(k)}} \nonumber \\
    &\qquad \quad + \underbrace{ f^x_j \, \mathrm M_A(y_j, q_j)(y_j - q_j) \cdot (\nabla f)^y_j }_{\text{Non-linear part } (y_j)}
        + \underbrace{ \mathrm M_A(x_j, q_j)(x_j - q_j) \cdot (\nabla f)^x_j \, f^y_j }_{\text{Non-linear part } (x_j)} + \underbrace{ \cancel{\sqrt{\hbar} \big(\nabla \cdot A(q_j)\big) f^x_j f^y_j} }_{\text{Divergence cancellation}} \Bigg) \nonumber \\
    &= \mathcal{I}_{\text{Orth}, j}^{(k)} + \mathcal{I}_{\text{NL}, j}^{(k)} .
\end{align}
Isolating the parts yields the distinct functional definitions:
\begin{align}
\label{eq:I_Orth_defined}
    \mathcal{I}_{\text{Orth}, j}^{(k)}
    &\coloneqq \tfrac{c_{N,k}}{(2\pi)^{3k}} \hbar^{-3k/2 - 1/2} \int_{\mathbb{R}^{6k}} \dZk \int_{\mathbb{R}^{6k}} \dXkYk \, \Phi(Z_k)
    \prod_{m \neq j}^k f^x_m f^y_m
    \nonumber \\
    &\quad \times e^{+}_\hbar \, \sqrt{\hbar} \Big[ \Theta_j\Big( \frac{x_j - q_j}{\sqrt{\hbar}} \Big) f^y_j + f^x_j \Theta_j\Big( \frac{y_j - q_j}{\sqrt{\hbar}} \Big) \Big] \dmat{k}(X;Y) ,
\end{align}
and
\begin{align}
\label{eq:I_NL_defined}
    \mathcal{I}_{\text{NL}, j}^{(k)} &\coloneqq \tfrac{c_{N,k}}{(2\pi)^{3k}} \hbar^{-3k/2 - 1/2} \int_{\mathbb{R}^{6k}} \dZk \int_{\mathbb{R}^{6k}} \dXkYk \, \Phi(Z_k) \prod_{m \neq j}^k f^x_m f^y_m
        e^{+}_\hbar \nonumber \\
    &\quad \times \Bigg[ f^x_j \, \mathrm M_A(y_j, q_j)(y_j - q_j) \cdot (\nabla f)^y_j
        + \mathrm M_A(x_j, q_j)(x_j - q_j) \cdot (\nabla f)^x_j \, f^y_j \Bigg] \dmat{k}(X;Y) .
\end{align}

\paragraph{Orthogonal part $\mathcal{I}_{\text{Orth}, j}^{(k)}$ in \eqref{eq:I_Orth_defined}.}

First, substitute $q_j = x_j - \sqrt{\hbar}z$. Second, apply the Fundamental Theorem of Calculus at $q_j = x_j - \sqrt{\hbar}z$ (and $q_j = y_j - \sqrt{\hbar}z$) to shift the observable together with the residual $\Theta_j$, which depends on $q_j$ through $\nabla A(q_j)$ and $\nabla \cdot A(q_j)$ in \eqref{eq:orthprojection_refined}. This dependence is Lipschitz by Assumption \ref{assump:A1}, so the identity holds for a.e.\ $x_j$:
\begin{align}
\label{eq:I_Orth_observable_shift}
    \big[ \Phi(Z_k) \, e^{+}_\hbar \, \Theta_j(z) \big]_{q_j = x_j - \sqrt{\hbar}z} &= \big[ \Phi(Z_k) \, e^{+}_\hbar \, \Theta_j(z) \big]_{q_j = x_j}
        - \sqrt{\hbar} z \cdot \int_0^1 d\sigma \, \nabla_{q_j} \big[ \Phi(Z_k) \, e^{+}_\hbar \, \Theta_j(z) \big]_{q_j = x_j - \sigma\sqrt{\hbar} z} .
\end{align}

\vspace{1ex}
\noindent
\begin{minipage}[c]{0.62\textwidth}
Third, the orthogonality \eqref{eq:correct_orthogonality_v2} enters through the cross-correlation of the residual
with the displaced profile, for $u \in \mathbb{R}^3$ and fixed $q_j \in \mathbb{R}^3$,
\begin{align}
\label{eq:Theta_correlation_def}
    \widehat{\Theta}_j(u) &= \int_{\mathbb{R}^3} dz \, \Theta_j(z) \, f(z - u) ,
\end{align}
so that $\widehat{\Theta}_j(0) = \langle \Theta_j, f \rangle_{L^2(\mathbb{R}^3)} = 0$ for every $q_j \in \mathbb{R}^3$. Since $\Theta_j \in L^2(\mathbb{R}^3)$ and $\nabla f \in L^2(\mathbb{R}^3)$ by Assumption \ref{assump:A4},
the convolution $\widehat{\Theta}_j$ is continuously differentiable with bounded gradient
\begin{align}
\label{eq:Theta_correlation_gradient}
\nabla \widehat{\Theta}_j(u) &= - \int_{\mathbb{R}^3} dz \, \Theta_j(z) \, (\nabla f)(z - u).
\end{align}
\end{minipage}
\hfill
\begin{minipage}[c]{0.34\textwidth}
\centering
\begin{tikzpicture}[>=Stealth, scale=0.8, transform shape]
    \useasboundingbox (-2.5,-1.9) rectangle (2.5,2.5);
    \fill[orange!12] (0,0) -- (1.75,1.66) -- (1.75,-1.66) -- cycle;
    \fill[orange!12] (0,0) -- (-1.75,1.66) -- (-1.75,-1.66) -- cycle;
    \draw[orange!75!black, thin] (-1.75,1.66) -- (0,0) -- (1.75,1.66);
    \draw[orange!75!black, thin] (-1.75,-1.66) -- (0,0) -- (1.75,-1.66);
    \draw[thick, ->, black!60] (-2.15,0) -- (2.15,0) node[right, text=black] {\small$u$};
    \draw[thick, ->, black!60] (0,-1.85) -- (0,2.0);
    \draw[very thick, green!60!black] plot[smooth, tension=0.7]
        coordinates {(-2.12,0) (-1.92,-0.005) (-1.72,-0.03) (-1.5,-0.11) (-1.25,-0.30)
                     (-0.95,-0.42) (-0.6,-0.40) (-0.3,-0.25) (0,0)
                     (0.35,0.32) (0.7,0.60) (1.05,0.74) (1.4,0.66) (1.7,0.38)
                     (1.9,0.15) (2.03,0.04) (2.12,0.005)};
    \node[green!45!black] at (1.15,1.12) {\small$\widehat{\Theta}_j$};
    \fill[black] (0,0) circle (1.6pt);
    \node[orange!75!black] at (0.05,2.28)
        {\small$\pm\,\|\Theta_j\|_{L^2}\|\nabla f\|_{L^2}\,|u|$};
\end{tikzpicture}

\vspace{0.5ex}
{\small Zero of $\widehat{\Theta}_j$ at the origin. TikZ code drafted with the assistance of AI.}
\end{minipage}
\vspace{1ex}
Applying the Fundamental Theorem of Calculus to $\widehat{\Theta}_j$ on the segment joining $0$ to
$(x_j-y_j)/\sqrt{\hbar}$, and using $\widehat{\Theta}_j(0) = 0$, yields
\begin{align}
\label{eq:I_Orth_FTC_shift}
\widehat{\Theta}_j\left(\frac{x_j-y_j}{\sqrt{\hbar}}\right) &= \frac{x_j-y_j}{\sqrt{\hbar}} \cdot \int_0^1 d\tau \, \nabla \widehat{\Theta}_j\left(\tau \frac{x_j-y_j}{\sqrt{\hbar}}\right).
\end{align}

Applying these steps in sequence expands the remainder into four parts:
\begin{align}
\label{eq:I_Orth_expanded_shifts}
    \mathcal{I}_{\text{Orth}, j}^{(k)} &= \tfrac{c_{N,k}}{(2\pi)^{3k}} \hbar^{-3(k-1)/2} \int_{\mathbb{R}^{3(k-1)}} (dq)^{\otimes k \setminus \{j\}} \int_{\mathbb{R}^{6k}} \dXkYk \int_{\mathbb{R}^{3k}} (dv)^{\otimes k}
        \prod_{m \neq j}^k f^x_m f^y_m \, \dmat{k}(X;Y) \nonumber \\
    &\quad \times \Bigg[ \underbrace{ \big[ \Phi(Z_k) \, e^{+}_\hbar \big]_{q_j=x_j} \left( \frac{x_j-y_j}{\sqrt{\hbar}} \cdot \int_0^1 d\tau \, \nabla \widehat{\Theta}_j\left(\tau\frac{x_j-y_j}{\sqrt{\hbar}}\right)\Big|_{q_j=x_j} \right) }_{\text{Convolution Shift } x_j} \nonumber \\
    &\qquad + \underbrace{ \big[ \Phi(Z_k) \, e^{+}_\hbar \big]_{q_j=y_j} \left( \frac{y_j-x_j}{\sqrt{\hbar}} \cdot \int_0^1 d\tau \, \nabla \widehat{\Theta}_j\left(\tau\frac{y_j-x_j}{\sqrt{\hbar}}\right)\Big|_{q_j=y_j} \right) }_{\text{Convolution Shift } y_j} \nonumber \\
    &\qquad - \underbrace{ \sqrt{\hbar} \int_{\mathbb{R}^3} dz \, f\left(z - \frac{x_j-y_j}{\sqrt{\hbar}}\right) z \cdot \int_0^1 d\sigma \, \nabla_{q_j} \big[ \Phi(Z_k) \, e^{+}_\hbar \, \Theta_j(z) \big]_{q_j = x_j - \sigma\sqrt{\hbar} z} }_{\text{Observable Shift } x_j} \nonumber \\
    &\qquad - \underbrace{ \sqrt{\hbar} \int_{\mathbb{R}^3} dz \, f\left(z - \frac{y_j-x_j}{\sqrt{\hbar}}\right) z \cdot \int_0^1 d\sigma \, \nabla_{q_j} \big[ \Phi(Z_k) \, e^{+}_\hbar \, \Theta_j(z) \big]_{q_j = y_j - \sigma\sqrt{\hbar} z} }_{\text{Observable Shift } y_j} \Bigg] .
\end{align}

Apply phase derivative identities, valid for every value of $q_j$, to transfer spatial differences into velocity gradients:
\begin{align}
\label{eq:phase_derivative_identities}
        \frac{x_j - y_j}{\sqrt{\hbar}} e^{+}_\hbar &= -i\sqrt{\hbar} \nabla_{v_j} e^{+}_\hbar, \qquad
    \nabla_{q_j} e^{+}_\hbar = - \big( \nabla A(q_j) \big)^T \nabla_{v_j} e^{+}_\hbar .
\end{align}

\emph{Convolution shift.} Integration by parts gives:
\begin{align}
\label{eq:profile_shift_IBP_xj}
    &\int_{\mathbb{R}^{3k}} (dv)^{\otimes k} \Phi(Z_k|_{q_j=x_j}) \Big( -i\sqrt{\hbar} \nabla_{v_j} e^{+}_\hbar \big|_{q_j=x_j} \Big)
    \cdot \int_0^1 d\tau \, \nabla \widehat{\Theta}_j\Big(\tau \frac{x_j-y_j}{\sqrt{\hbar}}\Big)\Big|_{q_j=x_j}
    \nonumber \\
    &\quad = i \sqrt{\hbar} \int_{\mathbb{R}^{3k}} (dv)^{\otimes k} e^{+}_\hbar \big|_{q_j=x_j}
    \nabla_{v_j}\Phi(Z_k|_{q_j=x_j}) \cdot \int_0^1 d\tau \, \nabla \widehat{\Theta}_j\Big(\tau \frac{x_j-y_j}{\sqrt{\hbar}}\Big)\Big|_{q_j=x_j} .
\end{align}
The $y_j$-term follows from the first identity in \eqref{eq:phase_derivative_identities}, multiplied by $-1$, with $x_j$ replaced by $y_j$ in $\Phi$, $e^{+}_\hbar$ and $\widehat{\Theta}_j$, and yields the prefactor $-i\sqrt{\hbar}$.

\emph{Observable shift.} The second identity of \eqref{eq:phase_derivative_identities} and an integration by parts over $v_j$ in the middle term give:
\begin{align}
\label{eq:observable_shift_IBP_xj}
    & - \sqrt{\hbar} \int_0^1 d\sigma \int_{\mathbb{R}^{3k}} (dv)^{\otimes k} \Big[ z \cdot \nabla_{q_j} \big( \Phi(Z_k) \, e^{+}_\hbar \, \Theta_j(z) \big) \Big]_{q_j = x_j - \sigma\sqrt{\hbar} z} \nonumber \\
    &= - \sqrt{\hbar} \int_0^1 d\sigma \int_{\mathbb{R}^{3k}} (dv)^{\otimes k} \Bigg[ \big( z \cdot \nabla_{q_j}\Phi(Z_k) \big) e^{+}_\hbar \, \Theta_j(z)
        - \big( \Phi(Z_k) \nabla A(q_j) z \big) \cdot \nabla_{v_j} e^{+}_\hbar \, \Theta_j(z) \nonumber \\
    &\hphantom{{}= - \sqrt{\hbar} \int_0^1 d\sigma \int_{\mathbb{R}^{3k}} (dv)^{\otimes k} \Bigg[} + \Phi(Z_k) \, e^{+}_\hbar \, z \cdot \nabla_{q_j} \Theta_j(z) \Bigg]_{q_j = x_j - \sigma\sqrt{\hbar} z} \nonumber \\
    &= - \sqrt{\hbar} \int_0^1 d\sigma \int_{\mathbb{R}^{3k}} (dv)^{\otimes k} \Bigg[ e^{+}_\hbar \Bigg( z \cdot \big( \nabla_{q_j}\Phi(Z_k) + \nabla A(q_j)^T \nabla_{v_j}\Phi(Z_k) \big) \Theta_j(z) \nonumber \\
    &\hphantom{{}= - \sqrt{\hbar} \int_0^1 d\sigma \int_{\mathbb{R}^{3k}} (dv)^{\otimes k} \Bigg[ e^{+}_\hbar \Bigg(} + \Phi(Z_k) \, z \cdot \nabla_{q_j} \Theta_j(z) \Bigg) \Bigg]_{q_j = x_j - \sigma\sqrt{\hbar} z} ,
\end{align}
and the $y_j$-term follows in the same way at $q_j = y_j - \sigma\sqrt{\hbar} z$, with the same sign.

\emph{Sum of the shifts.} Collect the four shift contributions into the unified complex test function:
\begin{align}
\label{eq:tilde_Phi_Orth_def}
    \tilde{\Phi}_{\text{Orth}}(Z_k) &\coloneqq i \, e^{-i \frac{A(x_j) \cdot (x_j-y_j)}{\hbar}} \Big(\nabla_{v_j}\Phi(Z_k|_{q_j=x_j})\Big) \cdot \int_0^1 d\tau \, \nabla \widehat{\Theta}_j\left(\tau \frac{x_j-y_j}{\sqrt{\hbar}}\right)\Big|_{q_j=x_j} \nonumber \\
    &\quad - i \, e^{-i \frac{A(y_j) \cdot (x_j-y_j)}{\hbar}} \Big(\nabla_{v_j}\Phi(Z_k|_{q_j=y_j})\Big) \cdot \int_0^1 d\tau \, \nabla \widehat{\Theta}_j\left(\tau \frac{y_j-x_j}{\sqrt{\hbar}}\right)\Big|_{q_j=y_j} \nonumber \displaybreak[0] \\
    &\quad - \int_{\mathbb{R}^3} dz \, f\left(z - \frac{x_j-y_j}{\sqrt{\hbar}}\right) \int_0^1 d\sigma \, e^{-i \frac{A(x_j - \sigma\sqrt{\hbar}z) \cdot (x_j-y_j)}{\hbar}} \nonumber \\
    &\qquad \times z \cdot \Big[ \Big( \big( \nabla_{q_j} + \nabla A(q_j)^T \nabla_{v_j} \big) \Phi(Z_k) \Big) \Theta_j(z) + \Phi(Z_k) \, \nabla_{q_j} \Theta_j(z) \Big]_{q_j=x_j-\sigma\sqrt{\hbar}z} \nonumber \\
    &\quad - \int_{\mathbb{R}^3} dz \, f\left(z - \frac{y_j-x_j}{\sqrt{\hbar}}\right) \int_0^1 d\sigma \, e^{-i \frac{A(y_j - \sigma\sqrt{\hbar}z) \cdot (x_j-y_j)}{\hbar}} \nonumber \\
    &\qquad \times z \cdot \Big[ \Big( \big( \nabla_{q_j} + \nabla A(q_j)^T \nabla_{v_j} \big) \Phi(Z_k) \Big) \Theta_j(z) + \Phi(Z_k) \, \nabla_{q_j} \Theta_j(z) \Big]_{q_j=y_j-\sigma\sqrt{\hbar}z} .
\end{align}
Each contribution carries $e^{+}_\hbar$ with $q_j$ at $p \in \{x_j, \, y_j, \, x_j - \sigma\sqrt{\hbar}z, \, y_j -
\sigma\sqrt{\hbar}z\}$, which is the exponential factor of the corresponding term in \eqref{eq:tilde_Phi_Orth_def} times
the phase $e^{\frac{i}{\hbar} ( v_j \cdot (x_j - y_j) + \sum_{n \neq j} (v_n - A(q_n)) \cdot (x_n - y_n) )}$, free of
$q_j$.

The Cauchy--Schwarz and Young inequalities decouple the trace, and the map $X \leftrightarrow Y$ on the second half generates the symmetric phase addition:
\begin{align}
\label{eq:I_Orth_resolved}
    \big| \mathcal{I}_{\text{Orth}, j}^{(k)} \big| &\le \tfrac{c_{N,k}}{(2\pi)^{3k}} \hbar^{1/2} \, \hbar^{-3(k-1)/2} \int_{\mathbb{R}^{3(k-1)}} (dq)^{\otimes k \setminus \{j\}} \int_{\mathbb{R}^{6k}} \dXkYk
        \prod_{m \neq j}^k |f^x_m f^y_m| \, \big| \dmat{k}(X;Y) \big| \nonumber \\
    &\quad \times \Big| \int_{\mathbb{R}^{3k}} (dv)^{\otimes k} \tilde{\Phi}_{\text{Orth}}(Z_k) \, e^{\frac{i}{\hbar} ( v_j \cdot (x_j - y_j) + \sum_{n \neq j} (v_n - A(q_n)) \cdot (x_n - y_n) )} \Big| \nonumber \\
    &\le \tfrac{c_{N,k}}{(2\pi)^{3k}} \hbar^{1/2} \, \hbar^{-3(k-1)/2} \int_{\mathbb{R}^{3(k-1)}} (dq)^{\otimes k \setminus \{j\}} \int_{\mathbb{R}^{6k}} \dXkYk
        \| a_X \Psi_{N,t} \|^2 \prod_{m \neq j}^k |f^x_m|^2 \nonumber \\
    &\quad \times \frac{1}{2} \bigg( \Big| \int_{\mathbb{R}^{3k}} (dv)^{\otimes k} \tilde{\Phi}_{\text{Orth}}(Z_k) \, e^{\frac{i}{\hbar} ( v_j \cdot (x_j - y_j) + \sum_{n \neq j} (v_n - A(q_n)) \cdot (x_n - y_n) )} \Big| \nonumber \\
    &\quad + \Big| \int_{\mathbb{R}^{3k}} (dv)^{\otimes k} \overline{\tilde{\Phi}_{\text{Orth}}(Z_k)} \, e^{-\frac{i}{\hbar} ( v_j \cdot (x_j - y_j) + \sum_{n \neq j} (v_n - A(q_n)) \cdot (x_n - y_n) )} \Big| \bigg) ,
\end{align}
since $\Phi$, $f$, $\Theta_j$ and $A$ are real-valued and $\tilde{\Phi}_{\text{Orth}}$ with $x_j$ and $y_j$ interchanged
equals $\overline{\tilde{\Phi}_{\text{Orth}}}$. The moduli in \eqref{eq:I_Orth_resolved} do not depend on $q_j$. Since $\|f\|_{L^2(\mathbb{R}^3)} = 1$, multiplying by
$1 = \hbar^{-3/2}\displaystyle \int_{\mathbb{R}^3} dq_j \, |f^x_j|^2$ reinstates $q_j$ as an integration variable,
restoring the full tensor measure $(dq)^{\otimes k}$ and the prefactor $\hbar^{1/2-3k/2}$. Inserting
$e^{\mp \frac{i}{\hbar} A(q_j) \cdot (x_j-y_j)}$, of modulus $1$ and independent of $v$, in the two moduli turns their phases
into $e^{\pm}_\hbar$:
\begin{align}
\label{eq:I_Orth_master_form}
    \big| \mathcal{I}_{\text{Orth}, j}^{(k)} \big|
    &\le \tfrac{c_{N,k}}{(2\pi)^{3k}} \hbar^{\frac{1}{2} - \frac{3k}{2}} \int_{\mathbb{R}^{3k}} (dq)^{\otimes k}
    \int_{\mathbb{R}^{6k}} \dXkYk \, \| a_X \Psi_{N,t} \|^2 |f^x_j|^2 \prod_{m \neq j}^k |f^x_m|^2
    \nonumber \\
    &\quad \times \frac{1}{2} \bigg( \Big| \int_{\mathbb{R}^{3k}} (dv)^{\otimes k} \tilde{\Phi}_{\text{Orth}}(Z_k) \, e^{+}_\hbar \Big|
    + \Big| \int_{\mathbb{R}^{3k}} (dv)^{\otimes k} \overline{\tilde{\Phi}_{\text{Orth}}(Z_k)} \, e^{-}_\hbar \Big| \bigg) .
\end{align}
The test function $\tilde{\Phi}_{\text{Orth}}$ carries the difference coordinate
$(x_j-y_j)/\sqrt{\hbar}$, so the associated constant requires a supremum over $x_j$ and $y_j$.

\paragraph{Finiteness of the orthogonal phase-space bound.}
By \eqref{eq:tilde_Phi_Orth_def}, $\tilde{\Phi}_{\text{Orth}}$ is linear in $\Phi$ and vanishes for $|v_n| > R_\Phi$, with $R_\Phi$ of Proposition \ref{thm:master_remainder_bound}.
Proposition \ref{thm:master_remainder_bound} needs the $v$-derivatives of $\tilde{\Phi}_{\text{Orth}}$ bounded uniformly
in $(q, x_j, y_j)$, so that the constant in \eqref{eq:master_final_exponent_annihilation} is finite. The Cauchy--Schwarz
inequality bounds the $\tau$- and $z$-integrals in \eqref{eq:tilde_Phi_Orth_def} by constants independent of $q_j$ and
$(x_j - y_j)/\sqrt{\hbar}$. Bounding the $\sigma$-integral by the spatial supremum decouples \eqref{eq:tilde_Phi_Orth_def}. In \eqref{eq:tilde_Phi_Orth_def} the $v$-derivatives fall on $\Phi$, and for every $s \in \mathbb{N}_0$ we find
\begin{align}
\label{eq:tilde_Phi_Orth_bound}
    \left| \nabla_v^{s} \tilde{\Phi}_{\text{Orth}}(Z_k) \right| &\le \C \Big( \left| \nabla_{v_j} \nabla_v^{s} \Phi(Z_k|_{q_j=x_j}) \right| + \left| \nabla_{v_j} \nabla_v^{s} \Phi(Z_k|_{q_j=y_j}) \right| \nonumber \\
    &\quad + \sup_{p \in \mathbb{R}^3} \left| \Big( \nabla_{q_j} \nabla_v^{s} \Phi + \nabla A(p)^T \nabla_{v_j} \nabla_v^{s} \Phi \Big)\big(Z_k|_{q_j = p}\big) \right|
        + \sup_{p \in \mathbb{R}^3} \left| \nabla_v^{s} \Phi\big(Z_k|_{q_j = p}\big) \right| \Big) ,
\end{align}
where $\C$ depends on $\|\nabla A\|_{L^\infty(\mathbb{R}^3)}$, $\|\nabla^2 A\|_{L^\infty(\mathbb{R}^3)}$ and on $f$ through
Assumption \ref{assump:A4}, and the right side is bounded, since $\Phi \in C_c^\infty(\mathbb{R}^{6k})$.

\paragraph{Non-linear part $\mathcal{I}_{\text{NL}, j}^{(k)}$ in \eqref{eq:I_NL_defined}.}
\label{par:trace_decoupling_NL}

For the non-linear part we take absolute values and apply the operator bound
\eqref{eq:M_components_bound},
\begin{align}
\label{eq:M_components_bound_NL}
    \| \mathrm M_A(z,q_j) \|_{\text{op}} &\le \frac{1}{2} \|\nabla^2 A\|_{L^\infty(\mathbb{R}^3)} |z - q_j|, \quad \forall z \in \{x_j, y_j\}.
\end{align}
Each quadratic deviation then
contributes $|z-q_j|^2 = \hbar \big| (z-q_j)/\sqrt{\hbar} \big|^2$, which results in an additional $\hbar$, modifying the prefactor $\hbar^{-3k/2-1/2}$ of \eqref{eq:I_NL_defined} into $\hbar^{1/2-3k/2}$. The Cauchy--Schwarz and Young inequalities decouple the trace, and the map $X \leftrightarrow Y$ on the
second half condenses the two terms onto the $X$-coordinates and generates the symmetric phase addition $\frac{1}{2}\sum_\pm$:
\begin{align}
\label{eq:I_NL_resolved}
    \big| \mathcal{I}_{\text{NL}, j}^{(k)} \big| &\le \tfrac{c_{N,k}}{(2\pi)^{3k}} \frac{1}{2} \|\nabla^2 A\|_{L^\infty(\mathbb{R}^3)} \hbar^{1/2 - 3k/2} \int_{\mathbb{R}^{3k}} (dq)^{\otimes k} \int_{\mathbb{R}^{6k}} \dXkYk
        \| a_X \Psi_{N,t} \| \| a_Y \Psi_{N,t} \| \nonumber \\
    &\quad \times \Big[ \Big| \frac{y_j - q_j}{\sqrt{\hbar}} \Big|^2 |(\nabla f)^y_j| |f^x_j| + \Big| \frac{x_j - q_j}{\sqrt{\hbar}} \Big|^2 |(\nabla f)^x_j| |f^y_j| \Big]
        \prod_{m \neq j}^k |f^x_m| |f^y_m| \, \Big| \int_{\mathbb{R}^{3k}} (dv)^{\otimes k} \Phi(Z_k) \, e^{+}_\hbar \Big|_\infty \nonumber \\
    &\le \tfrac{c_{N,k}}{(2\pi)^{3k}} \frac{1}{2} \|\nabla^2 A\|_{L^\infty(\mathbb{R}^3)} \hbar^{1/2 - 3k/2} \int_{\mathbb{R}^{3k}} (dq)^{\otimes k} \int_{\mathbb{R}^{6k}} \dXkYk \nonumber \\
    &\quad \times \| a_X \Psi_{N,t} \|^2 \Big[ \Big| \frac{x_j - q_j}{\sqrt{\hbar}} \Big|^4 |(\nabla f)^x_j|^2 + |f^x_j|^2 \Big] \prod_{m \neq j}^k |f^x_m|^2
        \frac{1}{2} \sum_{\pm} \Big| \int_{\mathbb{R}^{3k}} (dv)^{\otimes k} \Phi(Z_k) \, e^{\pm}_\hbar \Big|_\infty .
\end{align}

\subsection{Bounding the External Remainder Integral \texorpdfstring{$\mathcal{I}_{\text{ext}, j}^{(k)}$}{I-ext-j-k}}
\label{subsubsec:bound_I_ext_j_k}

The remainder $\mathcal{I}_{\text{ext}, j}^{(k)}$ \eqref{eq:unified_strong_weak_remainders} is defined with $\mathcal{R}_{\text{ext}, j}^{(k)}$ \eqref{eq:rkinextdef}. Resolving $\operatorname{Im}(\cdot)$ and swapping $x_n \leftrightarrow y_n$ on the conjugate cancels $\Vext(q_j)$, yielding the symmetrized remainder with $\Vext(x_j) - \Vext(y_j) - \nabla \Vext(q_j) \cdot (x_j-y_j)$. The Fundamental Theorem of Calculus eliminates the linear real-space component, the identity $(x_j - y_j) e^{+}_\hbar = -i\hbar \nabla_{v_j} e^{+}_\hbar$ turns $x_j - y_j$ into a velocity gradient, and the Mean Value Theorem bounds the gradient deviation of $\Vext$.
Velocity integration by parts transfers $-i\hbar \nabla_{v_j}$ onto $\Phi$. The Cauchy--Schwarz and Young inequalities decouple the trace, with amplitudes paired against $\mathbf{1}$, and the map $X \leftrightarrow Y$ on the second half converts $\left|\frac{y_j-q_j}{\sqrt{\hbar}}\right|$ into $\left|\frac{x_j-q_j}{\sqrt{\hbar}}\right|^2$ and generates the symmetric phase addition $\frac{1}{2}\sum_\pm$:
\begin{align}
\label{eq:I_ext_abs_bound}
    \big| \mathcal{I}_{\text{ext}, j}^{(k)} \big| &= \tfrac{c_{N,k}}{(2\pi)^{3k}} \frac{1}{\hbar} \hbar^{-\frac{3k}{2}} \Bigg| \int_{\mathbb{R}^{3k}} (dq)^{\otimes k} \int_{\mathbb{R}^{3k}} (dv)^{\otimes k} \int_{\mathbb{R}^{6k}} \dXkYk \, \Phi(Z_k) \nonumber \\
    &\quad \times \prod_{m = 1}^k f^x_m f^y_m \, e^{+}_\hbar \Big[ \Vext(x_j) - \Vext(y_j) - \nabla \Vext(q_j) \cdot (x_j - y_j) \Big] \dmat{k}(X;Y) \Bigg| \nonumber \\
    &= \tfrac{c_{N,k}}{(2\pi)^{3k}} \hbar^{-\frac{3k}{2}} \Bigg| \int_{\mathbb{R}^{3k}} (dq)^{\otimes k} \int_{\mathbb{R}^{3k}} (dv)^{\otimes k} \int_{\mathbb{R}^{6k}} \dXkYk
        \prod_{m = 1}^k f^x_m f^y_m \, e^{+}_\hbar \, \dmat{k}(X;Y) \nonumber \\
    &\quad \times \nabla_{v_j} \Phi(Z_k) \cdot \int_0^1 ds \, \Big[ \nabla \Vext\big(y_j + s(x_j - y_j)\big) - \nabla \Vext(q_j) \Big] \Bigg| \nonumber \\
    &\le \tfrac{c_{N,k}}{(2\pi)^{3k}} \frac{1}{2} \|\nabla^2 \Vext\|_{L^\infty(\mathbb{R}^3)} \hbar^{\frac{1}{2} - \frac{3k}{2}} \int_{\mathbb{R}^{3k}} (dq)^{\otimes k} \int_{\mathbb{R}^{6k}} \dXkYk \nonumber \\
    &\quad \times \| a_X \Psi_{N,t} \| \| a_Y \Psi_{N,t} \| \Big( \Big| \frac{x_j - q_j}{\sqrt{\hbar}} \Big| + \Big| \frac{y_j - q_j}{\sqrt{\hbar}} \Big| \Big) \prod_{m = 1}^k |f^x_m| |f^y_m|
        \Big| \int_{\mathbb{R}^{3k}} (dv)^{\otimes k} \nabla_{v_j} \Phi(Z_k) \, e^{+}_\hbar \Big|_\infty \nonumber \\
    &\le \tfrac{c_{N,k}}{(2\pi)^{3k}} \frac{1}{2} \|\nabla^2 \Vext\|_{L^\infty(\mathbb{R}^3)} \hbar^{\frac{1}{2} - \frac{3k}{2}} \int_{\mathbb{R}^{3k}} (dq)^{\otimes k} \int_{\mathbb{R}^{6k}} \dXkYk
        \| a_X \Psi_{N,t} \|^2 \Big( 1 + \Big| \frac{x_j - q_j}{\sqrt{\hbar}} \Big|^2 \Big) \prod_{m = 1}^k |f^x_m|^2 \nonumber \\
    &\quad \times \frac{1}{2} \sum_{\pm} \Big| \int_{\mathbb{R}^{3k}} (dv)^{\otimes k} \nabla_{v_j} \Phi(Z_k) \, e^{\pm}_\hbar \Big|_\infty .
\end{align}

\subsection{Bounding the Mean-Field Remainder Integral \texorpdfstring{$\mathcal{I}_{\text{mf}, j}^{(k)}$}{I-mf-j-k}}
\label{subsubsec:bound_I_mf_j_k}

The mean-field remainder $\mathcal{I}_{\text{mf}, j}^{(k)}$ \eqref{eq:unified_strong_weak_remainders} pairs $\mathcal{R}_{\text{mf}, j}^{(k)}$ \eqref{eq:R_mf_kin_def_prop} against the velocity gradient of the test function. Integrating over $v_{k+1}$ yields the Dirac distribution, which sets $x_{k+1} = y_{k+1}$. The Mean Value Theorem bounds the force deviation.
Resolving $\operatorname{Re}(\cdot)$ and swapping $x_n \leftrightarrow y_n$ symmetrizes the interaction deviation. The Cauchy--Schwarz and Young inequalities decouple the $(k+1)$-th marginal trace, and the map $X \leftrightarrow Y$ on the second half generates the symmetric phase addition $\frac{1}{2}\sum_\pm$:

\begin{align}
\label{eq:I_mf_abs_bound}
    & \big| \mathcal{I}_{\text{mf}, j}^{(k)} \big| = \Big| \int_{\mathbb{R}^{6k}} \dZk \, \nabla_{v_j} \Phi(Z_k) \cdot \mathcal{R}_{\text{mf}, j}^{(k)}(Z_k) \Big| \nonumber \\
    &= \tfrac{c_{N,k}}{(2\pi)^{3(k+1)}} \hbar^{-3(k+1)/2} \Bigg| \int_{\mathbb{R}^{6k}} \dZk \int_{\mathbb{R}^6} dq_{k+1} \, dv_{k+1} \int_{\mathbb{R}^{6(k+1)}} \dXkYkplus \nonumber \\
    &\quad \times \int_0^1 ds \, \nabla_{v_j} \Phi(Z_k) \cdot \operatorname{Re} \Bigg\{ \prod_{m = 1}^{k+1} f^x_m f^y_m \, e^{+}_\hbar \, e^{\frac{i}{\hbar} (v_{k+1} - A(q_{k+1})) \cdot (x_{k+1} - y_{k+1})} \nonumber \\
    &\quad \times \Big[ \nabla W\big(s y_j + (1-s) x_j - x_{k+1}\big) - \nabla W(q_j - q_{k+1}) \Big]
        \Big\langle \Psi_{N,t}, a_X^\dagger a_{x_{k+1}}^\dagger a_{y_{k+1}} a_Y \Psi_{N,t} \Big\rangle \Bigg\} \Bigg| \nonumber \\
    &= \tfrac{c_{N,k}}{(2\pi)^{3k}} \hbar^{3/2 - 3k/2} \Bigg| \int_{\mathbb{R}^{6k}} \dZk \int_{\mathbb{R}^3} dq_{k+1} \int_{\mathbb{R}^{3(k+1)}} (dx)^{\otimes (k+1)} \int_{\mathbb{R}^{3k}} (dy)^{\otimes k}
        \int_0^1 ds \prod_{m = 1}^k f^x_m f^y_m \, (f^x_{k+1})^2 \, e^{+}_\hbar \nonumber \\
    &\quad \times \Big\langle \Psi_{N,t}, a_X^\dagger a_{x_{k+1}}^\dagger a_{x_{k+1}} a_Y \Psi_{N,t} \Big\rangle
        \nabla_{v_j} \Phi(Z_k) \cdot \Big[ \nabla W\big(s y_j + (1-s) x_j - x_{k+1}\big) - \nabla W(q_j - q_{k+1}) \Big] \Bigg| \nonumber \\
    &\le \tfrac{c_{N,k}}{(2\pi)^{3k}} \frac{1}{2} \|\nabla^2 W\|_{L^\infty(\mathbb{R}^3)} \hbar^{2 - 3k/2} \int_{\mathbb{R}^{3k}} (dq)^{\otimes k} \int_{\mathbb{R}^3} dq_{k+1} \nonumber \\
    &\quad \times \int_{\mathbb{R}^{3(k+1)}} (dx)^{\otimes (k+1)} \int_{\mathbb{R}^{3k}} (dy)^{\otimes k} \big\| a_{x_{k+1}} a_X \Psi_{N,t} \big\| \big\| a_{x_{k+1}} a_Y \Psi_{N,t} \big\| \nonumber \\
    &\quad \times \Big( \Big| \frac{x_j - q_j}{\sqrt{\hbar}} \Big| + \Big| \frac{y_j - q_j}{\sqrt{\hbar}} \Big| + 2 \Big| \frac{x_{k+1} - q_{k+1}}{\sqrt{\hbar}} \Big| \Big) (f^x_{k+1})^2 \prod_{m = 1}^k |f^x_m| |f^y_m|
        \Big| \int_{\mathbb{R}^{3k}} (dv)^{\otimes k} \nabla_{v_j} \Phi(Z_k) \, e^{+}_\hbar \Big|_\infty \nonumber \\
    &\le \tfrac{c_{N,k}}{(2\pi)^{3k}} \frac{1}{2} \|\nabla^2 W\|_{L^\infty(\mathbb{R}^3)} \hbar^{2 - 3k/2} \int_{\mathbb{R}^{3k}} (dq)^{\otimes k} \int_{\mathbb{R}^3} dq_{k+1}
        \int_{\mathbb{R}^{3(k+1)}} (dx)^{\otimes (k+1)} \int_{\mathbb{R}^{3k}} (dy)^{\otimes k} \big\| a_{x_{k+1}} a_X \Psi_{N,t} \big\|^2 \nonumber \\
    &\quad \times \Big( 1 + \Big| \frac{x_j - q_j}{\sqrt{\hbar}} \Big|^2 \Big) \Big( 1 + \Big| \frac{x_{k+1} - q_{k+1}}{\sqrt{\hbar}} \Big|^2 \Big) (f^x_{k+1})^2 \prod_{m = 1}^k (f^x_m)^2
        \frac{1}{2} \sum_{\pm} \Big| \int_{\mathbb{R}^{3k}} (dv)^{\otimes k} \nabla_{v_j} \Phi(Z_k) \, e^{\pm}_\hbar \Big|_\infty .
\end{align}

The $(k+1)$-th marginal components are computed in turn:
\begin{align}
\label{eq:I_mf_background_resolve}
    \int_{\mathbb{R}^3} dq_{k+1} \Big( 1 + \Big| \frac{x_{k+1} - q_{k+1}}{\sqrt{\hbar}} \Big|^2 \Big) (f^x_{k+1})^2
        &= \hbar^{3/2} \Big( 1 + \big\| |\cdot| f \big\|_{L^2(\mathbb{R}^3)}^2 \Big), \nonumber
     \\
    \int_{\mathbb{R}^3} dx_{k+1} \big\| a_{x_{k+1}} a_X \Psi_{N,t} \big\|^2
    = \big\langle a_X &\Psi_{N,t}, \mathcal{N} \, a_X \Psi_{N,t} \big\rangle \le \hbar^{-3} \big\| a_X \Psi_{N,t} \big\|^2.
\end{align}
Substituting \eqref{eq:I_mf_background_resolve} into \eqref{eq:I_mf_abs_bound} yields the bound:
\begin{align}
\label{eq:I_mf_final_master_bound}
    \big| \mathcal{I}_{\text{mf}, j}^{(k)} \big| &\le \tfrac{c_{N,k}}{(2\pi)^{3k}} \frac{1}{2} \|\nabla^2 W\|_{L^\infty(\mathbb{R}^3)} \Big( 1 + \big\| |\cdot| f \big\|_{L^2(\mathbb{R}^3)}^2 \Big) \hbar^{1/2 - 3k/2}
        \int_{\mathbb{R}^{3k}} (dq)^{\otimes k} \int_{\mathbb{R}^{6k}} \dXkYk \, \| a_X \Psi_{N,t} \|^2 \nonumber \\
    &\quad \times \Big( 1 + \Big| \frac{x_j - q_j}{\sqrt{\hbar}} \Big|^2 \Big) (f^x_j)^2 \prod_{m \neq j}^k (f^x_m)^2
        \frac{1}{2} \sum_{\pm} \Big| \int_{\mathbb{R}^{3k}} (dv)^{\otimes k} \nabla_{v_j} \Phi(Z_k) \, e^{\pm}_\hbar \Big|_\infty .
\end{align}

\subsection{Bounding the Internal Correlation Remainder Integral \texorpdfstring{$\mathcal{I}_{\text{cor}, i, j}^{(k)}$}{I-cor-i-j-k}}
\label{subsubsec:bound_I_cor_j_k}

The remainder $\mathcal{I}_{\text{cor}, i, j}^{(k)}$ \eqref{eq:unified_strong_weak_remainders} is defined with $\mathcal{R}_{\text{cor}, i, j}^{(k)}$ \eqref{eq:R_intrn_kin_def_prop}. The Fundamental Theorem of Calculus for $W$ and the identity $\big[ (x_i - y_i) - (x_j - y_j) \big] e^{+}_\hbar = -i\hbar \, (\nabla_{v_i} - \nabla_{v_j}) \, e^{+}_\hbar$ prepare the bound.
$\operatorname{Im}(\cdot)$ vanishes in $|\cdot|$. The Cauchy--Schwarz and Young inequalities decouple the trace, and the map $X \leftrightarrow Y$ on the second half generates the symmetric phase addition $\frac{1}{2}\sum_\pm$:

\begin{align}
\label{eq:I_cor_abs_bound}
    \big| \mathcal{I}_{\text{cor}, i, j}^{(k)} \big| &= \tfrac{c_{N,k}}{(2\pi)^{3k}} \Bigg| \int_{\mathbb{R}^{6k}} \dZk \int_{\mathbb{R}^{6k}} \dXkYk \, \Phi(Z_k) \prod_{m = 1}^k \mathfrak{f}_m(x_m) \overline{\mathfrak{f}_m(y_m)} \nonumber \\
    &\quad \times \Big[ W(y_i - y_j) - W(x_i - x_j) \Big] \dmat{k}(X;Y) \Bigg| \nonumber \\
    &= \tfrac{c_{N,k}}{(2\pi)^{3k}} \hbar^{1 - 3k/2} \Bigg| \int_{\mathbb{R}^{3k}} (dq)^{\otimes k} \int_{\mathbb{R}^{3k}} (dv)^{\otimes k} \int_{\mathbb{R}^{6k}} \dXkYk
        \prod_{m = 1}^k f^x_m f^y_m \, e^{+}_\hbar \, \dmat{k}(X;Y) \nonumber \\
    &\quad \times \Big( \int_0^1 ds \, \nabla W\big(s(y_i - y_j) + (1-s)(x_i - x_j)\big) \Big) \cdot \Big[ (\nabla_{v_i} - \nabla_{v_j}) \Phi(Z_k) \Big] \Bigg| \nonumber \\
    &\le \tfrac{c_{N,k}}{(2\pi)^{3k}} \hbar^{1 - 3k/2} \|\nabla W\|_{L^\infty(\mathbb{R}^3)} \int_{\mathbb{R}^{3k}} (dq)^{\otimes k} \int_{\mathbb{R}^{6k}} \dXkYk \nonumber \\
    &\quad \times \| a_X \Psi_{N,t} \| \| a_Y \Psi_{N,t} \| \prod_{m = 1}^k |f^x_m| |f^y_m| \, \Big| \int_{\mathbb{R}^{3k}} (dv)^{\otimes k} (\nabla_{v_i} - \nabla_{v_j}) \Phi(Z_k) \, e^{+}_\hbar \Big|_\infty \nonumber \\
    &\le \tfrac{c_{N,k}}{(2\pi)^{3k}} \hbar^{1 - 3k/2} \|\nabla W\|_{L^\infty(\mathbb{R}^3)} \int_{\mathbb{R}^{3k}} (dq)^{\otimes k} \int_{\mathbb{R}^{6k}} \dXkYk \nonumber \\
    &\quad \times \| a_X \Psi_{N,t} \|^2 \prod_{m = 1}^k (f^x_m)^2 \, \frac{1}{2} \sum_{\pm} \Big| \int_{\mathbb{R}^{3k}} (dv)^{\otimes k} (\nabla_{v_i} - \nabla_{v_j}) \Phi(Z_k) \, e^{\pm}_\hbar \Big|_\infty .
\end{align}

\appendix

\section{Slater Determinants as Initial Data}
\label{app:slater_determinants}

This appendix proves the claims of Remark \ref{rem:initial_data_ground_states}: the Slater determinant of a trapped
ground state satisfies the hypotheses of Theorems \ref{thm:magnetic_vlasov_limit} and
\ref{thm:w2_energy_criterion_noV}. 

\paragraph{Preparation of Slater determinant in a trap.} Let $\Vext_0$ satisfy Assumption \ref{assump:A2} and
$\Vext_0(x) \ge c_0 |x|^2 - C_0$, with $c_0 > 0$ and $C_0 \ge 0$. Let $\varphi_1, \dots, \varphi_N$ be orthonormal eigenfunctions of
$- \tfrac12 \hbar^2 \Delta + \Vext_0$ to its $N$ lowest eigenvalues $\lambda_1 \le \dots \le \lambda_N$,
\begin{align}
\label{eq:ground_state_eigen}
    \big( - \tfrac12 \hbar^2 \Delta + \Vext_0 \big) \varphi_n &= \lambda_n \varphi_n ,
    \qquad
    \langle \varphi_n, \varphi_{n'} \rangle = \delta_{n n'} ,
    \qquad
    n, n' \le N .
\end{align}
The initial state $\Psi_{N,0}$ is their Slater determinant, a ground state of $H_N^{(0)}$,
\begin{align}
\label{eq:ground_state_slater}
    \Psi_{N,0}(x_1, \dots, x_N) \coloneqq (N!)^{-1/2} \det \big( \varphi_j(x_i) \big)_{i,j \le N} ,
    \\
    \label{eq:ground_state_H0}
    H_N^{(0)} \coloneqq \sum_{j=1}^N \Big( - \tfrac12 \hbar^2 \Delta_{x_j} + \Vext_0(x_j) \Big) ,
\qquad
    \langle \Psi_{N,0}, H_N^{(0)} \Psi_{N,0} \rangle &= \sum_{n=1}^N \lambda_n .
\end{align}
The magnetic Husimi measures of $\Psi_{N,0}$ are
\begin{align}
\label{eq:ground_state_wick}
    \mhus^{(k)}_{N,0}(Z_k)
    &= \tfrac{c_{N,k}}{(2\pi)^{3k}} \det \big( \langle \Psi_{N,0}, a^\dagger(\mathfrak{f}_j) \, a(\mathfrak{f}_i) \Psi_{N,0} \rangle \big)_{i,j \le k} .
\end{align}

For $N\rightarrow \infty$, we will prove that the limit of $\mhus^{(1)}_{N,0}$ can be identified by the following probability density with compact support.

\paragraph{Identification of $\mathbf{f}_0$.} By $\Vext_0(q) \ge c_0 |q|^2 - C_0$, the sets
$\{ |v - A(q)|^2 + 2\Vext_0(q) \le s \}$ are compact and their boundaries are null sets, so there is $s \in \mathbb{R}$
with
\begin{align}
\label{eq:ground_state_f0}
    \mathbf{f}_0(q,v) &\coloneqq \frac{1}{(2\pi)^3} \, \One\Big( |v - A(q)|^2 + 2\Vext_0(q) \le s \Big) ,
    \quad
    \iint_{\mathbb{R}^6} dq \, dv \, \mathbf{f}_0(q,v) = 1 .
\end{align}
By the bathtub principle \cite[Thm.~1.14]{liebloss2001analysis}, $\mathbf{f}_0$ is the minimizer of
$$\iint_{\mathbb{R}^6} dq \, dv \, \big( \tfrac12 |v - A(q)|^2 + \Vext_0(q) \big) \mathbf{f}$$
under restrictions
$0 \le \mathbf{f} \le (2\pi)^{-3}$ and $\iint_{\mathbb{R}^6} dq \, dv \, \mathbf{f} = 1$, unique up to null sets.

\vskip4mm
The following two propositions show separately that all the conditions for initial data in the two main results are satisfied by the above constructed Slater determinant.

\begin{proposition}[Hypotheses of Theorem \ref{thm:magnetic_vlasov_limit} for Slater Determinants]
\label{prop:slater_hypotheses}
Let Assumptions \ref{assump:A1}--\ref{assump:A4} hold and let $\Psi_{N,0}$ be the Slater determinant of
\eqref{eq:ground_state_slater}. Then:
\begin{enumerate}
\renewcommand{\theenumi}{(\roman{enumi})}
\renewcommand{\labelenumi}{\theenumi}
\item \label{item:slater_support}
Assumption \ref{assump:I0} holds for $\mathbf{f}_0$ of \eqref{eq:ground_state_f0}.
\item \label{item:slater_I1}
Assumption \ref{assump:I1} holds.
\item \label{item:slater_I2}
Assumption \ref{assump:I2} holds with the initial data $\mathbf{f}_0$ of \eqref{eq:ground_state_f0}, with the convergence
in the strong $L^1$ sense:
\begin{align}
\label{eq:ground_state_I2}
    \lim_{N \to \infty} \big\| \mhus^{(k)}_{N,0} - \mathbf{f}_0^{\otimes k} \big\|_{L^1(\mathbb{R}^{6k})} &= 0 ,
    \qquad k \ge 1 .
\end{align}
\end{enumerate}
\end{proposition}

\begin{proof}[Proof of Proposition \ref{prop:slater_hypotheses}]

\emph{Assumption \ref{assump:I0}.} By construction \eqref{eq:ground_state_f0}, $\mathbf{f}_0$ satisfies Assumption
\ref{assump:I0}.

\emph{Assumption \ref{assump:I1}.} For $\mathfrak{f}_1$ of \eqref{eq:coherent_shifted},
a direct computation gives
\begin{align}
\label{eq:ground_state_trial_energy}
    \big\langle \mathfrak{f}_1, \big( - \tfrac12 \hbar^2 \Delta + \Vext_0 \big) \mathfrak{f}_1 \big\rangle &= \tfrac12 |v - A(q)|^2 + \tfrac{\hbar}{2} \|\nabla f\|_{L^2(\mathbb{R}^3)}^2 + \int_{\mathbb{R}^3} du \, f^2(u) \, \Vext_0\big( q + \sqrt{\hbar} u \big) \nonumber \\
    &\le \tfrac12 |v - A(q)|^2 + \Vext_0(q) + \C \sqrt{\hbar} \big( 1 + |\nabla \Vext_0(q)| \big) .
\end{align}
We apply the variational principle for the sum of the $N$ lowest eigenvalues of
$- \tfrac12 \hbar^2 \Delta + \Vext_0$ to the operator
$$
\gamma \coloneqq N \iint_{\mathbb{R}^6} dq \, dv \, \mathbf{f}_0(q,v) \, |\mathfrak{f}_1\rangle\langle\mathfrak{f}_1|.
$$ 
With
\eqref{eq:ground_state_H0}, \eqref{eq:ground_state_trial_energy} and $\nabla \Vext_0$ bounded on the support of
$\mathbf{f}_0$, this gives
\begin{align}
\label{eq:ground_state_minmax}
    \tfrac{1}{N} \big\langle \Psi_{N,0}, H_N^{(0)} \Psi_{N,0} \big\rangle
    &\le \tfrac{1}{N} \operatorname{Tr} \Big( \big( - \tfrac12 \hbar^2 \Delta + \Vext_0 \big) \gamma \Big)
    \nonumber \\
    &= \iint_{\mathbb{R}^6} dq \, dv \, \mathbf{f}_0(q,v) \,
    \big\langle \mathfrak{f}_1, \big( - \tfrac12 \hbar^2 \Delta + \Vext_0 \big) \mathfrak{f}_1 \big\rangle
    \nonumber \\
    &\le \iint_{\mathbb{R}^6} dq \, dv \, \Big( \tfrac12 |v - A(q)|^2 + \Vext_0(q) \Big) \mathbf{f}_0(q,v) + \C \sqrt{\hbar} .
\end{align}
Both bounds of \eqref{eq:assumption_initial_data} reduce to the kinetic energy and the second spatial moment of
$\Psi_{N,0}$. We control both by $\Vext_0(x) \ge c_0 |x|^2 - C_0$ and \eqref{eq:ground_state_minmax}:
\begin{align}
\label{eq:ground_state_apriori}
    &\tfrac{1}{N} \Big\langle \Psi_{N,0}, \sum_{j=1}^N \big( - \tfrac12 \hbar^2 \Delta_{x_j} \big) \Psi_{N,0} \Big\rangle
    + \tfrac{c_0}{N} \int_{\mathbb{R}^3} dx \, |x|^2 \, \langle \Psi_{N,0}, a_x^\dagger a_x \Psi_{N,0} \rangle
    \nonumber \\
    &\le \tfrac{1}{N} \Big\langle \Psi_{N,0}, \sum_{j=1}^N \big( - \tfrac12 \hbar^2 \Delta_{x_j} \big) \Psi_{N,0} \Big\rangle
    + \tfrac{1}{N} \int_{\mathbb{R}^3} dx \, \Vext_0(x) \, \langle \Psi_{N,0}, a_x^\dagger a_x \Psi_{N,0} \rangle + C_0
    \nonumber \\
    &= \tfrac{1}{N} \big\langle \Psi_{N,0}, H_N^{(0)} \Psi_{N,0} \big\rangle + C_0
\le \C .
\end{align}
For the first bound of \eqref{eq:assumption_initial_data}, \eqref{eq:ground_state_apriori},
$|A|^2 + |\Vext| \le \C (1 + |x|^2)$ and the operator inequality
$\tfrac12 (-i\hbar \nabla + A)^2 \le (-i\hbar \nabla)^2 + |A|^2$ yield
\begin{align}
\label{eq:ground_state_energy_bound}
    &\tfrac{1}{N} \big\langle \Psi_{N,0}, \mathcal{H} \Psi_{N,0} \big\rangle \nonumber \\ 
    &\le \tfrac{1}{N} \Big\langle \Psi_{N,0}, \sum_{j=1}^N \big( - \hbar^2 \Delta_{x_j} \big) \Psi_{N,0} \Big\rangle + \tfrac{1}{N} \int_{\mathbb{R}^3} dx \, \big( |A(x)|^2 + |\Vext(x)| \big) \langle \Psi_{N,0}, a_x^\dagger a_x \Psi_{N,0} \rangle + \tfrac12 \|W\|_{L^\infty(\mathbb{R}^3)} \le \C .
\end{align}
For the second bound, by \eqref{eq:husimi_spatial_density} and \eqref{eq:ground_state_apriori},
\begin{align}
\label{eq:ground_state_moment_bound}
    \iint_{\mathbb{R}^6} dq \, dv \, |q|^2 \, \mhus^{(1)}_{N,0}(q,v) &= \hbar^3 \int_{\mathbb{R}^3} dx \, \langle \Psi_{N,0}, a_x^\dagger a_x \Psi_{N,0} \rangle \int_{\mathbb{R}^3} du \, f^2(u) \, |x - \sqrt{\hbar} u|^2 \nonumber \\
    &\le \tfrac{2}{N} \int_{\mathbb{R}^3} dx \, |x|^2 \, \langle \Psi_{N,0}, a_x^\dagger a_x \Psi_{N,0} \rangle + 2 \hbar \, \big\| |\cdot| f \big\|_{L^2(\mathbb{R}^3)}^2 \le \C .
\end{align}
Hence Assumption \ref{assump:I1} holds.

\emph{Assumption \ref{assump:I2} in the sector $k = 1$.} Computing as in \eqref{eq:proof_kin_husimi_derivation1} at
$t = 0$, with $|v - A(q)|^2$ in place of $v^2$ and $-i\hbar \nabla_x$ in place of $-i\hbar \nabla_x + A(q)$, we obtain
\begin{align}
\label{eq:ground_state_kinetic_identity}
    &\tfrac12 \iint_{\mathbb{R}^6} dq \, dv \, |v - A(q)|^2 \, \mhus^{(1)}_{N,0}(q,v)\nonumber \\
    &= \tfrac{\hbar^{3/2}}{2} \int_{\mathbb{R}^3} dx \int_{\mathbb{R}^3} dq \, \Big( f^2\big( \tfrac{x-q}{\sqrt{\hbar}} \big) \, \big\| \hbar \nabla_x a_x \Psi_{N,0} \big\|^2 + \hbar \, \Big| (\nabla f)\big( \tfrac{x-q}{\sqrt{\hbar}} \big) \Big|^2 \| a_x \Psi_{N,0} \|^2 \Big) \nonumber \\
    &= \tfrac{1}{N} \Big\langle \Psi_{N,0}, \sum_{j=1}^N \big( - \tfrac12 \hbar^2 \Delta_{x_j} \big) \Psi_{N,0} \Big\rangle + \tfrac{\hbar}{2} \|\nabla f\|_{L^2(\mathbb{R}^3)}^2 .
\end{align}
Next, by Assumption \ref{assump:I1}, \eqref{eq:energy_dec_external_bound_noV} at $t = 0$, with $\Vext_0$ in place of
$\Vext$, gives
\begin{align}
\label{eq:ground_state_trap_identity}
    \Big| \iint_{\mathbb{R}^6} dq \, dv \, \Vext_0(q) \, \mhus^{(1)}_{N,0}(q,v)
    - \tfrac{1}{N} \int_{\mathbb{R}^3} dx \, \Vext_0(x) \, \langle \Psi_{N,0}, a_x^\dagger a_x \Psi_{N,0} \rangle \Big|
    &\le \C \sqrt{\hbar} .
\end{align}
Adding and subtracting $\tfrac{1}{N} \int_{\mathbb{R}^3} dx \, \Vext_0(x) \langle \Psi_{N,0}, a_x^\dagger a_x \Psi_{N,0} \rangle$, \eqref{eq:ground_state_kinetic_identity} and \eqref{eq:ground_state_trap_identity} give
\begin{align}
\label{eq:ground_state_husimi_energy}
    &\Bigg| \iint_{\mathbb{R}^6} dq \, dv \, \Big( \tfrac12 |v - A(q)|^2 + \Vext_0(q) \Big) \mhus^{(1)}_{N,0}(q,v)
   - \tfrac{1}{N} \big\langle \Psi_{N,0}, H_N^{(0)} \Psi_{N,0} \big\rangle \Bigg|
    \le \C \sqrt{\hbar} .
\end{align}
With $0 \le \mhus^{(1)}_{N,0} \le \tfrac{1}{(2\pi)^3}$ and
$\iint_{\mathbb{R}^6} dq \, dv \, ( \mhus^{(1)}_{N,0} - \mathbf{f}_0 ) = 0$ by \eqref{eq:normalized_husimi_measure},
both factors $\tfrac12 |v - A(q)|^2 + \Vext_0(q) - \tfrac{s}{2}$ and $\mhus^{(1)}_{N,0} - \mathbf{f}_0$ are nonpositive on
$\{ |v - A(q)|^2 + 2\Vext_0(q) \le s \}$ and nonnegative on its complement. Together with the bounds
\eqref{eq:ground_state_husimi_energy} and \eqref{eq:ground_state_minmax}, this gives
\begin{align}
\label{eq:ground_state_energy_chain}
    &\iint_{\mathbb{R}^6} dq \, dv \, \Big| \tfrac12 |v - A(q)|^2 + \Vext_0(q) - \tfrac{s}{2} \Big| \,
    \big| \mhus^{(1)}_{N,0}(q,v) - \mathbf{f}_0(q,v) \big|
    \nonumber \\
    &= \iint_{\mathbb{R}^6} dq \, dv \, \Big( \tfrac12 |v - A(q)|^2 + \Vext_0(q) - \tfrac{s}{2} \Big)
    \big( \mhus^{(1)}_{N,0}(q,v) - \mathbf{f}_0(q,v) \big)
    \nonumber \\
    &= \iint_{\mathbb{R}^6} dq \, dv \, \Big( \tfrac12 |v - A(q)|^2 + \Vext_0(q) \Big) \mhus^{(1)}_{N,0}(q,v)
    - \tfrac{1}{N} \big\langle \Psi_{N,0}, H_N^{(0)} \Psi_{N,0} \big\rangle
    \nonumber \\
    &\quad + \tfrac{1}{N} \big\langle \Psi_{N,0}, H_N^{(0)} \Psi_{N,0} \big\rangle
    - \iint_{\mathbb{R}^6} dq \, dv \, \Big( \tfrac12 |v - A(q)|^2 + \Vext_0(q) \Big) \mathbf{f}_0(q,v)
    \nonumber \\
    &\le \C \sqrt{\hbar} .
\end{align}
For $\delta > 0$, we split the integral at
$\big| \tfrac12 |v - A(q)|^2 + \Vext_0(q) - \tfrac{s}{2} \big| = \delta$ and use \eqref{eq:ground_state_energy_chain}:
\begin{align}
\label{eq:ground_state_L1_split}
    \iint_{\mathbb{R}^6} dq \, dv \, \big| \mhus^{(1)}_{N,0}(q,v) - \mathbf{f}_0(q,v) \big| &\le \tfrac{1}{\delta} \iint_{\mathbb{R}^6} dq \, dv \, \Big| \tfrac12 |v - A(q)|^2 + \Vext_0(q) - \tfrac{s}{2} \Big| \, \big| \mhus^{(1)}_{N,0}(q,v) - \mathbf{f}_0(q,v) \big| \nonumber \\
    &\quad + \tfrac{1}{(2\pi)^3} \, \Big| \Big\{ \Big| \tfrac12 |v - A(q)|^2 + \Vext_0(q) - \tfrac{s}{2} \Big| < \delta \Big\} \Big| \nonumber \\
    &\le \tfrac{\C \sqrt{\hbar}}{\delta} + \tfrac{1}{(2\pi)^3} \, \Big| \Big\{ \Big| \tfrac12 |v - A(q)|^2 + \Vext_0(q) - \tfrac{s}{2} \Big| < \delta \Big\} \Big| .
\end{align}

Taking first $N \to \infty$ and then $\delta \to 0$ gives
\begin{align}
\label{eq:ground_state_L1_limit}
    \lim_{N \to \infty} \iint_{\mathbb{R}^6} dq \, dv \, \big| \mhus^{(1)}_{N,0}(q,v) - \mathbf{f}_0(q,v) \big| &= 0 .
\end{align}

\emph{Assumption \ref{assump:I2} in every sector.} By \eqref{eq:ground_state_wick} and  Hadamard's inequality
\begin{align}
\label{eq:ground_state_hadamard}
    0 \le \mhus^{(k)}_{N,0}(Z_k)
    &= \tfrac{c_{N,k}}{(2\pi)^{3k}} \det \big( \langle \Psi_{N,0}, a^\dagger(\mathfrak{f}_j) \, a(\mathfrak{f}_i) \Psi_{N,0} \rangle \big)_{i,j \le k}
    \nonumber \\
    &\le \tfrac{c_{N,k}}{(2\pi)^{3k}} \prod_{i=1}^k \langle \Psi_{N,0}, a^\dagger(\mathfrak{f}_i) \, a(\mathfrak{f}_i) \Psi_{N,0} \rangle= c_{N,k} \prod_{i=1}^k \mhus^{(1)}_{N,0}(q_i, v_i) .
\end{align} 
This inequality fixes the sign of the difference, so its $L^1$ norm is its integral
\begin{align}
\label{eq:ground_state_exchange_L1}
    \Big\| \mhus^{(k)}_{N,0} - c_{N,k} \prod_{i=1}^k \mhus^{(1)}_{N,0}(q_i, v_i) \Big\|_{L^1(\mathbb{R}^{6k})}
    &= c_{N,k} - 1 .
\end{align}
Since $\mhus^{(1)}_{N,0}$ and $\mathbf{f}_0$ are probability densities and $c_{N,k} \ge 1$,
\begin{align}
\label{eq:ground_state_telescoping}
    \Big\| c_{N,k} \prod_{i=1}^k \mhus^{(1)}_{N,0}(q_i, v_i) - \mathbf{f}_0^{\otimes k} \Big\|_{L^1(\mathbb{R}^{6k})}
    &\le c_{N,k} - 1 + k \, \big\| \mhus^{(1)}_{N,0} - \mathbf{f}_0 \big\|_{L^1(\mathbb{R}^6)} ,
\end{align}
and $c_{N,k} \to 1$, \eqref{eq:ground_state_L1_limit} and \eqref{eq:ground_state_exchange_L1} give
\eqref{eq:ground_state_I2}.
\end{proof}

\begin{proposition}[Assumption \ref{assump:I4} for Slater Determinants]
\label{prop:slater_I4}
Let Assumptions \ref{assump:A1}--\ref{assump:A4} hold and let $\Psi_{N,0}$ be the Slater determinant of
\eqref{eq:ground_state_slater}. Then Assumption \ref{assump:I4} holds, with
\begin{align}
\label{eq:ground_state_I4}
    \lim_{N \to \infty} \iint_{\mathbb{R}^6} dq \, dv \, |q|^2 \, \mhus^{(1)}_{N,0}(q,v)
    &= \iint_{\mathbb{R}^6} dq \, dv \, |q|^2 \, \mathbf{f}_0(q,v) ,
\end{align}
and the energy inequality \eqref{eq:criterion_energy_noV}.
\end{proposition}

\begin{proof}[Proof of Proposition \ref{prop:slater_I4}]
By Assumption \ref{assump:A1}, $|A(q)|^2 \le \C (1 + |q|^2)$, so
$|v|^2 \le 2 |v - A(q)|^2 + \C (1 + |q|^2)$, and with $|q|^2 \le c_0^{-1} ( \Vext_0(q) + C_0 )$,
\begin{align}
\label{eq:ground_state_weight_bound}
    1 + |q|^2 + |v|^2 &\le \C \Big( 1 + \tfrac12 |v - A(q)|^2 + \Vext_0(q) + C_0 \Big) \le \C \Big( 1 + \Big| \tfrac12 |v - A(q)|^2 + \Vext_0(q) - \tfrac{s}{2} \Big| \Big),
\end{align}
where $\C$ in the last inequality depends on $C_0$ and $s$.
Multiplying \eqref{eq:ground_state_weight_bound} by $| \mhus^{(1)}_{N,0} - \mathbf{f}_0 |$ and integrating, we obtain
\begin{align}
\label{eq:ground_state_weighted_L1}
    &\iint_{\mathbb{R}^6} dq \, dv \, \big( 1 + |q|^2 + |v|^2 \big)
    \big| \mhus^{(1)}_{N,0}(q,v) - \mathbf{f}_0(q,v) \big|
\\
    &\le \C \iint_{\mathbb{R}^6} dq \, dv \, \big| \mhus^{(1)}_{N,0}(q,v) - \mathbf{f}_0(q,v) \big|
 + \C \iint_{\mathbb{R}^6} dq \, dv \, \Big| \tfrac12 |v - A(q)|^2 + \Vext_0(q) - \tfrac{s}{2} \Big| \,
    \big| \mhus^{(1)}_{N,0}(q,v) - \mathbf{f}_0(q,v) \big|,     \nonumber 
\end{align}
which tends to $0$ by \eqref{eq:ground_state_energy_chain} and \eqref{eq:ground_state_L1_limit}. Since
$\big| \iint_{\mathbb{R}^6} dq \, dv \, |q|^2 ( \mhus^{(1)}_{N,0} - \mathbf{f}_0 ) \big|$ is bounded by the left side,
\eqref{eq:ground_state_I4} follows.

By Proposition \ref{prop:slater_hypotheses}, Assumption \ref{assump:I1} holds. With
Lemma \ref{lemma:energy_decomposition_noV} at $t = 0$ and
$\tfrac12 |v|^2 + |\Vext(q)| \le \C ( 1 + |q|^2 + |v|^2 )$ by Assumption \ref{assump:A2}, we have
\begin{align}
\label{eq:ground_state_energy_limit}
    & \Big| \tfrac{1}{N} \big\langle \Psi_{N,0}, \mathcal{H} \Psi_{N,0} \big\rangle - \mathcal{E}[\mathbf{f}_0] \Big| \nonumber \\
    &\le \Big| \tfrac{1}{N} \big\langle \Psi_{N,0}, \mathcal{H} \Psi_{N,0} \big\rangle - \mathcal{V}_N(0) - \mathcal{W}_N(0) \Big|
        + \Big| \mathcal{V}_N(0) - \iint_{\mathbb{R}^6} dq \, dv \, \Big( \tfrac12 |v|^2 + \Vext(q) \Big) \mathbf{f}_0(q,v) \Big| \nonumber \\
    &\quad + \Big| \mathcal{W}_N(0) - \tfrac12 \int_{\mathbb{R}^{12}} (dq \, dv)^{\otimes 2} \, W(q_1 - q_2) \, \mathbf{f}_0(q_1,v_1) \, \mathbf{f}_0(q_2,v_2) \Big| \nonumber \\
    &\le \C \sqrt{\hbar} + \C \iint_{\mathbb{R}^6} dq \, dv \, \big( 1 + |q|^2 + |v|^2 \big) \big| \mhus^{(1)}_{N,0}(q,v) - \mathbf{f}_0(q,v) \big|
        + \tfrac12 \|W\|_{L^\infty(\mathbb{R}^3)} \big\| \mhus^{(2)}_{N,0} - \mathbf{f}_0^{\otimes 2} \big\|_{L^1(\mathbb{R}^{12})} ,
\end{align}
which tends to $0$ by \eqref{eq:ground_state_weighted_L1} and \eqref{eq:ground_state_I2} at $k = 2$. This gives
\eqref{eq:criterion_energy_noV}, with equality in the limit.
\end{proof}

\section{Convergence in Wasserstein \texorpdfstring{$\mathbb{W}_p$}{Wp} and of the Current Density}
\label{app:corollary_proofs}

This appendix proves the two corollaries of Theorem \ref{thm:magnetic_vlasov_limit}: the convergence in
$\mathbb{W}_p$, $p \in [1,2)$, from the second moment bound \eqref{eq:prop_k_moment_bound}, and the convergence of
the current density from the first velocity moment of $\mhus_{N,t}^{(1)}$ and the convergence in $\mathbb{W}_1$.

\begin{proof}[Proof of Corollary \ref{cor:wasserstein_p_convergence}]
In the proof of theorem \ref{thm:magnetic_vlasov_limit}, we have obtained the uniform bounds for second moments. By applying \cite[Thm.~6.9]{villani2008optimal}, the proof is reduced to show that the follwoing convergence for $1\leq p<2$.
\begin{align}
\label{eq:uniform_tail_p_moment}
    &\limsup_{N \to \infty} \int_{\{ \sum_{j=1}^k ( |q_j|^p + |v_j|^p ) \ge R^p \}} \dZk \,
    \sum_{j=1}^k \Big( |q_j|^p + |v_j|^p \Big) \, \mhus_{N,t}^{(k)}(Z_k) \nonumber \\
    &\le R^{p-2} \limsup_{N \to \infty} \int_{\mathbb{R}^{6k}} \dZk \,
    \Big( \sum_{j=1}^k \big( |q_j|^p + |v_j|^p \big) \Big)^{2/p} \, \mhus_{N,t}^{(k)}(Z_k) \nonumber \\
    &\le R^{p-2} \, (2k)^{2/p - 1} \limsup_{N \to \infty} \int_{\mathbb{R}^{6k}} \dZk \,
    \sum_{j=1}^k \Big( |q_j|^2 + |v_j|^2 \Big) \, \mhus_{N,t}^{(k)}(Z_k) \nonumber \\
    &\le R^{p-2} \, (2k)^{2/p - 1} \, \C k \big( 1 + t \big)^2
    \; \xrightarrow[R \to \infty]{} \; 0 .
\end{align}
\end{proof}

\begin{proof}[Proof of Corollary \ref{cor:quantitative_transport_rate_noV}]
Expanding with the definition of $\mathbf{J}_t$ in
\eqref{eq:quantum_current_def}
\begin{align}
\label{eq:current_split_noV}
    \left| \int_{\mathbb{R}^3} dq \, \phi(q) \Big( \tfrac{1}{N} \mathbf{j}_{\Psi_{N,t}}(q) - \mathbf{J}_t(q) \Big) \right| &\le \left| \tfrac{1}{N} \int_{\mathbb{R}^3} dq \, \phi(q) \, \mathbf{j}_{\Psi_{N,t}}(q) - \iint_{\mathbb{R}^6} dq \, dv \, \phi(q) \, v \, \mhus_{N,t}^{(1)}(q,v) \right| \nonumber \\
    &\quad + \left| \iint_{\mathbb{R}^6} dq \, dv \, \phi(q) \, v \, \big( \mhus_{N,t}^{(1)}(q,v) - \mathbf{f}_t(q,v) \big) \right| = I_1+I_2.
\end{align}
Notice that $\nabla_x \| a_x \Psi_{N,t} \|^2 = 2 \operatorname{Re} \langle \nabla_x a_x \Psi_{N,t}, a_x \Psi_{N,t} \rangle$ and
\eqref{eq:quantum_current_def}, we obtain
\begin{align}
\label{eq:current_imaginary_part}
    \langle \mpi_x a_x \Psi_{N,t}, a_x \Psi_{N,t} \rangle &= i\hbar \, \langle \nabla_x a_x \Psi_{N,t}, a_x \Psi_{N,t} \rangle + A(x) \, \| a_x \Psi_{N,t} \|^2 \nonumber \\
    &= \Big( - \hbar \operatorname{Im} \langle \nabla_x a_x \Psi_{N,t}, a_x \Psi_{N,t} \rangle + A(x) \, \| a_x \Psi_{N,t} \|^2 \Big) + i\hbar \operatorname{Re} \langle \nabla_x a_x \Psi_{N,t}, a_x \Psi_{N,t} \rangle \nonumber \\
    &= \operatorname{Re} \langle a_x \Psi_{N,t}, \mpi_x a_x \Psi_{N,t} \rangle + \tfrac{i\hbar}{2} \nabla_x \| a_x \Psi_{N,t} \|^2 \nonumber \\
    &= \mathbf{j}_{\Psi_{N,t}}(x) + \tfrac{i\hbar}{2} \nabla_x \| a_x \Psi_{N,t} \|^2 .
\end{align}
Therefore, it follows
\begin{align}
\label{eq:current_transfer_bound_noV}
    I_1 &= \left| \tfrac{1}{N} \int_{\mathbb{R}^3} dq \, \phi(q) \, \mathbf{j}_{\Psi_{N,t}}(q) - \iint_{\mathbb{R}^6} dq \, dv \, \phi(q) \, v \, \mhus_{N,t}^{(1)}(q,v) \right| \nonumber \\
    &= \Big| \tfrac{1}{N} \int_{\mathbb{R}^3} dq \, \phi(q) \, \mathbf{j}_{\Psi_{N,t}}(q) - \tfrac{\hbar^{-3/2}}{(2\pi)^3} \int_{\mathbb{R}^3} dq \, \phi(q) \iint_{\mathbb{R}^6} dx \, dy \, f\left(\tfrac{x-q}{\sqrt{\hbar}}\right) f\left(\tfrac{y-q}{\sqrt{\hbar}}\right) \nonumber \\
    &\qquad \times \langle \Psi_{N,t}, a_x^\dagger a_y \Psi_{N,t} \rangle \int_{\mathbb{R}^3} dv \, v \, e^{\frac{i}{\hbar} (v-A(q)) \cdot (x-y)} \Big| \nonumber \\
    &= \Big| \tfrac{1}{N} \int_{\mathbb{R}^3} dq \, \phi(q) \, \mathbf{j}_{\Psi_{N,t}}(q) - \hbar^{3/2} \int_{\mathbb{R}^3} dq \, \phi(q) \int_{\mathbb{R}^3} dx \, f\left(\tfrac{x-q}{\sqrt{\hbar}}\right) \nonumber \\
    &\qquad \times \Big\langle \big( -i\hbar \nabla_x + A(q) \big) \Big[ f\left(\tfrac{x-q}{\sqrt{\hbar}}\right) a_x \Psi_{N,t} \Big], a_x \Psi_{N,t} \Big\rangle \Big| \nonumber \\
    &= \Big| \tfrac{1}{N} \int_{\mathbb{R}^3} dq \, \phi(q) \, \mathbf{j}_{\Psi_{N,t}}(q) - \hbar^{3/2} \int_{\mathbb{R}^3} dq \, \phi(q) \int_{\mathbb{R}^3} dx \, \Big( f^2\left(\tfrac{x-q}{\sqrt{\hbar}}\right) \Big( \mathbf{j}_{\Psi_{N,t}}(x) \nonumber \\
    &\qquad + \big( A(q) - A(x) \big) \langle \Psi_{N,t}, a_x^\dagger a_x \Psi_{N,t} \rangle \Big) + \tfrac{i\hbar}{2} \nabla_x \Big( f^2\left(\tfrac{x-q}{\sqrt{\hbar}}\right) \| a_x \Psi_{N,t} \|^2 \Big) \Big) \Big| \nonumber \\
    &= \Big| \hbar^{3/2} \iint_{\mathbb{R}^6} dq \, dx \, f^2\left(\tfrac{x-q}{\sqrt{\hbar}}\right) \Big( \big( \phi(x) - \phi(q) \big) \mathbf{j}_{\Psi_{N,t}}(x)
        - \phi(q) \big( A(q) - A(x) \big) \langle \Psi_{N,t}, a_x^\dagger a_x \Psi_{N,t} \rangle \Big) \Big| \nonumber \\
    &\le \hbar^{3/2} \iint_{\mathbb{R}^6} dq \, dx \, |x - q| \, f^2\left(\tfrac{x-q}{\sqrt{\hbar}}\right) \nonumber \\
    &\qquad \times \Big( \|\nabla \phi\|_{L^\infty(\mathbb{R}^3)} \big| \mathbf{j}_{\Psi_{N,t}}(x) \big| + \|\phi\|_{L^\infty(\mathbb{R}^3)} \|\nabla A\|_{L^\infty(\mathbb{R}^3)} \langle \Psi_{N,t}, a_x^\dagger a_x \Psi_{N,t} \rangle \Big) \nonumber \\
    &\le \sqrt{\hbar} \, \big\| \sqrt{|\cdot|} f \big\|_{L^2(\mathbb{R}^3)}^2 \Big( \|\nabla \phi\|_{L^\infty(\mathbb{R}^3)} \Big( \tfrac{2}{N} \langle \Psi_{N,t}, \mathcal{K}_A \Psi_{N,t} \rangle \Big)^{1/2}
        + \|\nabla A\|_{L^\infty(\mathbb{R}^3)} \|\phi\|_{L^\infty(\mathbb{R}^3)} \Big) \nonumber \\
    &\le \sqrt{\hbar} \, \big\| \sqrt{|\cdot|} f \big\|_{L^2(\mathbb{R}^3)}^2 \Big( \C + \|\nabla A\|_{L^\infty(\mathbb{R}^3)} \Big) \|\phi\|_{W^{1,\infty}(\mathbb{R}^3)} ,
\end{align}
where in the fourth equation we change variables $q \to x$ to get another $\mathbf{j}_{\Psi_{N,t}}(x)$ and insert $1 = \hbar^{-3/2} \int_{\mathbb{R}^3} dq \, f^2\big( \frac{x-q}{\sqrt{\hbar}} \big)$,
in the last line we use the kinetic energy bound \eqref{eq:kinetic_energy_uniform}.

For the second term in \eqref{eq:current_split_noV}, fix an optimal transport plan $\gamma_t^{(1)} \in \Gamma\big(\mhus_{N,t}^{(1)}, \mathbf{f}_t\big)$ for $\mathbb{W}_1$
\cite[Def.~6.1, Thm.~4.1]{villani2008optimal}. Since
$\iint_{\mathbb{R}^6} dq \, dv \, |v| \, \mhus_{N,t}^{(1)}(q,v) < \infty$ by \eqref{eq:prop_k_moment_bound} and $|\tilde v| \le R(t)$
on $\operatorname{supp}\mathbf{f}_t$ by \eqref{eq:vlasov_support}, both integrals are well-defined
against $\gamma_t^{(1)}$. With
$\phi(q)v - \phi(\tilde q)\tilde v = \phi(q)(v - \tilde v) + \big(\phi(q) - \phi(\tilde q)\big)\tilde v$, we have
\begin{align}
\label{eq:current_primal_expansion_noV}
    I_2 &= \left| \iint_{\mathbb{R}^6} dq \, dv \, \phi(q) \, v \, \big( \mhus_{N,t}^{(1)}(q,v) - \mathbf{f}_t(q,v) \big) \right| \nonumber \\
    &= \left| \iint_{\mathbb{R}^{12}} d\gamma_t^{(1)}(q,v,\tilde{q},\tilde{v}) \, \phi(q) v
    - \iint_{\mathbb{R}^{12}} d\gamma_t^{(1)}(q,v,\tilde{q},\tilde{v}) \, \phi(\tilde{q}) \tilde{v} \right| \nonumber \\
    &= \left| \iint_{\mathbb{R}^{12}} d\gamma_t^{(1)}(q,v,\tilde{q},\tilde{v}) \,
    \Big( \phi(q)v - \phi(\tilde{q})\tilde{v} \Big) \right| \nonumber \\
    &\le \iint_{\mathbb{R}^{12}} d\gamma_t^{(1)}(q,v,\tilde{q},\tilde{v}) \,
    \Big( \|\phi\|_{L^\infty(\mathbb{R}^3)} |v - \tilde{v}|
        + R(t) \|\nabla \phi\|_{L^\infty(\mathbb{R}^3)} |q - \tilde{q}| \Big) \nonumber \\
    &\le \Big( \|\phi\|_{L^\infty(\mathbb{R}^3)} + R(t) \|\nabla \phi\|_{L^\infty(\mathbb{R}^3)} \Big)
    \iint_{\mathbb{R}^{12}} d\gamma_t^{(1)}(q,v,\tilde{q},\tilde{v}) \,
    \Big( |q - \tilde{q}| + |v - \tilde{v}| \Big) \nonumber \\
    &\le \sqrt{2} \Big( \|\phi\|_{L^\infty(\mathbb{R}^3)} + R(t) \|\nabla \phi\|_{L^\infty(\mathbb{R}^3)} \Big)
    \mathbb{W}_1\big(\mhus_{N,t}^{(1)}, \mathbf{f}_t\big) ,
\end{align}
where we used $|q - \tilde{q}| + |v - \tilde{v}| \le \sqrt{2} \, |(q,v) - (\tilde{q},\tilde{v})|$ and the optimality of $\gamma_t^{(1)}$ in the last line. With \eqref{eq:current_transfer_bound_noV} and $1 + R(t) \le \C (1+t)^2$ by \eqref{eq:vlasov_support},
\eqref{eq:current_split_noV} gives \eqref{eq:current_convergence_bound_cor_noV},
and Corollary \ref{cor:wasserstein_p_convergence} at $p = 1$ gives the limit.
\end{proof}

\section*{Acknowledgements}
The research of L. Chen is partly supported by the German Research Foundation (No. 547277619). 

\section*{AI Disclosure}
The authors used artificial intelligence models for 
language formulation, typesetting and formatting, standard algebraic manipulations and derivations, consistency checks, for diagrams, and the search for literature and references. The authors reviewed and verified all AI-assisted content, including the code of the figures and the table, and take full responsibility for it as authors.

\section*{Competing Interests}
The authors have no competing interests to declare that are relevant to the content of this article.

\bibliographystyle{plain}
\bibliography{references}

\end{document}